\documentclass[a4paper, 12pt]{amsart}
\usepackage[all,cmtip]{xy}
\usepackage{amssymb}
\usepackage{amsfonts}
\usepackage{amsthm}
\usepackage{amsmath}
\usepackage{amscd}
\usepackage{color} 
\usepackage[utf8]{inputenc}
\usepackage{hyperref}
\usepackage{fullpage}
\usepackage{epic}
\usepackage{mathtools}
\usepackage{tikz}
\usetikzlibrary{chains}

\tikzset{node distance=2em, ch/.style={circle,draw,on chain,inner sep=2pt},chj/.style={ch,join},every path/.style={shorten >=4pt,shorten <=4pt},line width=1pt,baseline=-1ex}

\usepackage{graphicx}
\usepackage[all]{xy}
\usepackage{array}
\usepackage{youngtab}

\newcommand{\li}{{\,|\,}}

\newcommand{\fb}{{\mathfrak{b}}}

\newcommand{\fg}{{\mathfrak{g}}}
\newcommand{\fh}{{\mathfrak {h}}}

\newcommand{\fn}{{\mathfrak {n}}}

\newcommand{\fp}{{\mathfrak {p}}}

\newcommand{\ft}{{\mathfrak{t}}}

\newcommand{\bbC}{{\mathbb {C}}}

\newcommand{\bbN}{{\mathbb {N}}}

\newcommand{\bbP}{{\mathbb {P}}}
\newcommand{\bbQ}{{\mathbb {Q}}}
\newcommand{\bbR}{{\mathbb {R}}}

\newcommand{\bbZ}{{\mathbb {Z}}}

\newcommand{\CA}{{\mathcal {A}}}
\newcommand{\CB}{{\mathcal {B}}}
\newcommand{\CC}{{\mathcal {C}}}

\newcommand{\CF}{{\mathcal {F}}}

\newcommand{\CH}{{\mathcal {H}}}
\newcommand{\CI}{{\mathcal {I}}}

\newcommand{\CL}{{\mathcal {L}}}

\newcommand{\CP}{{\mathcal {P}}}

\newcommand{\CV}{{\mathcal {V}}}

\newcommand{\tr}{{\mathrm {tr}}}

\newcommand{\Ch}{{\mathrm{ch}}}

\newcommand{\Hom}{{\mathrm{Hom}}}

\newcommand{\Gr}{{\mathrm{Gr}}}
\newcommand{\Fl}{{\mathrm{Fl}}}

\newcommand{\SL}{{\mathrm{SL}}}

\newcommand{\ev}{{\mathrm{ev}}}

\newcommand{\incl}{\hookrightarrow}

\newcommand{\C}{\mathbb{C}}

\newcommand{\la}{\langle}
\newcommand{\ra}{\rangle}

\newcommand{\Wmin}{W^{\dagger}_{\ell+\check{h}}}
\newcommand{\Wl}{W_{\ell+\check{h}}}

\newcommand{\be}{\begin{equation}}
\newcommand{\ee}{\end{equation}}
\newcommand{\bt}{\begin{theorem}}
\newcommand{\et}{\end{theorem}}
\newcommand{\bd}{\begin{definition}}
\newcommand{\ed}{\end{definition}}
\newcommand{\bp}{\begin{proposition}}
\newcommand{\ep}{\end{proposition}}
\newcommand{\bl}{\begin{lemma}}
\newcommand{\el}{\end{lemma}}
\newcommand{\bco}{\begin{corollary}}
\newcommand{\eco}{\end{corollary}}
\newcommand{\br}{\begin{remark}}
\newcommand{\er}{\end{remark}}
\newcommand{\bex}{\begin{example}}
\newcommand{\eex}{\end{example}}
\newcommand{\ben}{\begin{enumerate}}
\newcommand{\een}{\end{enumerate}}
\newcommand{\bc}{\begin{cases}}
\newcommand{\ec}{\end{cases}}

\newcommand{\bpf}{\begin{proof}}
\newcommand{\epf}{\end{proof}}

\newcommand{\bma}{\begin{bmatrix}}
\newcommand{\ema}{\end{bmatrix}}

\newcommand{\lwedge}{\mbox{\Large $\wedge$}}

\theoremstyle{Theorem}

\theoremstyle{Theorem}

\theoremstyle{Theorem}

\theoremstyle{Definition}

\newtheorem{theorem}{Theorem}[section]

\newtheorem{lemma}[theorem]{Lemma}
\newtheorem{proposition}[theorem]{Proposition}
\newtheorem{corollary}[theorem]{Corollary}

\newtheorem{remark}[theorem]{Remark}

\theoremstyle{definition}
\newtheorem{definition}[theorem]{Definition}

\newtheorem{example}[theorem]{Example}

\begin{document}

\title{ Conformal blocks, Verlinde formula and diagram automorphisms}

\author{Jiuzu Hong}
\address{Department of Mathematics, University of North Carolina at Chapel Hill, Chapel Hill, NC 27599-3250, U.S.A.}
\email{jiuzu@email.unc.edu}
\maketitle
\begin{abstract}
The Verlinde formula computes the dimension of the space of conformal blocks associated to simple Lie algebras and stable pointed curves.
 If a simply-laced simple Lie algebra admits a nontrivial diagram automorphism, then this automorphism acts on the space of conformal blocks naturally.
We prove an analogue of the Verlinde formula for the trace of the diagram automorphism on the space of conformal blocks. Along the way, we get an analogue of the Kac-Walton formula for the trace of the diagram automorphism. We also get a twining type formula between the conformal blocks for the pair $(sl_{2n+1},sp_{2n})$.
\end{abstract}

\tableofcontents

\section{Introduction}

The Verlinde formula computes the dimension of the space of conformal blocks. It is fundamentally important in conformal field theory and algebraic geometry. The formula was originally conjectured by Verlinde \cite{V}  in conformal field theory. It was mathematically  derived by combining the efforts of  mathematicians including  Tsuchiya-Ueno-Yamada \cite{TUY}, Faltings\cite{Fa}. It was  proved by  Beauville-Laszlo\cite{BL}, Kumar-Narasimhan-Ramanathan \cite{KNR}, Faltings\cite{Fa}, that conformal blocks can be identified with the generalized theta functions on the moduli stack of $G$-bundles on projective curves where $G$ is a simply-connected simple algebraic group. Therefore the Verlinde formula also computes the dimension of the spaces of generalized theta functions. For a survey on Verlinde formula, see Sorger's Bourbaki talk \cite{So}.

Let $(C,\vec{p})$ be a stable $k$-pointed curve. Let $\fg$ be a  simple Lie algebra over $\bbC$. Let $\ell$ be a positive integer. Put
\begin{equation}
\label{weight_level}
P_\ell=\{ \lambda\in  P^+ \li   \la \lambda, \check{\theta} \ra\leq   \ell   \ra      \},
\end{equation}
where $\theta$ is the highest root of $\fg$ and $\check{\theta}$ is the coroot of $\theta$. Given a tuple of dominant weights $\vec{\lambda}=(\lambda_1,\lambda_2,\cdots, \lambda_k)$ such that $\lambda_i\in P_\ell$ for each $i$. We can attach the space $V_{\fg, \ell, \vec{\lambda}}(C,\vec{p})$ of  conformal blocks of level $\ell$  to $(C,\vec{p})$ and $\vec{\lambda}$. We will recall the definition of conformal blocks  in Section \ref{conformal_block}.

Let $\sigma$ be a  diagram automorphism on a  simple Lie algebra  $\fg$. One can attach another  simple Lie algebra $\fg_\sigma$ as the orbit Lie algebra of $\fg$ (see Section \ref{pair_groups} for  details). If $\sigma$ is trivial, then $\fg=\fg_\sigma$. Let $\Phi$ (resp.\,$\Phi_\sigma$ ) be the set of roots of $G$ (resp.\,$G_\sigma$). We put
\be
\label{Delta}
  \Delta= \prod_{\alpha\in \Phi}  (e^{\alpha}-1), \quad    \Delta_\sigma= \prod_{\alpha\in \Phi_\sigma}  (e^{\alpha}-1). 
\ee

 There is a natural correspondence between $\sigma$-invariant weights (resp.\,dominant weights) of $\fg$ and weights (resp.\,dominant weights) of $\fg_\sigma$ (see Section \ref{root_system}). In this introduction we will identify them if no confusion occurs. 
  For any dominant weight $\lambda$ of $\fg$ (resp.\,$\fg_\sigma$), we denote by $V_\lambda$ (resp.\,$W_\lambda$) the irreducible representation of $\fg$ (resp.\,$\fg_\sigma$) of highest weight $\lambda$. Let $\check{h}$ (resp.\,$\check{h}_\sigma$) be the dual Coxeter number of $\fg$ (resp.\,$\fg_\sigma$).
 Let $G$ (resp.\,$G_\sigma$) be the associated simply-connected simple algebraic group of $\fg$ (resp.\,$\fg_\sigma$). Let $T$ (resp.\,$T_\sigma$) be a maximal torus of $G$ (resp.\,$G_\sigma$). Let $W$ (resp.\,$W_\sigma$) denote the Weyl group of $G$ (resp.\,$G_\sigma$). 

Throughout this paper, we denote by $\tr(A|V)$ the trace of an operator $A$ on a finite dimensional vector space $V$. 
The following is the celebrated Verlinde formula. 
\begin{theorem}[Verlinde formula]
\label{verlinde_formula}
Let $(C,\vec{p})$ be a stable $k$-pointed curve of genus $g$.
Given any tuple $\vec{\lambda}=(\lambda_1,\lambda_2,\cdots, \lambda_k)$ of dominant weights  of $\fg$ such that $\lambda_i\in P_\ell$ for each $i$, we have 
\be
\label{verlinde}
 \dim V_{\fg, \ell, \vec{\lambda}}(C,\vec{p})=|T_\ell|^{g-1}    \sum_{t\in  T_\ell^{\rm reg}/W }   \frac{  \tr (t|V_{\vec{\lambda}})  }{ \Delta(t)^{g-1}   }, 
 \ee
where $V_{\vec{\lambda}}$ denotes the tensor product  $V_{\lambda_1}\otimes \cdots \otimes V_{\lambda_k}$ of representations of $\fg$  and
\[T_\ell=\{t\in  T \li e^{\alpha}(t)=1, \, \alpha\in (\ell+\check{h})Q_l   \}  \]
 is a finite abelian subgroup in the maximal torus $T$, $T^{\rm reg}_{\ell}$ denotes the set of regular elements in $T_\ell$ and $T_\ell^{\rm reg}/W$ denotes the set of $W$-orbits in $T^{\rm reg}_{\ell}$. Here 
   $Q_l$ denotes the lattice spanned by long roots of $\fg$, and for any $\alpha\in Q_l$, $e^{\alpha}$ is the associated character of $T$.
 \end{theorem}

From now on we always assume $\sigma$ is  nontrivial.
When the tuple $\vec{\lambda}$ of dominant weights  of $\fg$ is $\sigma$-invariant, one can define a natural operator on the space $V_{\fg, \ell, \vec{\lambda}}(C,\vec{p})$ of conformal blocks, which we still denote by $\sigma$, see Section \ref{conformal_block}. A natural question is how to compute the trace of $\sigma$ as an operator on the space of the conformal blocks. In this paper, we derive a formula  for the trace of $\sigma$, which is very similar to the Verlinde formula for the dimension of the space of conformal blocks. Very surprisingly, in the formula the role of $\fg$ is replaced by  $\fg_\sigma$.
 The following is the main theorem of this paper. 
\bt  
\label{main}
Let $(C,\vec{p})$ be a stable $k$-pointed curve of genus $g$. Let $\sigma$ be a nontrivial diagram automorphism on a simple Lie algebra $\fg$ which has dual Coxeter number $\check{h}$.
Given a tuple $\vec{\lambda}=(\lambda_1,\lambda_2,\cdots, \lambda_k)$  of $\sigma$-invariant dominant weights of $\fg$  such that for each $i$, $\lambda_i\in P_\ell$, 
 we have the following formula
\be
\label{main_formula}
\tr(\sigma| V_{\fg, \ell, \vec{\lambda}}(C,\vec{p}) )=|T_{\sigma, \ell}|^{g-1} \sum_{t\in T^{\rm reg}_{\sigma, \ell}/W_\sigma  }  \frac{\tr(t| W_{\vec{\lambda}} )   }{ \Delta_\sigma(t) ^{g-1}},
\ee
where $W_{\vec{\lambda}}$ denotes the tensor product  $W_{\vec{\lambda}}:=W_{\lambda_1}\otimes \cdots  \otimes W_{\lambda_k}$  of  representations of $\fg_\sigma$ and  
\[ T_{\sigma,\ell}= \{ t\in T_\sigma  \li   e^{\alpha}(t)=1, \alpha\in  (\ell+\check{h}) Q^\sigma   \}. \]
Here $T^{\rm reg}_{\sigma,\ell}$ denotes the set of regular elements in $T_{\sigma,\ell}$, and  $T^{\rm reg}_{\sigma, \ell}/W_\sigma$ denotes the set of $W$-orbits in   $T^{\rm reg}_{\sigma,\ell}$ and
\[ Q^\sigma=\begin{cases}   
\text{ root lattice of } \fg_\sigma  \quad  \text{ if } \fg\neq A_{2n}  \\
\text{ weight lattice of } \fg_\sigma   \quad   \text{ if }  \fg=A_{2n}.
\end{cases}  \]
\et
Since the space of conformal blocks can be identified with the space of generalized theta functions, Theorem \ref{main} implies the same formula for the trace of diagram automorphisms on the space of generalized theta functions.

\br
By the basic representation theory of finite groups, we  have the following formula   
\[\dim  V_{\fg, \ell, \vec{\lambda}}(C,\vec{p})^\sigma =\frac{1}{r}  \sum_{i=1}^r  \tr(\sigma^i |  V_{\fg, \ell, \vec{\lambda}}(C,\vec{p})),\]
where $r$ is the order of $\sigma$, and $ V_{\fg, \ell, \vec{\lambda}}(C,\vec{p})^\sigma$ denotes the space of $\sigma$-invariants in  $V_{\fg, \ell, \vec{\lambda}}(C,\vec{p})$. Combining  Theorem \ref{verlinde_formula} and Theorem \ref{main}, we immediately get a formula for the dimension of  $V_{\fg, \ell, \vec{\lambda}}(C,\vec{p})^\sigma$. 
\er

The  proof of  Theorem \ref{main} will be completed in Section \ref{proof_section}. Our proof closely follows  \cite{Fa,Be, Ku2} for the  derivation of the usual Verlinde formula, where the fusion ring plays essential role. 
In the standard approach to the Verlinde formula for general stable pointed curves, the factorization rules for conformal blocks and degeneration of projective smooth curves allow a reduction to projective line with three points case. Our basic idea is that  we  replace the dimension by the trace of the diagram automorphism everywhere. In our taking trace setting, we explain in Section \ref{conformal_block} that factorization rules for conformal blocks and degeneration of curves are compatible well with the trace operation on the space of conformal blocks. 

By replacing the dimension by the trace, we introduce $\sigma$-twisted fusion rings $R_\ell(\fg,\sigma)$ in Section \ref{twisted_fusion_sect}. 
We also introduce the $\sigma$-twisted representation ring $R(\fg,\sigma)$ of $\fg$ (see Section \ref{twist_ring_sect}). For the usual fusion ring $R_\ell(\fg)$ and the representation ring $R(\fg)$, it is important to establish a ring homomorphism from $R(\fg)$ to $R_\ell(\fg)$. Similarly, we establish a ring homomorphism from $R(\fg,\sigma)$ to $R_\ell(\fg,\sigma)$ in Section \ref{ring_hom_sect}. One of important technical tools is that we interpret  $\sigma$-twisted fusion product via affine analogue of Borel-Weil-Bott theorem, where  the new product is introduced in Section \ref{new_def_fusion_sect}. 
   A vanishing theorem of Lie algebra cohomology by Teleman \cite{Te} plays a key role in our arguments as in the dimension setting (cf.\,\cite[Chapter 4]{Ku2} ). We describe all characters of the ring $R_\ell(\fg, \sigma)$ in Section \ref{character_fusion_sect}. The Verlinde formula for the trace of diagram automorphism will be a consequence of the characterization of the ring $R_\ell(\fg,\sigma)$ and the determination of the Casimir element in $R_\ell(\fg,\sigma)$. As a byproduct we obtain an analogue of Kac-Walton formula (Theorem \ref{trace_conformal_tensor}) in Section \ref{new_def_fusion_sect}.

In the process of proving the coincidence of two products in the ring $R_\ell(\fg, \sigma)$  and establishing the ring homomorphism from $R(\fg,\sigma)$ to $R_\ell(\fg,\sigma)$, some interesting sign problems occur on the higher cohomology groups of vector bundles on affine Grassmannian and affine flag variety, also in affine BBG-resolution and affine Kostant homologies. The resolution of these sign problems is very crucial for the characterization of the ring $R_\ell(\fg,\sigma)$.


 Let   $\CL_\ell(V_\lambda)$ be the vector bundle on the affine Grassmannian $\Gr_G$ of $G$ associated to the level $\ell$ and the representation $V_\lambda$ of $G$. By affine Borel-Weil-Bott theorem (cf.\,\cite{Ku1}) there is  only one nonzero cohomology $H^{i}(\Gr_G, \CL_\ell(V_\lambda))$  and the restricted dual $H^i(\Gr_G, \CL_\ell(V_\lambda))^\vee$ is the irreducible integrable representation $\CH_\lambda$ of the affine Lie algebra $\hat{\fg}$ of level $\ell$. The action of $\sigma$ on the highest weight vectors of $H^i(\Gr_G, \CL_\ell(V_\lambda))^\vee$ is determined in Section \ref{sign_BWB_sect} and Section \ref{sign_affine}. This problem is closely related to similar problem on  the cohomology of line bundle on affine flag variety. The answer is very similar to the finite-dimenisonal situation which is due to Naito \cite{N1} where  Lefschetz fixed point formula is used. In the affine setting, we don't know how to apply Lefschetz fixed point formula since the affine Grassmannian and affine flag variety are infinite-dimensional. Instead, our method is inspired by Lurie's short proof of Borel-Weil-Bott theorem  \cite{Lu}. Our method should be applicable to  similar problems for general symmetrizable Kac-Moody groups. Similar sign problems also appear in BGG resolution and the Kostant homology for affine Lie algebras. See the discussions in Section \ref{affine_BBG_Kostant}.

Our starting point of this work is the Jantzen's twining  formula (cf.\,\cite{Ja,Ho,FSS,KLP,N1,N2})  relating representations of $\fg$ and $\fg_\sigma$, where the term ``twining" is coined by Fuchs-Schellekens-Schweigert  \cite{FSS}. Given a $\sigma$-invariant dominant weight $\lambda$ of $\fg$ where $\sigma$ is the diagram automorphism as above. There is a unique operator $\sigma$ on $V_\lambda$ such that $\sigma$  preserves the highest weight vector $v_\lambda\in V_\lambda$ and for any $u\in \fg$ and $v\in V_\lambda$, $\sigma(u\cdot v)=\sigma(u)\cdot \sigma(v)$. For any $\sigma$-invariant weight $\mu$, Jantzen \cite{Ja} proved the following formula
\[ \tr(\sigma|  V_\lambda(\mu))=\dim W_\lambda(\mu),\]
where $\lambda$ and $\mu$ are also regarded as (dominant) weights of $\fg_\sigma$.
Given a tuple $\vec{\lambda}$ of $\sigma$-invariant dominant weights of $\fg$. Let $V_{\vec{\lambda}}^\fg$ (resp.\,$W_{\vec{\lambda}}^{\fg_\sigma}$) be the tensor invariant space of $\fg$ (resp.\,$\fg_\sigma$). Induced from the action of $\sigma$ on each $V_{\lambda_i}$, $\sigma$ acts on $V_{\vec{\lambda}}^\fg$ diagonally. Shen and the author obtained the following twining formula in the setting of tensor invariant spaces in \cite{HS},
\begin{equation}
\label{HS_0}
 \tr(\sigma|  V_{\vec{\lambda}}^\fg  )= \dim W_{\vec{\lambda}}^{\fg_\sigma}.
 \end{equation}
 A consequence of (\ref{HS_0}) is that the $\sigma$-twisted representation ring $R(\fg,\sigma)$ of $\fg$ is isomorphic to the representation ring $R(\fg_\sigma)$ of $\fg_\sigma$. This is how we are able to express the trace of $\sigma$ on the space of conformal blocks by the data associated to $\fg_\sigma$.

It is well-known that given a tuple $\vec{\lambda}$ of dominant weights of $\fg$, the space  $V_{\fg, \ell, \vec{\lambda}}(\bbP^1, \vec{p})$  of conformal blocks on $(\bbP^1, \vec{p})$ stabilizes to the tensor co-invariant space $(V_{\vec{\lambda}})_{\fg}$ when the level $\ell$ increases. From Formula (\ref{HS_0}), it is natural to hope that the conformal blocks associated to $\fg$ and the conformal blocks associated to $\fg_\sigma$ are related and have a twining formula with a fixed level. Unfortunately this is not the case. We found the following counter-example using \cite{Sw}  (joint with P.\,Belkale). 
\bex
  We have 
\[ \dim V_{sl_6, 4, \lambda, \mu, \nu }(\bbP^1, 0,1,\infty)=4,\]
 where $\lambda=\omega_2+\omega_3+\omega_4$, $\mu=\omega_1+\omega_3+\omega_5$ and $\nu=\omega_1+2\omega_3+\omega_5$. Here $\omega_1,\omega_2,\omega_3,\omega_4,\omega_5$  denote  the fundamental weights of $sl_6$. 
 Since the order of $\sigma$ on $sl_{6}$ is 2, this forces that the trace $\tr(\sigma| V_{sl_6, 4, \lambda, \mu, \nu }(\bbP^1, 0,1,\infty))$ is even. 
 On the other hand, we have 
   \[ \dim V_{so_{7}, 4, \lambda, \mu,\nu }(\bbP^1, 0,1,\infty)=3,\] 
where $\lambda= \omega_{\sigma,2}+\omega_{\sigma,3} $, $\mu= \omega_{\sigma,1}+\omega_{\sigma,3}$ and  $\nu=\omega_{\sigma,1}+2\omega_{\sigma,3}$. Here $\omega_{\sigma,1}, \omega_{\sigma,2}, \omega_{\sigma,3}$ denotes the fundamental weights of $so_{7}$. 
\eex

Actually from formula (\ref{main_formula}) in Theorem \ref{main}, it is quite clear that  $\tr(\sigma |V_{\fg, \ell, \vec{\lambda}}(C, \vec{p}))$  should  not be the same as  $\dim V_{\fg_\sigma, \ell, \vec{\lambda}}(C, \vec{p})$. Nevertheless, for the special pair $(sl_{2n+1}, sp_{2n})$ we do have a twining formula where one needs to take different levels on both sides.
\bt
\label{trace_dim}
If $\ell$ is an odd positive integer, then the following formula holds
\begin{equation}
\label{trace_dim_formula}
  \tr(\sigma| V_{sl_{2n+1}, \ell, \vec{\lambda}}(C,\vec{p}))=\dim V_{sp_{2n}, \frac{\ell-1}{2}, \vec{\lambda}}(C,\vec{p}). 
\end{equation}
\et
This theorem is a corollary of Theorem \ref{main}, and  the proof will be given in  Section \ref{stange_corollary_sect}.
It has following interesting numerical consequences where $\ell$ is assumed to be odd. 
\begin{itemize}    
\item The trace $\tr(\sigma| V_{sl_{2n+1}, \ell, \vec{\lambda}}(C,\vec{p}))$ is non-negative. 
\item If  $\dim V_{sl_{2n+1}, \ell, \vec{\lambda}}(C,\vec{p})$ is $1$, then   $\dim V_{sp_{2n}, \frac{\ell-1}{2}, \vec{\lambda}}(C,\vec{p})$ is also $1$. 
\item   If $ V_{sp_{2n}, \frac{\ell-1}{2}, \vec{\lambda}}(C,\vec{p})$ is nonempty, then $V_{sl_{2n+1}, \ell, \vec{\lambda}}(C,\vec{p})$ is nonempty. 
\end{itemize}
Theorem \ref{trace_dim} establishes a bridge between the conformal blocks for $sl_{2n+1}$ and $sp_{2n}$. 

The failure of the formula (\ref{trace_dim_formula}) in general is not really the end of the story. The combinatoric data  appearing in the  formula in Theorem \ref{main} actually suggests a close connection with  twisted affine Lie algebras. It is very natural from the point of view of the twining formula for affine Lie algebras by Fuchs-Schellekens-Schweigert \cite{FSS}. Moreover the $\sigma$-twisted fusion ring $R_\ell(\fg,\sigma)$ defined in this paper is  closely related to Kac-Peterson formula  for S-matrices of twisted affine Lie algebras (cf.\,\cite{Ka}). The analogue of Kac-Walton formula  obtained in this paper is also a strong hint. In fact this perspective has recently been clarified in  \cite{Ho2} by the author. 
   The connection on the trace of diagram automorphism on the space of conformal blocks and certain  conformal field theory related to twisted affine Lie algebra was predicted by Fuchs-Schweigert \cite{FS}. It seems to the author that this work also has confirmed Conjecture 2 in \cite{FS2} when the automorphism is induced from a diagram automorphism of $\fg$.  This work should also be closely related to the fusion rules for the orbifold conformal field theory that is developed by Birke-Fuchs-Schweigert \cite{BFS} and Ishikawa-Tani \cite{IT}. 
   
    A general theory of twisted conformal blocks has been developed recently by S.\,Kumar and the author \cite{HK}. It would be interesting to investigate the connection between this paper and \cite{HK}.
\vspace{0.1 in}
  
\noindent{\bf Acknowledgments}
 The author would like to express his  gratitude to P.\,Belkale for introducing him into the theory of conformal blocks, and for many helpful and stimulating discussions throughout this work. He  would like to thank S.\,Kumar for  helpful discussions and for his careful reading on the first draft of the paper, and also for sharing his unpublished book on Verlinde formula \cite{Ku2}. He also  wants to thank I.\,Cherednik for his interest and many helpful comments. This work was partially supported by the Simons Foundation collaboration grant 524406.

Finally, the author would like to thank the anonymous referee for very careful reading and the help in improving the exposition of the paper. 
\section{The root systems and affine Weyl group of orbit Lie algebras}
\label{pair_groups}
\subsection{Root systems}
\label{root_system}
Let $\fg$ be a  simple Lie algebra over $\bbC$. Let $I$ be the set of vertices of the Dynkin diagram of $\fg$. 
  For each $i\in I$, let $\alpha_i$ (resp.\,$\omega_i$) be the corresponding simple root ( resp.\,fundamental weight). Let $P$ be the weight lattice of $\fg$ and let $P^+$ be the set of dominant weights of $\fg$. Let $\Phi$ (resp.\,$\check{\Phi}$) be the set of roots (resp.\,coroots) of $\fg$, and let $Q$ (resp.\,$\check{Q}$) be the root lattice (resp.\,coroot lattice) of $\fg$.
  For each root $\alpha\in \Phi$, let $\check{\alpha}\in \check{\Phi}$ be the associated coroot of $\alpha$. Let $\la, \ra : P\times  \check{Q}\to \bbZ$ be the perfect pairing between weight lattice and coroot lattices. Note that  the matrix  $(\la \alpha_i, \check{\alpha}_j \ra)$ is the Cartan matrix of $\fg$.

We denote by $e_i,f_i,h_i$ the corresponding Chevalley generators in $\fg$ for each $i\in I$. 
Let $\sigma$ be a nontrivial diagram automorphism of the Dynkin diagram of $\fg$. Note that $\fg$ can only be of types $A_n, D_n, E_6$ when $\sigma$ is nontrivial. 
  The automorphism $\sigma$ acts on $P$, such that $\sigma(\alpha_i)=\alpha_{\sigma(i)}$ and $\sigma(\omega_i)=\omega_{\sigma(i)}$ for each $i\in I$. Clearly $\sigma$ maps each dominant weight to another dominant weight. 

The diagram automorphism $\sigma$ defines an automorphism $\sigma$ of the Lie algebra $\fg$  such that $\sigma(e_i)=e_{\sigma(i)}$, $\sigma(f_i)=f_{\sigma(i)}$, $\sigma(h_i)=h_{\sigma(i)}$ for each $i\in I$. Here  we use the same notation $\sigma$ to denote these automorphisms if no confusion occurs. 

Let $I_\sigma$ be the set of orbits of $\sigma$ on $I$. There exists a unique  simple Lie algebra  $\fg_\sigma$ over $\bbC$ whose vertices of Dynkin diagram is indexed by $I_\sigma$ (cf.\,\cite[Section 2.2]{HS}), and the Cartan matrix is given as follows, 
\[a_{\imath \jmath}= \begin{cases}
\frac{ |\imath|}{2}  a_{ij}, \quad  \fg \text{ is of type } A_{2n} \text{ and }  \imath \text{ is disconnected}  \\
 |\imath|  a_{ij}, \quad \text{ otherwise}
 \end{cases}
  \]
    for any $ \imath\not= \jmath\in I_\sigma$, 
where $i\in \imath, j\in \jmath$ and $|\imath|$ is the cardinality of the $\imath$. The Lie algebra $\fg_\sigma$ is called the orbit Lie algebra of $(\fg,\sigma)$ in literature.

Let $\alpha_\imath$ (resp.\,$\check{\alpha}_\imath$) be the corresponding simple root (resp.\,simple coroot) for $\imath\in I_\sigma$. Let $P_{\sigma}$ be the weight lattice of $\fg_\sigma$. There exists a bijection of lattices $\iota: P^{\sigma}\simeq P_\sigma$ such that $\iota^{-1}(\omega_\imath)=\sum_{i\in \imath} \omega_i$ for each $\imath\in I_\sigma$, where $P^\sigma$ is the fixed point lattice of $\sigma$ on $P$. Let $\rho$ (resp.\,$\rho_\sigma$) be the sum of all fundamental weights of $\fg$ (resp.\,$\fg_\sigma$). Note that $\rho\in P^\sigma$, and  $\iota(\rho)=\rho_\sigma$. Moreover, 
 \be
 \label{root}
  \iota^{-1}(  \alpha_\imath)=\begin{cases} 
 \sum_{i\in \imath}  \alpha_i    \quad \text{ for any } i\not=j\in  \imath, a_{ij}=0\\
2(\alpha_i+\alpha_j )   \quad   \imath=\{i,j\}, a_{ij}=-1
 \end{cases}.
   \ee
  Let $Q_\sigma$ (resp.\,$\check{Q}_\sigma$) be the root lattice (resp.\,coroot lattice) of $\fg_\sigma$. 
There is  a projection map $\check{\iota}: \check{Q}\to  \check{Q}_{\sigma}$. 
 Under this projection, we have 
 \[ \check{\iota}(\check{\alpha}_i) =  \check{\alpha}_\imath, \quad \text{ for any }  i\in \imath. \] 
  For any $\lambda\in P_\sigma$ and $x\in \check{Q}_\sigma$, we have the following compatibility
  \be
  \label{compatibility}
    \la \iota(\lambda), \check{\iota} (x)  \ra=\la \lambda, x  \ra_{\sigma}, \ee
 where $\la, \ra_{\sigma}: P_{\sigma}\times \check{Q}_{\sigma}\to \bbZ$ is the perfect pairing between the weight lattice and dual root lattice for $\fg_\sigma$.
  The following is a table of  $\fg$ and $\fg_\sigma$ for nontrivial $\sigma$ (\cite[Section 2.2]{HS} or \cite[6.4]{Lus}):
\begin{enumerate}
\item If $\fg=A_{2n-1}$ and $\sigma$ is of order $2$, then $\fg_\sigma=B_n$, $n\geq 2$.
\item If $\fg=A_{2n}$ and $\sigma$ is of order $2$, then $\fg_\sigma=C_n$, $n\geq 1$, where $C_1$ by convention means $A_1$.
\item If $\fg=D_{n}$ and $\sigma$ is of order $2$, then $\fg_\sigma=C_{n-1}$, $n\geq 4$.
\item If $\fg=D_4$ and $\sigma$ is of order $3$, then $\fg_\sigma=G_2$.
\item If $\fg=E_6$ and $\sigma$ is of order $2$, then $\fg_\sigma=F_4$.
\end{enumerate}

 Let $\theta$ be the highest root of $\fg$. It is clear that $\sigma(\theta)=\theta$. 
 \bl
 \label{highest_root}
We have
\be
\label{long_root}
\iota(\theta)= \begin{cases}  
\theta_{\sigma,s} \quad (\fg, \fg_\sigma)\not= (A_{2n}, C_n)\\
\frac{1}{2}\theta_{\sigma}\quad  (\fg, \fg_\sigma)=(A_{2n}, C_n)\\
\end{cases}
\ee   
where $\theta_\sigma$ is the highest root of $\fg_\sigma$ and $\theta_{\sigma,s}$ is the highest short root of $\fg_\sigma$. Moreover, 

\be
\label{long_coroot}
\check{\iota}(\check{\theta})=\begin{cases}  
\check{\theta}_{\sigma}    \quad (\fg, \fg_\sigma)\not= (A_{2n}, C_n)\\
2\check{\theta}_{\sigma,s}  \quad  (\fg, \fg_\sigma)=(A_{2n}, C_n)
\end{cases}
\ee

where $\check{\theta}_{\sigma}$ (resp.\,$\check{\theta}_{\sigma,s}$) is the highest  coroot (the coroot of the highest root)  of $\fg_\sigma$.
\el
\bpf
We first determine $\check{\iota}(\check{\theta})$. Let $\check{\fg}$ be the Lie algebra with root system dual to the root system of $\fg$. We still denote by $\sigma$ the diagram automorphism on $\check{\fg}$ induced from the diagram automorphism $\sigma$ on $\fg$. It is well-known that the root system of $\fg_\sigma$ is dual to the root system of the fixed Lie algebra $\check{\fg}^{\sigma}$ (cf.\,\cite{Ho,HS}).

By \cite[Lemma 4.3]{Ho}, $\sigma$ acts on the highest root subspace $\check{\fg}_{\check{\theta}}$ by $1$ if $\fg$ is not of type $A_{2n}$; otherwise, $\sigma$ acts on $\check{\fg}_{\check{\theta}}$ by $-1$. It follows that if $\fg$ is not $A_{2n}$, then $\check{\fg}_{\check{\theta}}$ is the highest root subspace of the fixed point Lie subalgebra $\check{\fg}^\sigma$. Thus, in this case $\check{\iota}(\check{\theta})=\check{\theta}_{\sigma}$. When $\fg$ is of type $A_{2n}$, by \cite[Prop. 8.3]{Ka}  $\check{\iota}(\check{\theta})=2\check{\theta}_{\sigma,s}$. 

Finally, we can determine $\iota(\theta)$ from (\ref{root}) and \cite[Table 2, p.88]{Hu2}, and we get the formula (\ref{long_root}).
\epf

Note that $\check{\theta}_{\sigma}$ is the coroot of $\theta_{\sigma,s}$ and $\check{\theta}_{\sigma,s}$ is the coroot of $\theta_{\sigma}$.

\bl
\label{type_C}
Let $I_k$ be the Dynkin diagram of type $C_k$ with $k\geq 2$, where $I_k$ consists of vertices $i_1,i_2,\cdots, i_k$ such that the simple root $\alpha_{i_1}$ is a long root. Then the long root lattice $Q_{l}$ of $C_k$ is spanned by $\alpha_{i_1}, 2\alpha_{i_2},\cdots, 2\alpha_{i_k}$. 
\el
\bpf
For any $k\geq 1$, let $I_k$ be the Dynkin diagram of $C_k$  (where $C_1=A_1$), there exists a natural embedding $I_k\incl I_{k+1}$. Assume that  $I_{k}$ consists of vertices  $i_1,i_2,\cdots, i_k$, where the simple root $\alpha_{i_1}$ is the long root. 
  Let $\theta_k$ be the highest long root of $C_k$. Then $\theta_{k+1}-\theta_k=2\alpha_{k+1}$. Therefore the lattice of long roots of $C_k$ for $k\geq 2$, is spanned by $\alpha_{i_1}, 2\alpha_{i_2},\cdots, 2\alpha_{i_k}$. 
\epf

Let $Q^\sigma$ denote the lattice of $\sigma$-invariant elements in the root lattice $Q$ of $\fg$. 

\bl
\label{lattice_tran}

With respect to the isomorphism $\iota:  P_\sigma\simeq  P^\sigma$, we have
\[   
\iota(Q^\sigma)=\begin{cases}
Q_\sigma      \quad     \text{ if }  \fg \text{ is not of type } A_{2n}\\
P_\sigma=  \frac{1}{2} Q_{\sigma,l}  \quad  \text{ if }  \fg  \text{ is } A_{2n}
\end{cases}
   \]

where  $Q_{\sigma,l} $ is the lattice spanned by the long roots of $\fg_\sigma$. 

\el
\bpf
Clearly $Q^\sigma$ has a basis $\{ \sum_{i\in \imath} \alpha_i  \li   \imath\in I_\sigma   \}$. In view of  (\ref{root}), it is easy to see  when $\fg$ is not of type $A_{2n}$, $\iota(Q^\sigma)$ is the root lattice $Q_\sigma$ of $\fg_\sigma$. 

If $\fg$ is of type $A_2$, this is just a direct simple calculation. Otherwise, if $\fg= A_{2n}$ with $n\geq 2$, then  $\iota(Q^\sigma)=\sum_{\imath\in I_\sigma} a_\imath \bbZ \alpha_\imath$, where 
\[ a_\imath=\begin{cases}
1   \quad  \text{ if } \imath \text{ is not connected}\\
\frac{1}{2}  \quad \text{ if }  \imath \text{ is connected}
\end{cases}.
  \]

Let $\imath_0$ be the connected $\sigma$-orbit in $I$. Note that $\imath_0$ corresponds to the long root of $C_n$. By Lemma \ref{type_C}, our lemma follows.

\epf

  

\subsection{Affine Weyl groups  and  diagram automorphisms}
\label{affine_weyl_sect}
In this subsection, we refer to \cite{Hu1} the basics of affine Weyl groups and alcoves. 

Let $W$ be the Weyl group of $\fg$. The group $W$ acts on the weight lattice $P$. Let $P_\bbR$ be the space $P \otimes_\bbZ \bbR$. For each root $\alpha \in \Phi$, let $s_\alpha$ be the corresponding  reflection in $W$, i.e. for any $\lambda\in P_\bbR$, 
$s_\alpha(\lambda)=\lambda-\la \lambda, \check{\alpha} \ra \alpha $.
  

Let $W_\ell$ be the affine Weyl group $W\ltimes   \ell Q$ for any $\ell\in \bbQ$. Since $\fg$ is simply-laced, the Coxeter number is equal to the dual Coxeter number, moreover all roots have the same length. For any $\ell\in \bbN$, $W_\ell$ is the Weyl group of the affine Lie algebra $\hat{\fg}$ of level $\ell$. Let $s_0$ be the affine reflection $s_{\theta,1}$, i.e
\be
\label{ref_0}
 s_{\theta,1}(\lambda)=\lambda- (\la \lambda, \check{\theta} \ra -\ell) \theta, \ee
where $\theta$ is the highest root of $\fg$.
The affine Weyl group $W_\ell$ is a Coxeter group generated by  $\{s_i\, \li  \, i\in  \hat{I}\}$. 
For any $\alpha\in \Phi$ and $n\in \ell\bbZ$, the  hyperplane
\[ H_{\alpha,n}=\{\lambda \in P_\bbR  \li     \la \lambda, \check{\alpha}   \ra=n   \}  \]
is an affine wall of $W_\ell$.
  Every component of the complements of affine walls in $P_\bbR$ is an alcove. The affine Weyl group $W_\ell$ acts on the set of alcoves simply and transitively. Let $A_0$ be  the fundamental alcove, and it can be described as follows,
  \[   \{  \lambda\in P_\bbR  \li    \la \lambda, \check{\alpha}_i  \ra >0, \, \text{ for any } i\in I, \, \text{ and }  \la \lambda, \check{\theta}  \ra< \ell           \}. \]

The diagram automorphism $\sigma$ acts on $W$. Let  $W^\sigma$ be the fixed point group of $\sigma$ on $W$. Let $W_\sigma$ be the Weyl group of $\fg_\sigma$ with simple reflections $\{  s_\imath \li   \imath\in I_\sigma\}$. 
Then there exists an isomorphism $\iota: W^\sigma\simeq W_\sigma$ such that for any $\imath\in I_\sigma$, 
\be
\label{reflection}
\iota^{-1}( s_\imath)=\begin{cases} \prod_{i\in \imath} s_i    \quad   \text{ any  } i\not= j\in \imath, \, a_{ij}=0  \\
s_is_js_i  \quad   \text{ if } \imath=\{i,j\} \text{ and } a_{ij}=-1
      \end{cases}.
      \ee

The following lemma is obvious.
\bl
The group action of $W_\sigma$ on $P_\sigma$ is compatible with the action of $W^\sigma$ on $P^\sigma$, with respect to the isomorphisms $\iota: P^\sigma\simeq P_\sigma$ and $\iota: W^\sigma\simeq W_\sigma$.
\el

The diagram automorphism $\sigma$ also  acts naturally on $W_\ell$. Let $W^\sigma_\ell$ denote the fixed point group of $\sigma$ on $W_\ell$. 
It is easy to see that
\begin{equation}
\label{affine_Weyl_inv}
  W^\sigma_\ell= W^\sigma \ltimes  \ell Q^\sigma.
\end{equation}
Let $\hat{I}_\sigma$ be the set $I_\sigma\sqcup \{0\}$. We have the following lemma (cf.\,\cite[Section 5.2]{FSS}).
\bl
\label{affine_reflection}
$W^\sigma_\ell$ is a Coxeter group generated by $\{ \iota^{-1}(s_\imath)  \,|\, \imath\in \hat{I}_\sigma   \}$.
\el

 The group $W^\sigma_\ell$ naturally acts on $P^\sigma_\bbR$. Let $\CA$ denote the set of alcoves of $W_\ell$ in $P_\bbR$. There exists a natural action of $\sigma$ on $\CA$. 
   Let $\CA^\sigma$ be the set of $\sigma$-stable alcoves.

\bl
\label{fixed_point}
\ben
\item For any $A\in \CA^\sigma$, the set $A^\sigma$ is not empty, where $A^\sigma$ is the set of $\sigma$-invariant elements in $A$. 
\item For any two $\sigma$-stable alcoves $A$ and $A'$ in $\CA$, there exists a unique $w\in W^\sigma_\ell$ such that $w(A)=A'$.
\een
\el
\bpf
We first prove (1). For any $\lambda\in A$, $\lambda, \sigma(\lambda), \cdots  \sigma^{r-1}(\lambda)\in A$, where $r$ is the order of $\sigma$. By the convexity of $A$, 
\[ \frac{\lambda+ \sigma(\lambda)+\cdots  +  \sigma^{r-1}(\lambda)}{r} \in  A, \]
which is $\sigma$-invariant.

Now we prove (2). The affine Weyl group $W_\ell$ acts simply and transitively on $\CA$ (cf.\,\cite[\S 4.5]{Hu1}). Hence,
given any two elements $A,A'\in  \CA^\sigma$, there exists a unique $w\in W_\ell$ such that $w(A)=A'$. In particular, we have 
\[\sigma(w)(A)= \sigma w \sigma^{-1} (A)=\sigma w(A)=\sigma (A')=A'=w(A). \]
By the uniqueness of $w$, we have $\sigma(w)=w$. 

\epf

Let $P_{\sigma, \bbR}$ be the space $P_\sigma\otimes_\bbZ  \bbR$. We still denote by $\iota:   W^\sigma_\ell\simeq  W_\sigma \ltimes \iota(\ell Q^\sigma)$ the natural isomorphism of groups. By Lemma \ref{lattice_tran}, $W_\sigma \ltimes \iota(\ell Q^\sigma)$ is an affine Weyl group. In view of Lemma \ref{highest_root} and Lemma \ref{lattice_tran},
\[ A_{0,\sigma} =\{  \lambda\in  P_{\sigma,\bbR}  \li      \la\lambda, \check{\alpha}_\imath  \ra_\sigma>0\, \text{ for any } \imath\in I_\sigma, \, \text{ and }  \la\lambda, \check{\iota}(\check{\theta}) \ra_\sigma<\ell   \}        \]
is the fundamental alcove of $W_\sigma\ltimes \iota(\ell Q^\sigma) $. 

Let $\CA_\sigma$ be the set of alcoves of $W_\sigma \ltimes \iota(\ell Q^\sigma)$ in $P_{\sigma,\bbR}$.
\bp
\label{alcove_all_diagram}
\ben
\item The isomorphism $\iota:  P^\sigma_\bbR\simeq P_{\sigma, \bbR}$ induces a bijection  $\iota: (A_0)^\sigma\simeq A_{0,\sigma}$.
\item There exists a bijection  $\CA^{\sigma}\simeq \CA_\sigma$ with the map given by 
\[ A\mapsto   \iota(A^\sigma). \]
\item  For any $\lambda\in P_\bbR^\sigma$, $\lambda$ is in an affine wall of $W_\ell$ if and only if $\iota(\lambda)\in P_{\sigma,\bbR}$ is in an affine wall of $W_\sigma\ltimes \iota(\ell Q^\sigma)$.
\een
\ep
\bpf
We first prove (1). For any $\lambda\in P^\sigma_{\bbR}$, $\lambda\in  (A_0)^\sigma$ if and only if $\iota(\lambda)\in (A_0)^\sigma$, since 
\[ \la \lambda, \check{\alpha}_i \ra=\la \iota( \lambda), \check{\iota}(\check{\alpha}_i) \ra_\sigma=\la \iota( \lambda), \check{\alpha}_\imath \ra_\sigma>0, \]
for any $\imath\in I_\sigma$ and $i\in \imath$, and 
$ \la \lambda, \check{\theta}\ra=  \la  \iota( \lambda), \check{\iota}(\check{\theta}) \ra_\sigma<\ell. $  

The second part (2) of proposition follows from Lemma \ref{fixed_point} and the first part of the proposition. 
The third part (3) of the proposition follows from the first and second part  of proposition.
\epf

Let $\ell_\sigma:  W^\sigma_\ell\to  \bbN$ denote the length function on the Coxeter group $W^\sigma_\ell$. 
For any $\lambda\in \ell Q^\sigma$, let $\tau_\lambda$ be the translation on $P^\sigma_{ \bbR} $ by $\lambda$. The following lemma will be used in the proofs of Proposition \ref{wall_case} and Lemma \ref{kernel_describ} in Section \ref{Sect5}.
\bl
\label{even}
The length $\ell_\sigma(\tau_\lambda) $ is even. 
\el
\bpf
When $\fg$ is not $A_{2n}$, by Lemma \ref{lattice_tran} $\iota(Q^\sigma)=Q_\sigma$. The Coxeter group $W_\ell^\sigma$ is isomorphic to the affine Weyl group $W_\sigma \ltimes  Q_\sigma$. The problem is reduced to show that for any $\lambda\in Q_\sigma$, $\tau_\lambda$ has even length in $W_\sigma\ltimes  Q_\sigma$. 
If $\lambda$ is dominant, then $\ell_\sigma(\tau_\lambda)= \la\lambda, 2\check{\rho}_\sigma  \ra$ (cf.\,\cite{IM}), where $\check{\rho}_\sigma$ is the sum of all fundamental weights of $\fg_\sigma$. Hence $\ell_\sigma(\tau_\lambda)$ is even. For general $\lambda$, $\lambda=w(\lambda) $ for some $w\in W_\sigma$ and some dominant weight $\lambda^+\in Q_\sigma$. Then $\tau_\lambda=w \tau_{\lambda^+}w^{-1}$, and hence $\tau_\lambda$ is even. 

When $\fg$ is of type $A_{2n}$, by Lemma \ref{lattice_tran} $\iota(Q^{\sigma})= \frac{1}{2}Q_{\sigma,l}$. 
The normalized Killing form on $\fg_\sigma$ can identify $\frac{\ell}{2}Q_{\sigma,l}$ with $\frac{\ell}{2} \check{Q}_\sigma$ (cf.\,\cite[Proof of Lemma 9.3 (b)]{Be}), where $\check{Q}_\sigma$ is the coroot lattice of $\fg_\sigma$. This identification is compatible with the action of $W_\sigma$. Hence $W^\sigma_\ell$ is isomorphic to $W_\sigma \ltimes \check{Q}_\sigma$. By the same argument as above, the length $\ell_\sigma(\tau_\lambda) $ is also even.

\epf

\section{Conformal blocks }

\subsection{Affine Lie algebra}
\label{affine_lie_subsect}
Let $\fg$ be a simple Lie algebra. Let $\bbC((t))$ be the field of Laurent series over $\mathbb{C}$. 
Let $\tilde{\fg}$ be the associated affine Kac-Moddy algebra $\fg((t))\oplus \bbC c \oplus \bbC d$, where $\fg((t))$ denotes the loop Lie algebra $\fg\otimes_\bbC \bbC((t))$. The Lie bracket $[,]$ on $\tilde{\fg}$ is given by 
$$  [u\otimes f, v\otimes g]:=  [u,v]\otimes fg  +  (u|v) {\rm Res}_{t=0}(\frac{df}{dt} g)  c, $$
and  $[u\otimes t^n, d]=nu\otimes t^n, \quad   [d,c]=0,\quad  [u\otimes f, c]=0$,
for any $u,v \in \fg$ and $f, g\in \bbC((t))$, where $[u,v]$ is the Lie bracket on $\fg$ and $(\cdot|\cdot)$ is the normalized invariant bilinear form on $\fg$. For convenience, we identify $u\otimes 1$ with $u$ for any $u\in \fg$, and hence $\fg$ is naturally a Lie subalgebra of $\tilde{\fg}$. 

 Put $\hat{\fg}=\fg((t))\oplus \bbC c$. Clearly $\hat{\fg}$ is a Lie subalgebra of $\tilde{\fg}$.
 The affine Kac-Moody algebra $\tilde{\fg}$ corresponds to the extended Dynkin diagram $\hat{I}=I\sqcup\{0\}$ of $\fg$. The Cartan subalgebra $\tilde{\ft}$ associated to $\tilde{\fg}$ is $\ft\oplus \bbC c\oplus \bbC d$. 
  For any $\lambda\in P$  we view it as a weight of $\tilde{\fg}$ in the following way, $\lambda$ extends to $\tilde{\ft}$ such that $\lambda(d)=\lambda(c)=0$. Let $\delta$ be the linear functional on $\tilde{\ft}$ such that 
\[ \delta|_{\ft}=0,\quad   \delta(c)=0,\quad  \delta ( d)=1.\]  
 Let  $\alpha_0=-\theta+\delta$, where $\theta$ is the highest root of $\fg$. Then $ \{ \alpha_i\,|\, i\in \hat{I}\}$ is the set of simple roots of $\tilde{\fg}$. The fundamental weight $\Lambda_0$ of $\tilde{\fg}$ is given by  the linear functional on $\tilde{\ft}$ such that 
\[  \Lambda_0|_\ft=0,\quad \Lambda_0(c)=1,\quad \Lambda_0(d)=0.\]
\subsection{Affine Weyl groups and Weyl groups of affine Kac-Moody algebras}
 In the following we describe the relationship between the affine Weyl groups of  simple Lie algebras and the Weyl groups of affine Kac-Moody algebras. For more details, one can refer to \cite[\S 6]{Ka}.
  These two different perspectives are both crucial in this paper.

 Let $\hat{W}$ be the Weyl group of the affine Kac-Moody algebra $\tilde{\fg}$ (cf.\,\cite[$\S$3.7]{Ka}). Set $\tilde{\ft}^*_\bbR:=P_\bbR+\bbR \Lambda_0+\bbR \delta$. The Weyl group $\hat{W}$ acts on $\tilde{\ft}^*_\bbR$. Note that $\hat{W}$ keep $\delta$ invariant (cf.\,\cite[$\S$6.5]{Ka}). Hence the Weyl group $\hat{W}$ acts on $\hat{P}_{\bbR, \ell} $ for any $\ell\in \bbR$, where 
 \[\hat{P}_{\bbR, \ell} :=  \{ x\in \tilde{\ft}^*_\bbR  \li   \la x, c  \ra=\ell   \} / \bbR \delta. \]

With respect to the isomorphism  $P_\bbR \simeq  \hat{P}_{\bbR, \ell} $
 given by $\lambda\mapsto  \lambda+ \ell \Lambda_0$, we have the following lemma (cf.\,\cite[$\S$6.5,$\S$6.6]{Ka}).


 \bl
 \label{dot_lem}
 There exists an isomorphism ${\rm af}:\hat{W}\simeq  W_\ell$ of groups such that 
 for any  $\Lambda=\lambda+\ell \Lambda_0\in \hat{P}_{\bbR, \ell}$ and $w\in \hat{W}$, the following formula holds,
 \[  w\cdot \Lambda={\rm af}(w)\cdot \lambda+\ell \Lambda_0 \text{ in } \hat{P}_{\bbR, \ell}. \]
 \el
%
 
 Let $\hat{\rho}$ be the sum $\sum_{i\in \hat{I}} \Lambda_i$ of all fundamental weights of $\tilde{\fg}$. By \cite[\S 6.2.8]{Ka}, $\hat{\rho}=\rho+\check{h}\Lambda_0$  where $\rho$ is the sum $\sum_{i\in I}\omega_i$ of all fundamental weights of $\fg$, and $\check{h}$ is the dual Coxeter number of $\fg$. 
 
 We define $\star$ action of $\hat{W}$ on $ \hat{P}_{\bbR, \ell}$ as follows,
 \[ w\star \Lambda=w\cdot (\Lambda+\hat{\rho})-\hat{\rho}, \quad  w\in \hat{W}, \Lambda\in  \hat{P}_{\bbR, \ell}. \]
Similarly, we still denote by  $\star$ the following action of $W_\ell$ on $P_\bbR $,
 \[ w\star \lambda=w\cdot(\lambda+\rho)-\rho, \quad  w\in W_\ell, \lambda\in P_\bbR. \]
\bl
\label{star_leml}
Given $\Lambda=\lambda+\ell \Lambda_0\in  \hat{P}_{\bbR, \ell} $ and $w\in \hat{W}$, we have 
\[ w\star \Lambda={\rm af}(w)\star \lambda+ \ell \Lambda_0, \quad \text{ where ${\rm af}(w)\in W_{\ell+\check{h}}$. }\]
\el
\bpf
It follows from Lemma \ref{dot_lem} and the formula  $\hat{\rho}=\rho+\check{h}\Lambda_0$.
\epf

\subsection{Diagram automorphisms as intertwining operators of representations } 
\label{rep_lie}

We denote by $V_\lambda$ the irreducible representation of $\fg$ of highest weight $\lambda$ for each $\lambda\in P^+$. From now on we always fix a highest weight vector $v_\lambda\in V_\lambda$ for each $\lambda$. There exists a unique operator  $\sigma:V_\lambda\to V_{\sigma(\lambda)}$ such that 
  $\sigma(v_\lambda)=v_{\sigma(\lambda)}$, and   $\sigma(u\cdot v)=\sigma(u)\cdot\sigma(v)$ for any $u\in \fg$ and $v\in V_\lambda$.

When $\sigma(\lambda)=\lambda$, $\sigma$ acts on $V_\lambda$. Given any $\sigma$-invariant dominant weight of $\fg$ and any $r$-th root of unity $\xi\in  \bbC$ where $r$ is the order of $\sigma$, we denote by $V_{\lambda,\xi}$ the representation of $\fg\rtimes \la \sigma \ra$, i.e. it consists of  $V_\lambda$ as representation of $\fg$ and an operator $\sigma: V_\lambda\to V_\lambda$ such that  $\sigma$ acts on $v_\lambda$ by $\xi$, and $\sigma(u\cdot v)=\sigma(u)\cdot \sigma(v)$.

 Given a tuple $\vec{\lambda}=(\lambda_1,\cdots,\lambda_k)$ of dominant weights of $\fg$. We denote by $V_{\vec{\lambda}}$ the tensor product $V_{\lambda_1}\otimes \cdots   \otimes V_{\lambda_k}$. Denote by $V_{\vec{\lambda}}^\fg$ the  invariant space of $\fg$ on $V_{\vec{\lambda}}$. The collection of operators $\{  \sigma: V_\lambda\to  V_{\sigma(\lambda)}\}$ induce 
\[ \sigma:  V_{\vec{\lambda}}\to  V_{\sigma(\vec{\lambda})}, \quad    \sigma:  V^\fg_{\vec{\lambda}}\to  V^\fg_{\sigma(\vec{\lambda})}, \]
where $\sigma(\vec{\lambda})=(\sigma(\lambda_1), \cdots, \sigma(\lambda_k))$.

Recall the set $P_\ell$ defined in (\ref{weight_level}). 
The following lemma is well-known (cf.\,\cite[$\S$12.4]{Ka}).
 \bl
 \label{dominant_finite_affine}
 For any $\lambda\in P^+$ and $\ell\in  \bbN$,
   $\lambda+\ell \Lambda_0$ is a dominant weight of $\tilde{\fg}$ if and only if $\lambda\in P_\ell$.   
 \el

  For any $\lambda\in P_\ell$, let $\hat{M}(V_\lambda)$ dnote the generalized Verma module $U(\hat{\fg})\otimes_{U(\hat{\fp})} V_\lambda$  of $\hat{\fg}$, where $\hat{\fp}=\fg[[t]]\oplus \bbC\cdot c$ acts on $V_\lambda$ by evaluating $t=0$ and $c$ acts by $\ell$. Then the unique maximal irreducible quotient $\mathcal{H}_\lambda$ is an irreducible integrable representation of $\hat{\fg}$ of level $\ell$. The action of $\hat{\fg}$ on  $\CH_\lambda$ extends uniquely to the irreducible integrable representation of $\tilde{\fg}$ of highest weight $\lambda+\ell\Lambda_0$ by letting  $d$ act trivially on the highest weight vectors. 
  
  From the construction of $\CH_\lambda$, there exists a natural inclusion $V_\lambda\to  \CH_\lambda$. Denote by $\bar{v}_\lambda$ the image of $v_\lambda\in V_\lambda$ in $\CH_\lambda$, which is again a highest weight vector in $\CH_\lambda$.
 
 The diagram automorphism $\sigma:\fg\to \fg$ extends to an automorphism on $\hat{\fg}$ (by abuse of notation we still use $\sigma$) such that $\sigma(u\otimes f)=\sigma(u)\otimes f$ for any $u\in \fg$ and $f\in \bbC((t))$, and $\sigma(c)=c$.
As in the case of $V_\lambda$, there exists a unique operator $\sigma:  \CH_\lambda\to  \CH_{\sigma(\lambda)}$ such that  
  $\sigma(\bar{v}_\lambda)=\bar{v}_{\sigma(\lambda)}$, and  $\sigma(X\cdot v)=\sigma(X)\sigma(v)$ for any $X\in \hat{\fg}$ and $v\in \CH_\lambda$.
 In particular $\sigma$ acts on $\CH_\lambda$ when $\sigma(\lambda)=\lambda$. As in the case of $V_\lambda$, for any $\sigma$-invariant dominant weight $\lambda$ of $\fg$ and for any $r$-th root of unity $\xi$, we denote by $\CH_{\lambda,\xi}$ the representation of $\hat{\fg}\rtimes \la \sigma\ra $ which satisfies similar conditions for $V_{\lambda,\xi}$.
 
 Given a tuple $\vec{\lambda}$ of dominant weights, denote by $\CH_{\vec{\lambda}}$ the tensor product $\CH_{\lambda_1}\otimes \cdots  \otimes \CH_{\lambda_k}$. The operators $\{\sigma:  \CH_\lambda\to  \CH_{\sigma(\lambda)}\}$ induce the operator  $\sigma:  \CH_{\vec{\lambda}}\to  \CH_{\sigma(\vec{\lambda})}$ such that 
 \[ \sigma(v_1\otimes \cdots  \otimes v_k )= \sigma(v_1)\otimes \cdots  \otimes \sigma(v_k), \]
 for any $v_i\in \CH_{\lambda_k}$, $i=1,\cdots, k$.
 
 The inclusion $V_\lambda\incl \CH_\lambda$  is compatible with the diagram automorphism, i.e. 
 \begin{equation}
 \label{VH_equ}
 \xymatrix{
V_\lambda  \ar[r] \ar[d]^{\sigma}  & \CH_\lambda  \ar[d]^{\sigma}\\
V_{\sigma(\lambda)}  \ar[r]   &  \CH_{\sigma(\lambda)}   
 }.
 \end{equation}
 
 Let $\hat{\fg}^-$ denote the Lie subalgebra $t^{-1}\fg[t^{-1}]$. We denote by $(\CH_\lambda)_{\hat{\fg}^- }$  the coinvariant space of $\CH_\lambda$ with respect to the action of $\hat{\fg}^-$. The Lie algebra $\fg$  acts naturally on $(\CH_\lambda)_{\hat{\fg}^- }$.
 The following lemma is well-known.
 \bl
\label{affine_finite}
As representations of $\fg$, we have a natural isomorphism  $V_\lambda\simeq (\CH_\lambda)_{\hat{\fg}^- }$. Moreover the following diagram commutes
 \begin{equation}
 \label{VH_equ}
 \xymatrix{
V_\lambda  \ar[r] \ar[d]^{\sigma}  &  (\CH_\lambda)_{\hat{\fg}^-}  \ar[d]^{\sigma}\\
V_{\sigma(\lambda)}  \ar[r]   &  ( \CH_{\sigma(\lambda)})_{\hat{\fg}^- }   
 }.
 \end{equation}
 \el

Let $\tau$ be the Cartan involution of $\fg$ such that $\tau(e_i)=-f_i, \tau(f_i)=-e_i, \tau(h_i)=-h_i$, where $e_i,f_i,h_i$ for $i\in I$, are Chevalley generators of $\fg$. Then  $\tau$ is an automorphism on $\fg$. For any finite dimensional representation $V$ of $\fg$, by composing $\tau$ we can redefine a new representation structure on $V$, $ X * v := \tau(X)\cdot v $, 
for any $X\in \fg$ and $v\in V$. We denote  by $V^\tau$ this $\tau$-twisted representation. 

 For any dominant weight $\lambda$, let $\lambda^*$ be the dominant weight $-\omega_0(\lambda)$ where $\omega_0$ is the longest element in the Weyl group $W$. The space   $V_\lambda^\tau$ is isomorphic to $V_{\lambda^*}$ as representations of $\fg$.

 The Cartan involution  $\tau$ on $\fg$ extends to an automorphism on $\hat{\fg}$ (by abuse of notation we still use $\tau$) such that $\tau(u\otimes f)=\tau(u)\otimes f$ and $\tau(c)=c$ for any $u\in\fg, f\in \C((t))$. 
      Denote by $\CH_\lambda^\tau$ the representation of $\hat{\fg}$ by composing the automorphism $\tau:  \hat{\fg}\to  \hat{\fg}$. Then $\CH_\lambda^\tau\simeq \CH_{\lambda^*}$. 
      
     Summarizing the above discussions, we have the following lemma.
   \bl
   \ben
   \item There exists a unique $\bbC$-linear isomorphism $\tau_\lambda:  V_\lambda\to  V_{\lambda^*}$ such that 
   \[ \tau_\lambda(v_\lambda)=v_{\lambda^*},\quad   \tau_\lambda(u\cdot v)=\tau(u)\cdot \tau_\lambda(v), \quad \text{for any $u\in \fg$ and $v\in V_\lambda$.}   \]
   \item There exists a unique $\bbC$-linear isomorphism $\tau_\lambda:  \CH_\lambda\to  \CH_{\lambda^*}$ such that 
   \[ \tau_\lambda(v_\lambda)=v_{\lambda^*},\quad   \tau_\lambda(X\cdot v)=\tau(X)\cdot \tau_\lambda(v), \quad \text{   for any $X\in \hat{ \fg}$ and $v\in \CH_\lambda$.}   \]
   \een
   \el

   The isomorphism $\tau_\lambda: V_\lambda\to  V_{\lambda^*}$ for each $\lambda$ induces an isomorphism $\tau_{\vec{\lambda}}:  V_{\vec{\lambda}}^\fg\to V_{\vec{\lambda}^* }^\fg$ for any tuple of dominant weights $\vec{\lambda}$. Since for any weight $\lambda$, we have $\sigma(\lambda^*)=\sigma(\lambda)^*$, and $\sigma\circ\tau=\tau\circ \sigma$, we have the following lemma.
  \bl
  \label{inv_intertwiner} Let  $\vec{\lambda}^*$ denote $(\lambda_1^*,\cdots, \lambda_k^*)$.
  The following diagram commutes:
\[
\xymatrix{
V^\fg_{\vec{\lambda}  }  \ar[d]^{\sigma}  \ar[r]^{\tau }  & V^\fg_{ \vec{\lambda}^* }  \ar[d]^{\sigma} \\
V^\fg_{\sigma(\vec{\lambda} ) }  \ar[r]^{\tau}  &  V^\fg_{ \sigma( \vec{\lambda})^*  } 
},
 \]
  \el

\subsection{Conformal blocks and diagram automorphisms}

\label{conformal_block}

A $k$-pointed projective curve consists of a projective curve $C$ and $k$-distinct smooth points $\vec{p}=(p_1,\cdots, p_k)$ in $C$.
Given a  $k$-pointed projective curve $(C, \vec{p})$, we associate a dominant weight $\lambda_i\in P_\ell$ to each point $p_i$. 
Let $\fg(C\backslash \vec{p})$ be the space of $\fg$-valued regular functions on $C\backslash \vec{p}$. The space $\fg(C\backslash \vec{p})$ is naturally a Lie algebra induced from $\fg$. There exists a Lie algebra action of $\fg(C\backslash \vec{p})$ on $\CH_{\vec{\lambda}}$, and 
 the space  $V_{\fg, \ell, \vec{\lambda}}(C,\vec{p})$ of  conformal blocks  associated to $\vec{p}$ and $\vec{\lambda}$ is defined as follows: 
$$ V_{\fg, \ell, \vec{\lambda}}(C,\vec{p}):=  (\CH_{\vec{\lambda}})_{\fg(C\backslash \vec{p})}=\CH_{\vec{\lambda}}/ \fg(C\backslash \vec{p})\CH_{\vec{\lambda}}.$$

Let $\tau_{\vec{\lambda}}:  \CH_{\vec{\lambda}}\to  \CH_{\vec{\lambda}^*}$ be the $\bbC$-linear isomorphism $\tau_{\lambda_1}\otimes  \cdots \otimes \tau_{\lambda_k}$. 
  The map $\tau_{\vec{\lambda}}$ descends to an isomorphism on the space of  conformal blocks
\[  \tau_{\vec{\lambda}}: V_{\fg, \ell, \vec{\lambda}}(C,\vec{p})\to  V_{\fg, \ell, \vec{\lambda}^* }(C,\vec{p}). \]

\bl
\label{sigma_tau}
We have the following commutative diagram:
\[
\xymatrix{
V_{\fg, \ell, \vec{\lambda}}(C,\vec{p}) \ar[d]^{\sigma}  \ar[r]^{\tau_{\vec{\lambda}} }  & V_{\fg, \ell, \vec{\lambda}^* }(C,\vec{p})\ar[d]^{\sigma} \\
V_{\fg, \ell, \sigma(\vec{\lambda}) }(C,\vec{p}) \ar[r]^{\tau_{\vec{\lambda}}}  &  V_{\fg, \ell, \sigma(\vec{\lambda}^*)   }(C,\vec{p})
}.
 \]
\el
\bpf
The automorphism $\sigma$ commutes with the automorphism $\tau$ on $\fg$, i.e. $\tau\circ \sigma= \sigma\circ \tau$. Then commutativity easily follows.
\epf

\bp
\label{Propagation_prop}
Let $\vec{p}=\{p_1,p_2,\cdots, p_s \}$, $\vec{q}=\{ q_1,q_2,\cdots, q_t \}$ be two finite nonempty subsets smooth points of $C$, without common points; let $\lambda_1,\cdots, \lambda_s, \mu_1,\cdots, \mu_t$ be elements in $P_\ell$. We let $\fg(C\backslash  \vec{p})$ act on $V_{\mu_j}$ through the evaluation map $X\otimes f\mapsto  f(q_j)X$. The inclusions $V_{\mu_j}\incl \CH_{\mu_j}$ induce an isomorphism 
\begin{equation}
\label{propagation_equ}
  (\CH_{\vec{\lambda }}\otimes  V_{\vec{\mu}}   )_{\fg(C\backslash  \vec{p})} \simeq  (\CH_{\vec{\lambda}}\otimes  \CH_{\vec{\mu}}  )_{\fg(C\backslash \vec{p}\cup \vec{q} )} =V_{\fg, \ell, \vec{\lambda},\vec{\mu}}(C,\vec{p},\vec{q}), 
  \end{equation} 
and this isomorphism is  compatible with the diagram automorphism $\sigma$, i.e. the following diagram commutes
\begin{equation}
\label{propagation_com}
\xymatrix
{
 (\CH_{\vec{\lambda }}\otimes  V_{\vec{\mu}}   )_{\fg(C\backslash  \vec{p})}  \ar[r]^{\simeq}   \ar[d]^{\sigma} &   V_{\fg, \ell, \vec{\lambda},\vec{\mu}}(C,\vec{p},\vec{q})  \ar[d]^{\sigma}  \\
 (\CH_{\sigma(\vec{\lambda } ) } \otimes  V_{ \sigma(\vec{\mu} )  }   )_{\fg(C\backslash  \vec{p})}  \ar[r]^{\simeq}  \ar[r] & V_{\fg, \ell, \sigma(\vec{\lambda}),\sigma(\vec{\mu} ) }(C,\vec{p},\vec{q})
}.
\end{equation}
\ep
\bpf
Isomorphism (\ref{propagation_equ}) is a well-known theorem (cf.\,\cite[Proposition 2.3]{Be}). The commutativity of diagram (\ref{propagation_com}) follows from the commutativity (\ref{VH_equ}).
\epf

When $\vec{q}=q$ and $\mu=0$. Isomorphism (\ref{propagation_equ}) is the so-called ``propogation of vacua". Proposition \ref{Propagation_prop}  shows that the propagation of vacua   is  compatible with the action of the diagram automorphism. 


\bl
\label{basic}
\ben
\item  For any $p\in \bbP^1$, one has  $V_{\fg, \ell}(\bbP^1) \simeq V_{\fg,\ell, 0}(\bbP^1,p) \simeq \bbC$ by $1$, and the automorphism $\sigma$ acts on $V_{\fg, \ell}(\bbP^1)$ and $V_{\fg, \ell,0}(\bbP^1,p) $ by $1$.
\item For any $p\neq q$ in $\bbP^1 $, one has $V_{\fg,\ell,\lambda, \lambda^*}(\bbP^1,p,q )\simeq (V_\lambda\otimes V_{\lambda^*})_{\fg}=\bbC$. As a consequence, the automorphism $\sigma$ acts on $V_{\fg,\ell,\lambda, \lambda^*}(\bbP^1,p,q )$ by $1$ for any $\sigma$-invariant dominant weight $\lambda$.
\een
\el
\bpf
  By Proposition \ref{Propagation_prop}, there exists a map $\bbC\to  V_{\bbP^1}(p;0)$ compatible with the action of $\sigma$ where $\bbC$ is viewed as a trivial representation of $\fg$ and $\sigma$ acts on $\bbC$ trivially. By \cite[Corollary 4.4]{Be}, this map is an isomorphism. By  Proposition \ref{Propagation_prop} again, $V_{\fg, \ell}(\bbP^1) \simeq V_{\fg,\ell, 0}(\bbP^1,p)  \simeq \bbC$ and this isomorphism is also compatible with the action of $\sigma$. Hence, $\sigma$ acts on $V_{\fg, \ell}(\bbP^1)$  and   $V_{\fg,\ell, 0}(\bbP^1,p) $   by $1$. This proves $(1)$.
  
  Similarly, by  Proposition \ref{Propagation_prop} there exists a map  $(V_\lambda\otimes V_{\lambda^*})_\fg\to V_{\fg,\ell,\lambda, \lambda^*}(\bbP^1,p,q )$ which is compatible with the action of $\sigma$. This map is an isomorphism in view of \cite[Corollary 4.4]{Be}. On the other hand, it is easy to see that $\sigma$ acts on 
$(V_\lambda\otimes V_{\lambda^*})_\fg$ by $1$. Hence it also acts on $V_{\fg,\ell,\lambda, \lambda^*}(\bbP^1,p,q )$ by $1$.
\epf



Given a  stable $k$-pointed curve $(C,\vec{p})$. Assume that $q\in C$ is a nodal point in $C$. 
 Let $\pi:\tilde{C}\to C$ be the normalization of $C$ at $q$. Denote by $\{q_+,q_-\}$ the preimage of $q$ via $\pi$. 
  Without confusion, we will still denote by $p_1,\cdots, p_k$ the preimages of $p_1,p_2,\cdots,p_k\in C$ in $\tilde{C}$.

We  choose a system of $\fg$-equivariant maps $\kappa_\mu: \bbC\to  V_\mu\otimes V_{\mu^*} $ for $\mu\in P^+$, 
 such that  the following diagram commutes 
\[
\xymatrix{
\bbC   \ar[r]^{\kappa_\mu}  \ar[rd]_{\kappa_{\sigma(\mu)}}  &   V_\mu\otimes V_{\mu^*}   \ar[d]^{\sigma}  \\
&   V_{\sigma(\mu)}\otimes V_{\sigma(\mu^*) }   
},
\]
for any dominant weight $\mu$. 
Note that the map $\kappa_\mu$ induces the following map 
\[ \hat{\kappa}_\mu:  V_{\fg, \ell, \vec{\lambda}}(C,\vec{p})\simeq V_{\fg, \ell, \vec{\lambda},0}(C,\vec{p},q)\to V_{\fg, \ell, \vec{\lambda},\mu,\mu^*}(C,\vec{p},q_+,q_-).\]

Moreover, it is easy to see that  the following  diagram commutes 
\be
\label{fact_map}
\xymatrix{
V_{\fg, \ell, \vec{\lambda}}(C,\vec{p})   \ar[r]^<<<<<<{\hat{\kappa}_\mu}  \ar[d]^{\sigma}  & V_{\fg, \ell, \vec{\lambda},\mu,\mu^*}(C,\vec{p},q_+,q_-)\ar[d]^{\sigma}  \\
V_{\fg, \ell, \sigma(\vec{\lambda})  }(C,\vec{p})      \ar[r]^<<<<{\hat{\kappa}_{\sigma(\mu) } }  &   V_{\fg, \ell, \sigma(\vec{\lambda}),\sigma(\mu),\sigma(\mu^*)}(C,\vec{p},q_+,q_-)
}.
\ee

\bt
\label{normalization}

The  map
\be
\label{factorization_equ}
\xymatrix{
V_{\fg, \ell, \vec{\lambda}}(C,\vec{p})  \ar[r]^<<<<{(\hat{\kappa}_\mu )}  &  \bigoplus_{\mu \in P_\ell}  V_{\fg, \ell, \vec{\lambda},\mu,\mu^*}(\tilde{C},\vec{p},q_+,q_-)
}
\ee
is an isomorphism. Moreover the following diagram commutes,
\be
\label{factorizatio_com}
\xymatrix{
V_{\fg, \ell, \vec{\lambda}}(C,\vec{p})  \ar[r]^<<<<<<{(\hat{\kappa}_\mu )}   \ar[d]^{\sigma} &  \bigoplus_{\mu \in P_\ell}  V_{\fg, \ell, \vec{\lambda},\mu,\mu^*}(\tilde{C},\vec{p},q_+,q_-)  \ar[d]^{\sigma} \\
V_{\fg, \ell, \sigma(\vec{\lambda}) }(C,\vec{p})  \ar[r]^<<<<{(\hat{\kappa}_\mu )}    &        \bigoplus_{\mu \in P_\ell}  V_{\fg, \ell, \sigma(\vec{\lambda}),\sigma(\mu),\sigma(\mu^*)}(\tilde{C},\vec{p},q_+,q_-) 
}.
\ee
\et
\bpf
Isomorphism (\ref{factorization_equ}) is the well-known factorization theorem on conformal blocks (cf.\,\cite[Theorem 3.19]{Ue} ), and the commutativity (\ref{factorizatio_com}) follows from the commutativity (\ref{fact_map}).
\epf

Recall that $P_\ell^\sigma$ is the set of $\sigma$-invariant dominant weights in $P_\ell$.
\bco
\label{norm_trace}
With the same setup as above. If $\sigma(\vec{\lambda})=\vec{\lambda}$, then the following equality holds
\[ \tr(\sigma|V_{\fg, \ell, \vec{\lambda}}(C,\vec{p})  )=\sum_{\mu\in P_\ell^\sigma }\tr(\sigma| V_{\fg, \ell, \vec{\lambda},\mu,\mu^*}(\tilde{C},\vec{p},q_+,q_-)   ). \]
\eco
\bpf
This is an immediate consequence of Theorem \ref{normalization}. 
\epf

Given a family $(\pi:\CC\to X, \vec{\fp})$ of stable $k$-pointed curves  where $\pi$ is a family of projective curves with at most nodal singularities over a smooth variety $X$  and $\vec{\fp}=(\fp_1,\cdots,\fp_k)$ is a collection of sections $\fp_i: X\to \CC$ with disjoint images such that $\fp_i(x)$ is a smooth point in $C_x:=\pi^{-1}(x)$ for each $i$ and $x\in X$, one can attach a sheaf of conformal blocks 
 $\CV_{ \fg,\ell, \vec{\lambda}}(\CC,\vec{\fp})$ on $X$ which is locally free and of finite rank, see \cite{Loo} for the coordinate-free approach to the sheaf of conformal blocks. For each $x\in X$, the fiber
 $\CV_{ \fg,\ell, \vec{\lambda}}(\CC,\vec{\fp})|_x$  is the space of conformal blocks 
 $V_{ \fg,\ell, \vec{\lambda}}(C_x,\vec{\fp}(x))$, where $\vec{\fp}(x)=(\fp_1(x),\cdots,\fp_k(x))$ are the $k$-distinct smooth points in $C_x$ as the image of $x$ via $\vec{\fp}$.


 From the construction the sheaf of conformal blocks (cf.\,\cite{Loo}), one can see the diagram automorphism $\sigma$ acts algebraically on $\CV_{ \fg,\ell, \vec{\lambda}}(\CC,\vec{\fp})$. 
 Denote by $\la \sigma \ra$ the cyclic group generated by $\sigma$. Then the group $\la \sigma \ra$  is isomorphic to $\bbZ/r\bbZ$, where $r$ is the order of $\sigma$.

\bl
\label{constant_lem}
For any family $(\pi:\CC\to X, \vec{\fp})$  of stable pointed curves, 
the function  $x\in X\mapsto   \tr(\sigma|  V_{ \fg,\ell, \vec{\lambda}}(C_x,\vec{\fp}(x)) )$ is constant. 
\el
\bpf
Given any irreducible representation $\rho$ of $\la \sigma \ra$, we denote by ${\rm ch}(\rho)$ and ${\rm ch}(V_{ \fg,\ell, \vec{\lambda}}(C_x,\vec{\fp}(x)) )$ the characters of $\rho$ and $V_{ \fg,\ell, \vec{\lambda}}(C_x,\vec{\fp}(x))$ respectively as representations of $\la \sigma \ra$. For any two functions $\phi, \psi$ on $\la \sigma \ra$, we define the bilinear form
\[ (\phi, \psi)=\frac{1}{r}\sum_{i=0}^{r-1}  \phi(\sigma^i)\psi(\sigma^{-i}), \]
where $r$ is the order of $\sigma$. For any $x\in X$, let $m_{\rho}(x)$ be the multiplicity of $\rho$ appearing in $V_{ \fg,\ell, \vec{\lambda}}(C_x,\vec{\fp}(x))$. 
 By representation theory of finite groups, we have
  \[m_\rho(x)= ({\rm ch}\rho, {\rm ch} V_{ \fg,\ell, \vec{\lambda}}(C_x,\vec{\fp}(x)) ).\]
This  is a continuous function on $X$ with integer values. This is forced to be constant. Hence 
\[\tr(\sigma|V_{ \fg,\ell, \vec{\lambda}}(C_x,\vec{\fp}(x)))=\sum m_\rho(x) \tr(\sigma|\rho) \]
is constant along  $x\in X$.
\epf

The following theorem shows that the trace of the diagram automorphism on the space of  conformal blocks satisfies factorization properties. 
\bt
\label{factorization}
\ben
\item  
For any  stable $k$-pointed curve $(C,\vec{p})$, let $\vec{\lambda}$ be a tuple of dominant weights in $P_\ell^\sigma$ attached to $\vec{p}$. Then the value $\tr(\sigma|V_{ \fg,\ell, \vec{\lambda}}(C,\vec{p}))$ only depends on  $\vec{\lambda}$ and the genus of $C$.
\item   
Given a  stable $k$-pointed curve $(C,\vec{p})$ of genus $g\geq 1$ and a stable $(k+2)$-pointed curve $(C', \vec{q})$ of genus $g-1$. We have the following formula
\[  \tr(\sigma|V_{ \fg,\ell, \vec{\lambda}}(C,\vec{p}) )=\sum_{\mu\in P_\ell^\sigma} \tr(\sigma|V_{ \fg,\ell, \vec{\lambda},\mu,\mu^* }(C',\vec{q})),\]
where a tuple $\vec{\lambda}=(\lambda_1,\cdots, \lambda_k)$ of dominant weights in $P_\ell^\sigma$, is attached to $\vec{p}$ and the first $k$ points of $\vec{q}$.
\item
Given any tuples of dominant weights $\vec{\lambda}=(\lambda_1,\lambda_2,\cdots, \lambda_s)$ and $\vec{\mu}=(\mu_1,\cdots, \mu_t)$ in $P_\ell^\sigma$ where $s,t\geq 2$, 
we have the following equality
\[  \tr(\sigma| V_{ \fg,\ell, \vec{\lambda}, \vec{\mu} } ( \bbP^1, \vec{p}_1) ) = \bigoplus_{\nu\in P^\sigma_\ell}   \tr(\sigma|  V_{ \fg,\ell, \vec{\lambda}, \nu  } ( \bbP^1, \vec{p}_2) ) \tr(\sigma| V_{ \fg,\ell, \vec{\mu}, \nu^* } ( \bbP^1, \vec{p}_3) ), \]
where $\vec{p}_1$ is any tuple of $s+t$ distinct points, $\vec{p}_2$ is any tuple of $s+1$ distinct points and $\vec{p}_3$ is any tuple of $t+1$ distinct points in $\bbP^1$. 
\een
\et
\bpf
We first prove part (1). By the standard theory of moduli of curves (cf.\,\cite[Theorem 2.15]{HM}), there exists a chain of  families of stable $k$-pointed curves over  smooth bases connecting any two  stable $k$-pointed curves with the same genus. In view of Lemma \ref{constant_lem}, $(1)$ follows. 

From the theory of moduli of curves again (cf.\,\cite[Theorem 2.15]{HM}) and the dimension formula for the space of nodal curves with fixed nodal types (see the discussions after \cite[Theorem 2.15]{HM}), we know that any smooth pointed stable curve can be degenerated to  an irreducible stable pointed curve with only one nodal point. Then part (2) follows from part (1) and Corollary \ref{norm_trace}. 

We now proceed to prove  part (3).
Let $C$ be the union of two projective lines $C=C_1\cup C_2$ where $C_1$ and $C_2$ intersect at  one point $z$. Let $\vec{p}=(p_1,\cdots, p_s)$ be a set of  $s$ distinct points in $C_1\backslash \{z\}$ and $\vec{q}=\{q_1,\cdots, q_t\}$ be another set of  $t$ distinct points in $C_2\backslash \{z\}$ where $s, t\geq 2$. Clearly $(C,\vec{p}\cup \vec{q})$ is a stable $(s+t)$-pointed curve of genus zero. Again by the theory of moduli of curves, there exists a family  $\pi: \CC\to X$ of stable $(s+t)$-pointed curves  over a smooth variety  $X$ such that $C_{x_0}=C$ with $\vec{p}\cup \vec{q}$ and any other fiber is a projective line with a tuple $\vec{p}_1$ of $s+t$ points. By Lemma \ref{constant_lem}, 
\[ \tr(\sigma|V_{ \fg,\ell, \vec{\lambda}, \vec{\mu} } ( C, \vec{p},\vec{q}) )=\tr(\sigma| V_{ \fg,\ell, \vec{\lambda}, \vec{\mu} } ( \bbP^1, \vec{p}_1 ). \]    
Let $\pi:  \tilde{C}\to C$ be the normalization of $C$ at $z$ with the preimage $(z_+,z_-)$ of $z$. The pointed  curve $(\tilde{C}, \vec{p},\vec{q},z_+,z_-)=(\bbP^1, \vec{p},z_+)\sqcup (\bbP^1,\vec{q},z_-)$ is a  disjoint union of  a $(s+1)$-points projective line and a $(t+1)$-pointed projective line. Finally, part (3) follows from Corollary \ref{norm_trace} and Lemma \ref{constant_lem}. 
\epf

\begin{remark}
By Theorem \ref{factorization}, the computation of the trace of the diagram automorphism on the space of conformal blocks can be reduced to the trace of the diagram automorphism on the space of conformal blocks on the pointed curve $(\bbP^1, (0,1,\infty))$.
\end{remark}
\subsection{$\sigma$-twisted fusion ring}
\label{twisted_fusion_sect}
Let $J$ be a finite set with an involution $\lambda\mapsto \lambda^*$. We denote by $\bbN^J$ the free commutative monoid generated by $J$, that is, the set of sums $\sum_{\lambda\in J}  n_\lambda  \lambda$ with $n_\lambda\in \bbN$. The involution of $J$ extends by linearity to an involution $x\mapsto  x^*$ of $\bbN^J$. 
We first recall the definition of fusion rule (cf.\,\cite[\S 5]{Be}). 
\bd
\label{fusion_rule}
A fusion rule on $J$ is a map $N: \bbN^J\to  \bbZ$ satisfying the following conditions:  
\ben
\item  One has $N(0)=1$, and $N(\lambda)>0$ for some $\lambda\in J$;
\item  $N(x^*)=N(x)$ for every $x\in \bbN^J$;
\item  For $x,y\in \bbN^J$, one has $N(x+y)=\sum_{\lambda} N(x+\lambda)N(y+\lambda^*)$.
\een
\ed
The kernel of a fusion rule $N$ by definition is the set of elements $\lambda\in J$ such that $N(\lambda+x)=0$ for all $x\in \bbN^I$. A fusion rule on $J$ is called non-degenerate if the kernel is empty. 

\bl
\label{trace_integer}
 If $\sigma(\vec{\lambda})=\vec{\lambda}$, then the trace $\tr(\sigma| V_{\fg, \ell, \vec{\lambda}}(C,\vec{p}))$ is an integer.
\el
\bpf
When the order of $\sigma$ is $2$, this is obvious. In general, it follows from Theorem \ref{trace_conformal_tensor}, Formula (\ref{HS_0}) in the introduction and part (3) of Theorem \ref{factorization}.
\epf

\bt
\label{trace_fusion_rule}
The map $\tr_\sigma:  \bbN^{P_\ell^\sigma}\to  \bbZ$ given by 
\[ \sum \lambda_i \mapsto   \tr(\sigma| V_{\fg, \ell, \vec{\lambda}}(\bbP^1,\vec{p})), \]
where $\vec{\lambda}=(\lambda_1,\cdots, \lambda_k)$ and $\vec{p}=(p_1,\cdots, p_k)$ is the set of  any $k$-distinct points in $\bbP^1$,
 is a  non-degenerate fusion rule. Here the set $P^\sigma_\ell$ is equipped with the involution $\lambda\mapsto \lambda^*:=-w_0(\lambda)$, where $w_0$ is the longest element in the Weyl group $W$.
\et
\bpf
By Lemma \ref{trace_integer}, the trace map $\tr_\sigma$ indeed always takes integer values. 

Condition (1) of Definition \ref{fusion_rule} follows from part (1) of Lemma \ref{basic}. Condition (2) follows from Lemma \ref{sigma_tau}. Condition (3) follows from  part (3) of Theorem \ref{factorization}. 
 The non-degeneracy follows from part (2) of Lemma \ref{basic}.

\epf

Let $R_\ell(\fg,\sigma)$ be a free abelian group with the set $P_\ell^\sigma$ as a basis. 
As a consequence of Theorem \ref{trace_fusion_rule} and  \cite[Proposition 5.3]{Be}, we can define a ring structure on $R_\ell(\fg,\sigma)$  by putting 
\begin{equation} 
\label{fusion_product}
 \lambda \cdot \mu :=  \sum_{\nu\in P_\ell^\sigma}   \tr(\sigma| V_{\fg,\ell, \lambda, \mu,\nu^* } ( \bbP^1,0,1,\infty) ) \nu, \quad \text{for any $\lambda,\mu\in P^\sigma_\ell$.
}
 \end{equation}

Let $S_\sigma$ be the set of characters (i.e. ring homomorphisms) of $R_\ell(\fg,\sigma)$ into $\bbC$. The following proposition is a consequence of general facts on fusion ring by Beauville \cite[Corollary 6.2]{Be}. 
\bp 
\ben
\item $R_\ell(\fg,\sigma)\otimes \mathbb{C}$ is a reduced commutative ring. 
\item  The map $R_\ell(G,\sigma)\otimes \bbC \to  \bbC^{S_\sigma}$ given by $x\mapsto (\chi(x))_{x\in S_\sigma}$ is an isomorphism of $\bbC$-algebras. 
\item  We have $\chi(x^*)=\overline{\chi(x)}$, where $\overline{\chi(x)}$ denotes the complex conjugation of $\chi(x)$ for any $\chi\in S_\sigma$ and $x\in R_\ell(\fg, \sigma)$.
\een
\ep
Let $\omega_\sigma$ be the Casimir element in $R_\ell(\fg,\sigma)$ defined as follows
\be
\label{Casimir}
 \omega_\sigma=\sum_{\lambda\in P_\ell^\sigma}  \lambda\cdot \lambda^*.
 \ee
\bp
\label{machine}For any $k$-pointed stable curve $(C,\vec{p})$ and for any $\sigma$-invariant tuple $\vec{\lambda}$ of dominant weights in $P_\ell$,
we have the following formula
\[ \tr(\sigma|V_{\fg, \ell, \vec{\lambda}}(C,\vec{p}))=\sum_{\chi\in S_\sigma}    \chi(\lambda_1)\cdots \chi(\lambda_k) \chi(\omega_\sigma)^{g-1}, \]
where $g$ is the genus of $C$ and  $\chi(\omega_\sigma)=\sum_{\lambda\in P_\ell^\sigma} |\chi(\lambda)|^2$.
\ep
\bpf
This is a consequence of  part (2) of Theorem \ref{factorization} and \cite[Proposition 6.3]{Be}.
\epf
From this proposition, if we can determine the set $S_\sigma$ and the value $\chi(\omega_\sigma)$ for each $\chi\in S_\sigma$, then the trace $\tr(\sigma|V_{\fg, \ell, \vec{\lambda}}(C,\vec{p}))$ is known. 

\section{Sign problems}



\subsection{ Borel-Weil-Bott theorem on the affine flag variety}
\label{sign_BWB_sect}

Let $G$ be a connected and simply-connected  simple algebraic group associated to a simple Lie algebra $\fg$.
  Let $G((t))$ be the loop group of $G$, and let $\hat{G}$ be the nontrivial central extension of $G((t))$ by the center $\bbC^\times$. Then $\hat{\fg}$ is the Lie algebra of  $\hat{G}$. Let  $\tilde{G}$ be the group $\tilde{G}=\hat{G}\rtimes \bbC^\times $ whose Lie alegbra is the affine Kac-Moody algebra $\tilde{\fg}$.


Let $\CI$ be the Iwahori subgroup of $G((t))$, i.e. $\CI=\ev_0^{-1}(B)$, where $B$ is the Borel subgroup of $G$. 
Let $\Fl_G$ be the affine flag variety $G((t))/\CI$ of $G$. Let $\hat{\CI}$ be the group $ \CI\times \bbC^\times$, where $\bbC^\times$ is the center of $\hat{G}$. Let $\tilde{\CI}$ be the product $\hat{\CI}\rtimes  \bbC^\times$ as subgroup of $\tilde{G}$. 
Then we have 
\[ \Fl_G\simeq  \hat{G}/ \hat{\CI} \simeq \tilde{G}/\tilde{\CI}.\]

Given any algebraic representation $V$ of $\tilde{\CI}$, 
we can attach a $\tilde{G}$-equivariant vector bundle $\CL(V)$ on $\Fl_G$ as  $ \CL({V}):= \tilde{G}\times_{\tilde{\CI}} V^*$, 
where $V^*$ is the dual representation of $\tilde{\CI}$. Let $\Lambda$ be a character of $\tilde{\CI}$ and let $\bbC_\Lambda$ be the associated $1$-dimensional representation of $\tilde{\CI}$. We denote by $\CL(\Lambda)$ the $\tilde{G}$-equivariant line bundle $\CL(\bbC_\Lambda)$ on $\Fl_G$.

For any ind-scheme $X$ and any vector bundle $\CF$ on $X$, the cohomology groups $H^*(X,\CF)$ carry a topology. We put $H^*(X,\CF)^\vee$ the restricted dual of $H^*(X,\CF)$, i.e. $H^*(X,\CF)^\vee$  consists of continuous functional on $H^*(X,\CF)$ where we take discrete topology on $\bbC$. The affine flag variety $\Fl_G$  is an ind-scheme of ind-finite type. We refer the reader to \cite{Ku1} for the foundation of flag varieties of Kac-Moody groups. 

  Recall the following affine analogue of Borel-Weil-Bott theorem (cf.\,\cite[Theorem 8.3.11]{Ku1}). 
\bt
\label{BWB}
Given any dominant weight $\Lambda$ of $\tilde{G}$ and any  $w\in \hat{W}$, the space
$ H^{\ell(w)}(\Fl_G, \CL({w\star \Lambda }) )^\vee $
is naturally  the integrable irreducible representation $\CH_{\Lambda}$ of $\tilde{\fg}$ of highest weight $\Lambda$, where $w\star \Lambda=w\cdot(\Lambda+\hat{\rho})-\hat{\rho}$ and $ H^{\ell(w)}(\Fl_G, \CL({w\star \Lambda }) )$ is the cohomology of the line bundle $\CL({w\star \Lambda })$ on $\Fl_G$. Moreover, $ H^{i}(\Fl_G, \CL({w\star \Lambda }) )=0 $ if $i\neq  \ell(w)$.
\et

Let $\sigma$ be a diagram automorphism on $G$. It induces an action on $G((t))$ by acting trivially on $t$. It also  induces actions on $\hat{G}$ and $\tilde{G}$ by acting trivially on the the center and degree component.
 Note that $\sigma$ preserves $\tilde{\CI}$. 
For any $\sigma$-invariant character $\Lambda$ of $\tilde{\CI}$, we have a natural $\sigma$-equivariant structure on $\CL(\Lambda)$, since 
\[ \tilde{G}\rtimes \la \sigma \ra \times_{\tilde{\CI}\rtimes \la \sigma \ra } (\bbC_\Lambda)^*\simeq \tilde{G}\times_{\tilde{\CI}} ( \bbC_\Lambda)^*, \]
where the action of $\sigma$ on $\bbC_\Lambda$ is by the scalar $1$. 
Let $\xi$ be an $r$-th root of unity, where $r$ is the order of $\sigma$. We denote by  $\CL(\Lambda, \xi)$ the following $\tilde{G}\rtimes \la \sigma \ra$-equivariant line bundle, 
\[  \CL(\Lambda, \xi):=\tilde{G}\rtimes \la \sigma \ra\times_{ \tilde{\CI} \rtimes \la \sigma \ra  }  (\bbC_{\Lambda,\xi})^*   \]
where $\hat{\CI}$ acts on   $\bbC_{\Lambda, \xi}$ by $\Lambda$ and $\sigma$ acts on $\bbC_{\Lambda,\xi}$ by $\xi$. By this convention the natural $\tilde{G}\rtimes \la \sigma \ra$-equivariant structure on $\CL(\Lambda)$ is isomorphic to $\CL(\Lambda, 1)$.


For any $\sigma$-orbit $\imath$ in the affine Dynkin diagram $\hat{I}$, let $G_\imath$ be the  simply-connected algebraic group associated to the sub-diagram $\imath$ and let $B_\imath$ be the 
Borel subgroup of $G_\imath$. We have the following possibilities
\[ G_\imath=\begin{cases}
\SL_2  \quad  \quad    \imath=\{i\}\\
\SL_2\times \SL_2  \quad   \imath=\{i,j\} \text{ and } i, j \text{ are not connected }\\
\SL_2\times \SL_2\times \SL_2    \quad   \imath=\{i,j,k\}   \text{ and } i, j,k  \text{ are not connected }\\
\SL_3   \quad   \quad   \imath=\{i,j\}  \text{ and } i, j \text{ are  connected }
\end{cases}. \]

 We still denote by $\sigma$ the diagram automorphism on $G_\imath$ which preserves $B_\imath$. Any $\sigma$-invariant  weight $\lambda$ of $G_\imath$ can be written as $n\rho_\imath$ for some integer $n\in \bbZ$, where $\rho_\imath$ is the sum of all fundamental weights of $G_\imath$. Let $\CB_\imath:=G_\imath/B_\imath$ be the flag variety of $G_\imath$. Put  $d_\imath=\dim  G_\imath/B_\imath$. 
  
 As in the affine case for any $r$-th root of unity and any $\sigma$-invariant character $\lambda$ of $B_\imath$, we set 
 \[ \CL(\lambda, \xi)=G_\imath\times_{B_\imath} ( \bbC_{\lambda,\xi})^*   \]
 as a $G_\imath\rtimes \la \sigma \ra$-equivariant line bundle on $\CB_\imath$. Let $\Omega_\imath$ be the canonical bundle of $\CB_\imath$. 
   Note that the canonical  bundle  $\Omega_\imath$ is naturally a $G_\imath\rtimes \la \sigma \ra$-equivariant line bundle. 
  \bl
  \label{canonical_sign}
We have the following isomorphism of $G_\imath\rtimes \la \sigma\ra$-equivariant line bundles
$\Omega_{ \imath }\simeq  \CL(-2\rho_\imath, \epsilon_\imath)$,
 where  $\epsilon_\imath=(-1)^{d_\imath-1}$.
 \el
 \bpf
The canonical bundle $\Omega_{\imath}$ is naturally isomorphic to $G_\imath\times _{B_\imath}  (\lwedge^{d_\imath}( \fg_\imath/\fb_\imath) )^*$, where $\fg_\imath$ (resp.\,$\fb_\imath$) is the Lie algebra of $G_\imath$ (resp $B_\imath$). Hence, it suffices to determine the action of $T_\imath$ and $\sigma$ on $ \lwedge^{d_\imath}( \fg_\imath/\fb_\imath) $, where $T_\imath$ is the maximal torus of $G_\imath$ contained in $B_\imath$. Note that 
\[  \lwedge^{d_\imath}( \fg_\imath/\fb_\imath) \simeq  \lwedge^{d_\imath}  \fn^-_\imath, \]
where $\fn^-_\imath$ is the nilpotent radical of the negative Borel subalgebra of $\fg_\imath$. Hence, as $1$-dimensional representation of $T_\imath$, it is isomorphic to $-2\rho_\imath$, and by case-by-case analysis it is easy to check $\sigma$ acts on it exactly by $\epsilon_\imath$. This finishes the proof of the lemma.
 \epf
 
 Only when $\imath$ consists of two vertices and $\imath=\{i,j\}$ is not connected, $\epsilon_\imath=-1$; otherwise $\epsilon_\imath=1$.

\bl
\label{orbit_case}Given any $n\in \bbZ$ and any $r$-th root of unity $\xi$, 
there exists a unique isomorphism up to a scalar
\[   
H^{d_\imath}( \CB_\imath, \CL(n\rho_\imath, \xi))\simeq H^0(\CB_\imath, \CL((-n-2)\rho_\imath, \epsilon_\imath\cdot \xi)),
\]
as representations of $G_\imath\rtimes \la \sigma\ra$. Moreover,
\[ H^k(\CB_\imath, \CL(n\rho_\imath, \xi))=0  \quad  \text{ if }  k\neq 0, d_\imath. \]
\el
\bpf
By Borel-Weil-Bott theorem we have the following isomorphism of representations of $G_\imath\rtimes \la\sigma  \ra$
 \be
 \label{orbit_BWB}
  H^{0}( \CB_\imath, \CL(n\rho_\imath, \xi))^*= \bc V_{n\rho_\imath,\xi}   \quad n \geq 0   \\
 0  \quad \quad   n<0
  \ec, \ee
 for any $n\in \bbZ$ and $r$-th root of unity $\xi$, where $V_{n\rho_\imath,\xi}$ is the irreducible representation of $G_\imath$ of highest weight $n\rho_\imath$ with the compatible action of $\sigma$ which acts on the highest weight vectors by $\xi$. 

  By Serre duality  we have the following canonical isomorphism
 \be
 \label{serre_imath}
     H^{d_\imath}( \CB_\imath, \CL(n\rho_\imath,\xi))\simeq  H^0( \CB_\imath, \CL(-n\rho_\imath,\xi^{-1})\otimes  \Omega_{ \CB_\imath})^* 
 \ee
 as representations of $G_\imath\rtimes \la \sigma \ra$. 
In view of Lemma \ref{canonical_sign}, 
\[ H^0( \CB_\imath, \CL(-n\rho_\imath,\xi^{-1})\otimes  \Omega_{ \CB_\imath}) \simeq  H^0( \CB_\imath, \CL((-n-2)\rho_\imath, \epsilon_\imath\cdot \xi^{-1})  ).\]
In view of (\ref{orbit_BWB}), by Schur lemma there exists a unique isomorphism up to a scalar   
\begin{equation}
\label{cite_lem4.4}
 H^0( \CB_\imath, \CL((-n-2)\rho_\imath, \epsilon_\imath\cdot \xi^{-1})  )^* \simeq H^0( \CB_\imath, \CL((-n-2)\rho_\imath,\epsilon_\imath\cdot \xi)  ) \end{equation}
as representations of $G_\imath\rtimes \la \sigma \ra$.
Therefore we have an isomorphism
\be
\label{natural_isom}
 H^{d_\imath}( \CB_\imath, \CL(n\rho_\imath,\xi))\simeq H^0( \CB_\imath, \CL((-n-2)\rho_\imath, \epsilon_\imath\cdot \xi)  )  \ee
as representations of $G_\imath\rtimes \la \sigma \ra$. 

Now we prove the second part of the lemma. When $n\geq 0$, $n\rho_\imath$ is dominant, then Borel-Weil-Bott theorem implies that  $ H^k(\CB_\imath, \CL(n\rho_\imath, \xi))=0$   unless $k=0$. In view of isomorphism (\ref{natural_isom}), when $n\leq -2$, $H^{k}(\CB_\imath, \CL(n\rho_\imath, \xi))=0$ unless $k=d_\imath$. When $n=-1$, we have $s_i\star \rho_\imath=\rho_\imath$. Thus, $H^{k}(\CB_\imath, \CL(-\rho_\imath, \xi))=0$ for any $k$. 
\epf

Let $\tilde{\CP}_\imath$ be the parabolic subgroup of $\tilde{G}$ containing $\tilde{\CI}$ and $G_\imath$. We have an isomorphism of varieties $\tilde{\CP}_\imath/\tilde{\CI}\simeq \CB_\imath$. 
Let $\pi_\imath: \Fl_G\to  \tilde{G}/ \tilde{\CP}_\imath$ be the projection map. The fiber is isomorphic to $\CB_\imath$. There exists the following natural isomorphism as $\tilde{G}\rtimes \la \sigma \ra$-equivariant ind-schemes
\[\tilde{G}\rtimes \la \sigma  \ra \times_{\tilde{\CP}_\imath \rtimes \la \sigma \ra}  \CB_\imath\simeq  \Fl_G. \]

From the $G_\imath \rtimes \la\sigma \ra$-equivariant line bundle $\CL(n\rho_\imath,\xi)$ on $\CB_\imath$, by descent theory one can attach a $\tilde{G}\rtimes \la \sigma\ra$-equivariant line bundle $\CL_{\pi_\imath}(n\rho_\imath,\xi)$ on $\Fl_G$, i.e.
\[\CL_{\pi_\imath}(n\rho_\imath,\xi):=  \tilde{G}\rtimes \la \sigma  \ra \times_{\tilde{\CP}_\imath \rtimes \la \sigma \ra} \CL(n\rho_\imath, \xi),\] 
where the action of $\tilde{\CP}_\imath\rtimes \la \sigma \ra$ on $\CL(n\rho_\imath, \xi)$ factors through $G_\imath\rtimes \la \sigma \ra$. 
Let $\Omega_{\pi_\imath}$ be the relative canonical line bundle of $\Fl_G$ over $\tilde{G}/\tilde{\CP}_\imath$. By Lemma \ref{canonical_sign} as a $\tilde{G}\rtimes \la\sigma \ra$-equivariant bundle, we have 
\begin{equation}
\label{can_bun}
 \Omega_{\pi_\imath}\simeq  \CL_{\pi_\imath}(-2\rho_\imath,\epsilon_\imath). \end{equation}

Let $R^k(\pi_\imath)_*$ be the $k$-th derived functor of the pushforward functor $(\pi_\imath)_*$. The following lemma is a relative version of Lemma \ref{orbit_case}.
\bl
\label{family_orbit_case}
There exits a natural isomorphism of $\tilde{G}\rtimes \la \sigma \ra$-equivariant vector bundles
\[  R^{d_\imath}(\pi_\imath)_*(\CL_{\pi_\imath}(n\rho_\imath,\xi)) \simeq ( \pi_\imath)_*(\CL_{\pi_\imath}((-n-2)\rho_\imath,\xi\cdot \epsilon_\imath) ). \]
\el
\bpf
By relative Serre duality for the morphism $\pi_\imath: \Fl_G\to  \tilde{G}/ \tilde{\CP}_\imath$, there exists a canonical isomorphism of $\tilde{G}\rtimes \la \sigma  \ra$-equivariant sheaves on $ \tilde{G}/ \tilde{\CP}_\imath$,
\[   R^{d_\imath}(\pi_\imath)_*(\CL_{\pi_\imath}(n\rho_\imath,\xi)) \simeq  (\pi_\imath)_*(\CL_{\pi_\imath}(-n\rho_\imath,\xi^{-1} ) \otimes   \Omega_{\pi_\imath}  ) ^\vee, \]
where $\vee$ denotes the dual of coherent sheaf. From isomorphism (\ref{can_bun}), it gives rise to the following isomorphism of $\tilde{G}\rtimes \la \sigma  \ra$-equivariant sheaves on $ \tilde{G}/ \tilde{\CP}_\imath$,
\begin{equation}
\label{rel_serre}
   R^{d_\imath}(\pi_\imath)_*(\CL_{\pi_\imath}(n\rho_\imath,\xi)) \simeq  (\pi_\imath)_*(\CL_{\pi_\imath}((-n-2)\rho_\imath,\xi^{-1}  \epsilon_\imath )  ) ^\vee. \end{equation}
We look at the fiber of the sheaf $(\pi_\imath)_*(\CL_{\pi_\imath}((-n-2)\rho_\imath,\xi^{-1}  \epsilon_\imath )  ) ^\vee$  at the base point $e\tilde{\CP }_\imath\in  \tilde{G}/ \tilde{\CP}_\imath $. This is the representation of $\tilde{\CP }_\imath\ltimes \la \sigma\ra$ on $H^0( \CB_\imath, \CL((-n-2)\rho_\imath, \epsilon_\imath\cdot \xi^{-1})  )^* $
by factoring through the map $\tilde{\CP }_\imath\ltimes \la \sigma\ra\to G_\imath\ltimes \la \sigma\ra$. From isomorphism (\ref{cite_lem4.4}), the $\tilde{G}\rtimes \la \sigma  \ra$-equivariance gives rise to an isomorphism of $\tilde{G}\rtimes \la \sigma  \ra$-vector bundles
\[ (\pi_\imath)_*(\CL_{\pi_\imath}((-n-2)\rho_\imath,\xi^{-1}  \epsilon_\imath )  ) ^\vee\simeq  ( \pi_\imath)_*(\CL_{\pi_\imath}((-n-2)\rho_\imath,\xi\cdot \epsilon_\imath) ). \]
Combining with (\ref{rel_serre}), the lemma follows. 
\epf

By Lemma \ref{affine_reflection}, the affine Weyl group $(\hat{W})^\sigma$ consists of simple reflections $\{s_{\imath}\li  \imath\in \hat{I}_\sigma\}$. 

\bl
\label{orbit_refl}
For any $\sigma$-invariant weight $\Lambda$ of $\tilde{\fg}$ and  for any $\sigma$-orbit in $\hat{I}$, we have
\[
s_\imath \cdot\Lambda=  \begin{cases}   
\Lambda-\la \Lambda, \check{\alpha}_i \ra \sum_{i\in \imath} \alpha_i  \quad   \text{ if } a_{ij}=0 \text{ for any } i\not =j\in  \imath,\\
\Lambda- 2\la \Lambda, \alpha_i\ra   (\alpha_i+\alpha_j)      \quad   \text{ if }   \imath=\{i,j\} \text{ is connected}.\\
\end{cases}
\]
\el
\bpf
For any $\sigma$-orbit $\imath$ in $I$, this is routine to check, in particular we use the formula (\ref{root}). When $\imath=\{0\}$, this is simply the definition of $s_0$. \epf

\bp
\label{small_case}
For any $\sigma$-invariant  weight $\Lambda$, and for any $\sigma$-orbit in the affine Dynkin diagram $\hat{I}$  and any $r$-th root of unity $\xi$, we have the following isomorphism 
\[  H^{k+d_\imath}(\Fl_G, \CL(s_\imath \star \Lambda, \xi ) )\simeq  H^k (\Fl_G, \CL(\Lambda, \epsilon_\imath \cdot \xi))    \]
as representations of $\tilde{G}\rtimes \la\sigma \ra$, for all integer $k$.
\ep
\bpf

Note that 
the restriction  $\CL(\Lambda, \xi)|_{\CB_\imath}$ of the $\tilde{G}\rtimes \la \sigma \ra$-equivariant line bundle $\CL(\Lambda, \xi)$ to the fiber $\CB_\imath$ is isomorphic to $\CL(\la \Lambda, \check{\alpha}_i  \ra \rho_\imath, \xi)$ as a $G_\imath \rtimes \la \sigma \ra$-equivariant line bundle for any $i\in  \imath$. Note that for any $i,j\in \imath$, $\la \Lambda, \check{\alpha}_i  \ra=\la \Lambda, \check{\alpha}_j  \ra$.
In view of Lemma \ref{orbit_refl}, we have
\[ s_\imath\star \Lambda=\begin{cases}
\Lambda- (\la \Lambda, \check{\alpha}_i \ra+1) \sum_{i\in \imath} \alpha_i  \quad   \text{ if } \imath \text{ is not connected }\\
\Lambda- 2(\la \Lambda, \alpha_i\ra+1)   (\alpha_i+\alpha_j)      \quad   \text{ if }   \imath=\{i,j\} \text{ is connected}\\
\end{cases}.
 \]
Hence  for any $\sigma$-orbit $\imath$ in $\hat{I}$ and $i\in \imath$, we have
\[  \la s_\imath\star \Lambda, \check{\alpha}_i  \ra= - \la \Lambda, \check{\alpha_i}  \ra  - 2. \]
It follows that  
\[ \CL(s_\imath\star \Lambda, \xi)|_{\CB_\imath}=\CL( -(\la\Lambda, \check{\alpha}_i  \ra+2)  \rho_\imath, \xi  ). \]

By Lemma \ref{family_orbit_case}, we have the following natural isomorphism of $\tilde{G}\rtimes \la \sigma \ra$-equivariant vector bundles
\be
\label{eq_a}
 R^{d_\imath} \pi_* (  \CL(\Lambda, \xi   ) )\simeq  \pi_*  ( \CL( s_\imath\star \Lambda, \epsilon_\imath\cdot \xi ).
\ee

By Lemma \ref{orbit_case}, we have 
\be
\label{eq_b}
 R^k\pi_*  (   \CL(\Lambda, \xi   ))=0  \quad   \text{ if } k\neq 0, d_\imath. 
 \ee

In view of (\ref{eq_a}) and (\ref{eq_b}),
 Leray's spectral sequence implies that
\[  H^{k+d_\imath}(\Fl_G, \CL(\Lambda, \xi ) )\simeq  H^k(\Fl_G, \CL(\Lambda, \epsilon_\imath \cdot \xi))     \]
as representations of $\tilde{G}\rtimes \la \sigma \ra$.
\epf

For any $w\in(\hat{W})^\sigma$, put
\begin{equation}
\label{def_epsilon_w}
 \epsilon_{w}=(-1)^{\ell(w)-\ell_\sigma(w)}. \end{equation} 
 For any reduced expression $w=s_{\imath_k}s_{\imath_{k-1}}\cdots s_{\imath_1}$ of $w$ in the Coxeter group $(\hat{W})^\sigma$ where $\imath_1,\cdots, \imath_k$ are $\sigma$-orbits in $\hat{I}$  and each $s_{\imath}$ is defined in (\ref{reflection}) for any $\imath\in I_\sigma$ and $s_{\{0\}}=s_0$, we have 
$\epsilon_w=\epsilon_{\imath_k}\cdots \epsilon_{\imath_1}$,
where $\epsilon_\imath$ is introduced in Lemma \ref{canonical_sign}. 

Finally, we are now ready to prove the following theorem. 
\bt
\label{sign_thm1}
For any $w\in (\hat{W})^\sigma$  and for any $\sigma$-invariant dominant weight $\Lambda$ of $\tilde{G}$. 
We have the following isomorphism of  representations of  $\tilde{G}\rtimes \la \sigma \ra$
\[  H^{\ell(w)} (\Fl_G, \CL(  w\star \Lambda,\xi) )\simeq     H^0(\Fl_G, \CL( \Lambda, \epsilon_w\cdot \xi  )   )  . \] 
\et
\bpf
We can write $w=s_{\imath_k}s_{\imath_{k-1}}\cdots s_{\imath_1} $ as a reduced expression in the Coxeter group $(\hat{W})^\sigma$, where $\imath_1,\cdots, \imath_k$ are $\sigma$-orbits in $\hat{I}$. Then  
\[\Lambda, s_{\imath_1}\star \Lambda, (s_{\imath_2}s_{\imath_1})\star \Lambda, \cdots, w\star  \Lambda   \]
are all $\sigma$-invariant  weights of $\tilde{G}$.

Note that  as an element in $\hat{W}$, the length  $\ell(w)$ of $w$ is equal to $\sum_{i=1}^k  d_{\imath_i} $. 
In view of Proposition \ref{small_case}, we get a chain of isomorphisms of $\tilde{G}\rtimes \la \sigma \ra$-representations
\begin{align*}
 H^{\ell(w)}(\Fl_G, \CL(w\star \Lambda, \xi))  &\simeq H^{\ell(w)-d_{\imath_1}}(\Fl_G, \CL( (s_{\imath_1}  w)\star \Lambda, \epsilon_{\imath_1} \xi))\\
 &\simeq H^{\ell(w)-d_{\imath_1}-d_{\imath_{2}}}(\Fl_G, \CL((s_{\imath_2}s_{\imath_1}w )\star \Lambda, \epsilon_{\imath_2}\epsilon_{\imath_1} \xi))\\
 &\cdots \cdots \\
 &\simeq H^0(\Fl_G, \CL(\Lambda, \epsilon_{w} \cdot \xi)). 
 \end{align*}

This finishes the proof of the theorem.
\epf

For any dominant weight $\Lambda$ of $\tilde{\fg}$ and an $r$-th root of unity, as always we  denote by $\CH_{\Lambda,\xi}$ the irreducible integrable representation of $\tilde{\fg}$ of highest weight $\Lambda$ together with a compatible action of $\sigma$ which acts on the highest weight vectors of $\CH_{\Lambda,\xi}$ by $\xi$. 
\bco
\label{affine_rep_sign}
In the same setting as in Theorem \ref{sign_thm1}, we have the following isomorphism of representations of $\tilde{\fg}\rtimes \la \sigma \ra$,
\[   H^{\ell(w)}(\Fl_G, \CL(w\star \Lambda, \xi))^\vee\simeq \CH_{\Lambda,\epsilon_w\cdot \xi} . \]
\eco
\bpf
This is an immediate consequence of Theorem \ref{BWB} and Theorem  \ref{sign_thm1}.
\epf

\begin{remark}
For any $\sigma$-invariant weight $\lambda$ of $G$, let $\CL(\lambda)$ be the associated line bundle on $G/B$. By Borel-Weil-Bott theorem, 
$H^i(G/B, \CL(\lambda))$ carries an action of the diagram automorphism. The action was determined by Naito. Theorem \ref{sign_thm1} and Theorem \ref{sign_Gr} are the affine analogues of the results of Naito \cite{N1}. 
\end{remark}

\subsection{Borel-Weil-Bott theorem on affine Grassmannian}
\label{sign_affine}
For any weight $\lambda$ of $G$, 
let $\CL_\ell(\lambda)$ be the $\hat{G}$-equivariant line bundle on $\Fl_G$ defined as follows, 
\[ \CL_\ell(\lambda) := \hat{G}\times_{\hat{\CI}} I_\ell(\bbC_{\lambda})^*, \]
where $ I_\ell(\bbC_{\lambda})$ is the 1-dimensional representation of $\hat{\CI}$ such that $\CI$ factors through the character $\lambda: B\to \bbC^\times$ and the center $\bbC^\times$ acts by $t\mapsto  t^\ell$, and $I_\ell(\lambda)^*$ is the dual of $I_\ell(\lambda)$ as the representation of $\CI$. 

 For any character $\Lambda$ of $\tilde{\CI}$, if $\Lambda=\lambda+\ell \Lambda_0$ where $\Lambda$ is a  weight of $\tilde{G}$ and $\lambda$ is a  weight of $G$, then as $\hat{G}$-equivariant line bundles, $\CL(\Lambda)=  \CL_\ell(\lambda)$.

If $\lambda$ is $\sigma$-invariant, then $\CL_\ell(\lambda)$ has a natural $\sigma$-equivariant structure as in the case of $\CL(\Lambda)$. Similarly,  to an $r$-th root of unity $\xi$ where $r$ is the order of $\sigma$, we can associate a $\hat{G}\rtimes \la \sigma \ra$-equivariant line bundle $\CL_\ell(\lambda,\xi)$. 
 If $\Lambda=\lambda+\ell\Lambda_0$ where $\lambda\in P^\sigma$, then  $\CL(\Lambda, \xi)=\CL_\ell(\lambda,\xi)$ as $\hat{G}\rtimes \la \sigma  \ra$-equivariant line bundles. 
 
 Recall  from Lemma \ref{dominant_finite_affine}, the weight $\Lambda=\lambda+\ell \Lambda_0$ is dominant for $\tilde{G}$ if and only if $\lambda$ is dominant for $G$ and $\la \lambda, \check{\theta}\ra\leq \ell$. 
 Recall the affine Weyl group $W_{\ell+\check{h}}$ discussed in Section \ref{affine_weyl_sect}, the action of $W_{\ell+\check{h}}$ on the weight lattice $P$ of $G$ is compatible with the action of $\hat{W}$ on the space of weights of $\tilde{G}$ of level $\ell$, see Lemma \ref{dot_lem} and Lemma \ref{star_leml}. 
  Therefore we can translate Theorem \ref{sign_thm1} into the following equivalent theorem.

\bt
\label{sign_thm1_equiv}
For any $w\in W_{\ell+\check{h}}$ such that $\sigma(w)=w$  and for any $\sigma$-invariant dominant weight $\lambda\in P_\ell$, 
we have the following isomorphism
\[  H^{\ell(w)} (\Fl_G, \CL_\ell(  w\star \lambda,\xi) )\simeq     H^0(\Fl_G, \CL_\ell( \lambda, \epsilon_w\cdot \xi  )   )   \] 
as  representations of  $\hat{G}\rtimes \la \sigma \ra$.
\et

 Let $\hat{\CP}$ be the subgroup $G[[t]]\times \bbC^\times$ of $\hat{G}$ where $\bbC^\times$ is the center torus. 
The affine Grassmannian $\Gr_G:=G((t))/G[[t]]$ is isomorphic to the partial flag variety $\hat{G}/\hat{\CP}$. For any finite dimensional representation $V$ of $G$, 
let $I_\ell(V)$ be the representation of $\hat{\CP}$ such that $G[[t]]$ acts via the evaluation map $\ev_0: G[[t]]\to G$ given by  evaluating $t=0$, and the center $\bbC^\times$ acts by $t\mapsto t^\ell$. Let $\CL_\ell(V)$ be the induced $\hat{G}$-equivariant vector bundle on $\Gr_G$, i.e. 
$  \CL_\ell(V):=\hat{G}\times_{\hat{\CP}} I_\ell (V)^* $,
where $ I_\ell (V)^*$ is the dual of $I_\ell(V)$ as the representation of $\hat{\CP}$. 


The diagram automorphism $\sigma$  on $G$ induces an automorphism on $\hat{G}$  and it preserves $\hat{\CP}$. For any $\lambda\in (P^+)^\sigma$, the vector bundle  $\CL_\ell(V_\lambda)$ is naturally equipped with a $\sigma$-equivariant structure, since 
\[ \hat{G}\rtimes \la \sigma \ra\times _{\hat{\CP}\rtimes \la \sigma \ra}   I_\ell(V_\lambda)^* \simeq  \hat{G}\times_{\hat{\CP}} I_\ell(V_\lambda)^*. \]
Similarly, for any $r$-th root of unity $\xi$, we have the $\hat{G}\rtimes \la \sigma \ra$-equivariant vector bundle $\CL_\ell( V_{\lambda, \xi})$ on $\Gr_G$. 

The following lemma is well-known.
\bl
\label{Well-known}
Let $H_1$ be a linear algebraic group and $H_2$ be a subgroup of $H_1$. Let  $V_1$ be a finite dimensional representation of $H_1$ and let $V_2$ be  a finite dimensional representation of $H_2$. Then we have an isomorphism of $H_1$-equivariant vector bundles
\[ H_1\times_{H_2} (V_2\otimes  V_1|_{H_2}   )\simeq  (H_1\times_{H_2} V_2)\otimes V_1, \]
 given by 
$ (h_1, v_2\otimes v_1)\mapsto   (h_1,v_2)\otimes h_1\cdot v_1$, 
 where $h_1\in H_1$, $v_1\in V_1$ and $v_2\in V_2$.
\el

\bl
\label{gr_fl}
Let $\lambda$ be a $\sigma$-invariant dominant weight of $G$, and let $V$ be a finite dimensional representation of $G\rtimes \la  \sigma  \ra$. 
There is an isomorphism of  $\hat{G}\rtimes  \la \sigma \ra$-representations
\[  H^k(\Gr_G, \CL_\ell(V_{\lambda,\xi}\otimes V ) )\simeq  H^k(\Fl_G, \CL_\ell(  \bbC_{\lambda,\xi} \otimes  V|_{B\rtimes \la \sigma \ra}),\]
for any $k\geq 0$ and $\xi$ an $r$-th root of unity.
\el

\bpf
We have the following isomorphisms of $\hat{G}\rtimes \la \sigma \ra$-equivariant vector bundles
\begin{align*}
\CL_\ell(  \bbC_{\lambda,\xi} \otimes  V|_{B\rtimes \la \sigma \ra})
 &\simeq  \hat{G}\rtimes\la \sigma \ra \times_{\hat{\CI} \rtimes \la \sigma \ra } (\bbC_{\lambda,\xi}\otimes  V|_{B\rtimes \la \sigma \ra })^*\\
 &\simeq    \hat{G} \rtimes\la \sigma \ra    \times_{\hat{\CP}   \rtimes \la \sigma \ra }  (\hat{\CP}\rtimes \la \sigma \ra\times_{\hat{\CI}\rtimes \la \sigma \ra} I_\ell(\bbC_{\lambda,\xi}\otimes  V|_{B\rtimes \la \sigma \ra} ))^*\\
 &\simeq      \hat{G} \rtimes\la \sigma \ra   \times_{\hat{\CP}\rtimes \la \sigma \ra }  ( (\hat{\CP}\rtimes \la \sigma \ra\times_{\hat{\CI}\rtimes \la \sigma \ra} I_\ell(\bbC_{\lambda,\xi}))\otimes  V )^*,
\end{align*}
where the last isomorphism follows from Lemma \ref{Well-known}.

 It is a $\hat{G}\rtimes \la\sigma \ra$-equivariant vector bundle on $\Fl_G$. By Borel-Weil-Bott theorem for finite type algebraic group, we have 
 \[ R^k\pi_*  \CL_\ell(  \bbC_{\lambda,\xi} \otimes  V|_{B\rtimes \la \sigma \ra}) \simeq     
 \begin{cases}
 0  \quad   \quad k>0\\
\CL_\ell (V_{\lambda,\xi}\otimes  V)  \quad  k=0
 \end{cases},
 \]
 where $R^k\pi_*$ is the right derived functor of $\pi_*$.
By  Leray's spectral sequence, the lemma follows.
\epf

Let $\Wmin$ denote the set of  the  minimal  representatives of the left cosets of $W$ in $W_{\ell+\check{h}}$, then for any $w_1 \in W$ and $w_2\in \Wmin$, we have
$\ell(w_1w_2)=\ell(w_1)+\ell(w_2) $. 
Moreover for any $w\in W_{\ell+\check{h}}$ and $\lambda\in P_\ell$ 
\begin{equation}
\label{Affine_Weyl_Dom}
 w\star \lambda \in P^+  \text{ if and only if }  w \in \Wmin, 
 \end{equation}
see \cite[Remark 1.3]{Ko}.
Since $P_\ell$ is the set of integral points in the fundamental alcove of the affine Weyl group $W_{\ell+\check{h}}$, for any dominant weight $\lambda\in P^+$, there exists a unique $w\in \Wmin$ such that $w^{-1}\star \lambda\in P_\ell$. By Lemma \ref{fixed_point}, for any $\sigma$-invariant dominant weight $\lambda\in P^+$, there exists a unique  $w\in (\Wmin)^\sigma$ such that $w^{-1}\star \lambda\in  P_\ell^\sigma$.

Recall that we defined in  Section \ref{rep_lie} the representation  $V_{\lambda, \xi}$ of $\fg\rtimes \la\sigma \ra$ as the representation $V_\lambda$ of $\fg$ together with an operator  $\sigma$ such that $\sigma$ acts on the highest weight vectors by $\xi$, where $\lambda\in (P^+)^\sigma$ and $\xi$ is an $r$-th root of unity. Similarly, the representation $\CH_{\lambda,\xi}$ is the representation $\CH_\lambda$ of $\hat{\fg}\rtimes \la \sigma\ra$ of level $\ell$ together with an operator  $\sigma$ such that $\sigma$ acts on the highest weight vectors by $\xi$.
We  have the following theorem
\bt
\label{sign_Gr}
For any $w\in \Wmin$ such that $\sigma(w)=w$ and for any $\lambda\in P_\ell^\sigma$, 
we have the following isomorphism of  representations of  $\hat{G}\rtimes \la \sigma \ra$,
\[  H^{\ell(w)} (\Gr_G, \CL_\ell( V_{w\star \lambda,\xi} ) )\simeq     H^0(\Gr_G, \CL_\ell( V_{\lambda, \epsilon_w \xi}  )   )  . \] 
\et
\bpf
This follows from Theorem \ref{sign_thm1_equiv} and Lemma \ref{gr_fl}.
\epf

\bco 
\label{sign_cor}
With the same assumption as in Theorem \ref{sign_Gr}.
\begin{enumerate}
\item There exists an isomorphsm of representations of $\hat{\fg}\rtimes \la \sigma \ra$
\[H^{\ell(w)}(\Gr_G, \CL_\ell(V_{\lambda} ))^\vee \simeq  \CH_{\lambda,\epsilon_w }. \]
\item  There exists an isomorphism of representations of $\fg\rtimes \la \sigma \ra$
\[  (H^{\ell(w)}(\Gr_G, \CL_\ell(V_{w\star \lambda}))^\vee)_{\hat{\fg}^-}\simeq V_{\lambda, \epsilon_w}, \]
 where $\hat{\fg}^-=t^{-1}\fg[t^{-1}]$. 
\end{enumerate}
\eco
\begin{proof}
This proposition follows from Theorem \ref{sign_Gr}, combining with Corollary \ref{affine_rep_sign}, Lemma \ref{gr_fl} and Lemma \ref{affine_finite}.
\end{proof}

\subsection{ Affine analogues of BBG resolution and Kostant homology}
\label{affine_BBG_Kostant}

We first recall the construction of BGG resolution in the setting of affine Lie algebra, we refer  the reader to \cite[Section 9.1]{Ku1} for more details, in particular Theorem 9.1.3 therein. There exists a Koszul resolution of the trivial representation $\bbC$ of $\hat{\fg}$, 
\[ \cdots   \to   X_p   \xrightarrow{\delta_{p}}  \cdots   \xrightarrow{\delta_1}  X_0  \xrightarrow{\delta_0}  \mathbb{C}, \]
where $  X_p= U(\hat{\fg})\otimes_{U(\hat{\fp})} \lwedge^p (\hat{\fg}/  \hat{\fp}   ) $. 
From the construction of Koszul resolution, this complex is $\hat{\fg}\rtimes \la \sigma \ra $-equivariant. Given a $\sigma$-invariant dominant weight $\lambda\in P_\ell$.
Set $X_{\lambda, p}:=  U(\hat{\fg})\otimes_{U(\hat{\fp})} ( \lwedge^p (\hat{\fg}/  \hat{\fp}   )  \otimes  \CH_\lambda  ) $.
The complex $X_{\lambda, \bullet}$ is a resolution of $\CH_\lambda$. 
Set
\be
\label{F_lambda}
 F_{\lambda,p}:=   \bigoplus_{w\in \Wmin, \ell(w)=p} \hat{M}(V_{w\star \lambda}),
 \ee
 where $\hat{M}(V_{w\star \lambda})$ is the generalized Verma module introduced in Section \ref{rep_lie}. 
In fact $F_{\lambda,\bullet}$ is a $\sigma$-stable subcomplex of $X_{\lambda,\bullet}$, and moreover $X_{\lambda,\bullet}$ is quasi-isomorphic to $F_{\lambda, \bullet}$. 
Hence $F_{\lambda,\bullet}$ is a resolution of $\CH_\lambda$.


The proof of the following proposition heavily replies on the work of Naito \cite{N2}. 
\bp
\label{BGG_resolution}
Assume that $\sigma(\lambda)=\lambda$.
 Then the complex $F_{\lambda,\bullet}$ is a resolution of $\CH_\lambda$ as representations of $\hat{\fg}\rtimes \la \sigma \ra$, where $\sigma$ maps $\hat{M}(V_{w\star \lambda})$ to $\hat{M}(V_{\sigma(w)\star \lambda})$. In particular when $\sigma(w)=w$, $\sigma$ acts on the highest weight vectors of $\hat{M}(V_{w\star \lambda})$ by the scalar $\epsilon_w$, where $\epsilon_w=(-1)^{\ell(w)-\ell_\sigma(w)} $ as defined in (\ref{def_epsilon_w}).
 \ep
 
 \bpf
  First of all, we note that $\sigma$ maps $\hat{M}(V_{w\star \lambda})$ to $\hat{M}(V_{  \sigma(w)\star \lambda })$ for any $w\in W^\dagger_{\ell+\hat{h}}$, since $\sigma(\rho)=\rho$. In particular if $\sigma(w)=w$, $\sigma$ keeps $\hat{M}(w\star \lambda)$ stable. We need to determine the action of $\sigma$ at the highest weight vector $m_{w\star \lambda}$ of $\hat{M}(w\star \lambda)$. It is easy to see that $\sigma$ acts on $m_{w\star \lambda}$ by a scalar $\epsilon'_w$. In the following we will show that $\epsilon'_w=\epsilon_w$.
 
 Recall that $\hat{\fg}^-$ denote the nilpotent Lie algebra $t^{-1}\fg[t^{-1}]$. 
 It is standard that $\hat{M}(V_{w\star \lambda})$ is a free $U(\hat{\fg}^-)$-module, for each $w\in  \Wmin$. Thus, 
 the resolution $F_{\lambda, \bullet}$ can be used to compute the $\hat{\fg}^-$-homologies of $\CH_\lambda$, in other words, 
 \begin{equation}
 \label{Kostant_Koszul}
 H_p( \hat{\fg}^-, \CH_\lambda  )\simeq  H_p ( (F_{\lambda,\bullet} )_{\hat{\fg}^- } ), \end{equation}
 where the LHS is the $p$-th $ \hat{\fg}^-$-homology of $\CH_\lambda$, and the RHS is the $p$-th homology of the complex $(F_{\lambda,\bullet} )_{\hat{\fg}^- }$ obtained from taking $\hat{\fg}^- $-coinvariants on the complex $F_{\lambda,\bullet}$. Moreover, the isomorphism (\ref{Kostant_Koszul}) is $\fg\rtimes \la \sigma \ra$-equivariant. As a consequence,
 we get the following isomorphism of $\fg\rtimes \la \sigma \ra$-representations,
 \begin{equation}
 \label{Kost_hom1}
   H_p(\hat{\fg}^-, \CH_\lambda  )\simeq   \bigoplus_{w\in  W^\dagger_{\ell+\check{h}}, \ell(w)=p} V_{w\star \lambda} \end{equation}
for each $p\geq 0$, since $(\hat{M}(V_{w\star \lambda}))_{\hat{\fg}^-  }\simeq V_{w\star \lambda}$ as representations of $\fg$ (cf.\,Lemma \ref{affine_finite}). As mentioned above, $\sigma$ acts on $m_{w\star \lambda}\in \hat{M }( V_{w\star \lambda})$ by the scalar $\epsilon'_w$ if $\sigma(w)=w$. It follows that $\sigma$ acts on the highest weight vector $v_{w\star \lambda}$ of $V_{w\star \lambda}$ by $\epsilon'_w$ if $\sigma(w)=w$.

   Let $\fn^-$ be the  nilpotent radical of the negative Borel subalgebra $\fb^-$ of $\fg$. Put 
  \[ \hat{\fn}^-:=\hat{\fg}^-\oplus \fn^-. \] 
 Note that $ \hat{\fn}^-$ is  the nilpotent radical of the opposite affine Borel subalgebra $\hat{\fb}^-:=\hat{\fg}^-\oplus \fb^-$ of $\hat{\fg}$, and  $\hat{\fn}^-$ is $\sigma$-stable. Since $\hat{\fg}^-$ is an ideal in the Lie algebra $\hat{\fn }^-$, we have the following spectral sequence which is compatible with the actions of $\sigma$, 
 \begin{equation}
 \label{Kostant_hom_spec_sequ}
  H_i( \fn^-, H_j(\hat{\fg}^-, \CH_\lambda  ) )\Rightarrow  H_{i+j}(\hat{\fn}^-, \CH_\lambda  ). 
  \end{equation}
 Meanwhile, $H_i( \fn^-, H_j(\hat{\fg}^-, \CH_\lambda  ) )$ and $ H_{i+j}(\hat{\fn}^-, \CH_\lambda  )$ both carry the actions of the Cartan subalgebra $\fh\subset \fb^-$. 
 In fact the spectral sequence (\ref{Kostant_hom_spec_sequ}) degenerates at $E_2$, since we have the following sequence of isomorphisms of $\fh$-modules:
 \begin{align*}
 \bigoplus_{i+j=p} H_i( \fn^-, H_j(\hat{\fg}^-, \CH_\lambda  ) ) & \simeq \bigoplus_{i+j=p}   \bigoplus_{w\in  W^\dagger_{\ell+\check{h}}, \ell(w)=j} H_i(\fn^-, V_{w\star \lambda}) \\
       &\simeq \bigoplus_{i+j=p}  \bigoplus_{w\in  W^\dagger_{\ell+\check{h}}, \ell(w)=j}   \bigoplus_{y\in W, \ell(y)=i}   \mathbb{C}_{y\star(w\star \lambda)} \\
      &\simeq  \bigoplus_{w\in W_{\ell+\check{h}  }, \ell(w) =p} \mathbb{C}_{w\star \lambda}\\      
         &\simeq H_{p}(\hat{\fn}^-, \CH_\lambda ),
      \end{align*}      
where the first isomorphism follows from (\ref{Kostant_Koszul}), the second isomorphism follows from Kostant homology formula for $\fn^-$, the last isomorphism follows from the affine version of Kostant homology formula for $\hat{\fn}^-$ (cf.\,\cite{GL}), and 
the third isomorphism follows since $\Wmin$ is the set of minimal representatives of the left cosets of $W$ in $W_{\ell+\check{h}}$. The set $W^\dagger_{\ell+\check{h}}$ satisfies the following property: for any $u\in W_{\ell+\check{h}}$, there exist unique $w\in W^\dagger_{\ell+\check{h}}$ and $y\in W$ such that $u=yw$ and $\ell(u)=\ell(y)+\ell(w)$. 
 
 We now make a digression on twining characters. Let $V$ be a finite dimensional $\fh\rtimes \la \sigma \ra$-representation such that $\fh$ acts on $V$ semi-simply. Define 
 \[{\rm ch}_\sigma(V):=\sum_{\mu\in \fh^*, \sigma(\mu) =\mu} {\rm tr} (\sigma | V(\mu)) e^{\mu}, \]
 where $V(\mu)$ denotes the $\mu$-weight space in $V$. Then 
 \begin{align} 
 \sum_{i+j=p} {\rm ch}_\sigma ( H_i( \fn^-, H_j(\hat{\fg}^-, \CH_\lambda  ) ) ) 
 &=\sum_{i+j=p}    \sum_{w\in (\Wmin)^\sigma,\ell(w)=j } \epsilon'_w  {\rm ch}_\sigma (H_i(\fn^-, V_{w\star \lambda}) )  \\
 &=\sum_{i+j=p}  \sum_{w\in (\Wmin)^\sigma, \ell(w)=j } \epsilon'_w \sum_{y\in W^\sigma,\ell(y)=i } c_y(\sigma,V_{w\star \lambda}) e^{y\star (w\star \lambda) }\\
 &= \sum_{i+j=p}  \sum_{ \substack{w\in (\Wmin)^\sigma, \ell(w)=j \\ y\in W^\sigma,\ell(y)=i  }}   \epsilon'_w c_y(\sigma,V_{w\star \lambda}) e^{(yw)\star \lambda},
 \end{align}
 where  $c_y(\sigma,V_{w\star \lambda}):= {\rm tr}(\sigma| H_i(\fn^-, V_{w\star \lambda})_{(yw)\star \lambda} )$. Here $H_i(\fn^-, V_{w\star \lambda})_{ (yw)\star \lambda}$ denotes the $(yw)\star \lambda$-weight space in $H_i(\fn^-, V_{w\star \lambda})$. In the above sequence of equalities, 
 the first equality follows from  (\ref{Kost_hom1}) and the discussions after that, the second equality follows from \cite[Prop.3.2.1]{N2} for the Kostant homology of $\fn^-$. By \cite[Prop.3.2.1]{N2} for the Kostant homology of $\hat{\fn}^-$, we  have 
 \begin{equation}
 \label{Naito_2}
  {\rm ch}_\sigma(H_p(\hat{\fn}^-, \CH_\lambda ))=\sum_{u\in (W_{\ell+\check{h}})^\sigma, \ell(u)=p}  c_{u}(\sigma,\CH_\lambda) e^{u\star \lambda}, \end{equation}
 where  $ c_{u}(\sigma,\CH_\lambda):= {\rm tr}(\sigma| H_p(\hat{\fn}^-, \CH_\lambda )_{u\star \lambda} )$. Here $H_p(\hat{\fn}^-, \CH_\lambda )_{u\star \lambda}$ denotes the $u\star \lambda$-weight space in $H_p(\hat{\fn}^-, \CH_\lambda )$. Since the spectral sequence (\ref{Kostant_hom_spec_sequ}) degenerates at $E_2$, we have
 \begin{equation}
 \label{Naito_3}
   \sum_{i+j=p} {\rm ch}_\sigma ( H_i( \fn^-, H_j(\hat{\fg}^-, \CH_\lambda  ) ) ) ={\rm ch}_\sigma(H_p(\hat{\fn}^-, \CH_\lambda )).\end{equation}
 Comparing formulae (39) and (\ref{Naito_2}) via (\ref{Naito_3}), we see that for any $w\in (\Wmin)^\sigma$, $c_w(\sigma,\CH_\lambda)=\epsilon'_w c_e(\sigma,V_{w\star \lambda})$, where $e$ is the identity element in the Weyl group $W$. Clearly $c_e(\sigma,V_{w\star \lambda})=1$, hence $c_w(\sigma,\CH_\lambda)=\epsilon'_w $.
 We can read further from \cite[Corollary 3.2.3]{N2}, in fact $c_w(\sigma,\CH_\lambda)=\epsilon_w$. Hence $\epsilon'_w=\epsilon_w$. Thus, this finishes the proof.

 \epf

 For any finite dimensional representation  $V$ of $\fg$ and for any $z\in \bbC^\times$, we denote by  $V^z$ the
representation of $\hat{\fg}^-$ that is obtained by evaluating $t$ at $z$. Let $H_i(\hat{\fg}^-, \CH_\lambda\otimes  V_\mu^1)$ be the $i$-th  $\hat{\fg}^-$-homology on $ \CH_\lambda\otimes  V_\mu^z$ where  $\hat{\fg}^-$ acts on $\CH_\lambda\otimes  V_\mu^z$ diagonally. The following theorem will be used in the proof of Theorem \ref{new_old}.
\bt
\label{BGG_Kostant}
For any $\lambda\in P^\sigma_\ell$ and  $\mu\in (P^+)^\sigma$, 
The $\hat{\fg}^-$-homology groups $H_*(\hat{\fg}^-, \CH_\lambda\otimes  V_\mu^1)$ can be computed by the cohomology groups of a complex of $\fg\rtimes \la  \sigma  \ra$-representations,
   \[ \cdots   \to   D_p  \xrightarrow{\delta^p}  \cdots   D_1   \xrightarrow{\delta^1} D_0  \xrightarrow{\delta^0}   0, \]
where as representations of $\fg$,
$ D_p=\bigoplus_{w\in \Wmin, \, \ell(w)=p}   V_{w\star \lambda}\otimes V_\mu$, 
and $\sigma$ maps $V_{w\star \lambda}\otimes V_\mu $ to $V_{\sigma(w)\star \lambda}\otimes V_\mu $. In particular if $\sigma(w)=w$, then $\sigma$ acts on the highest weight vectors of  $V_{w\star \lambda}$
  by $\epsilon_w=(-1)^{\ell(w)-\ell_\sigma(w)}$.
\et

\bpf
 From the resolution  $F_{\lambda, \bullet}\to \CH_\lambda$, by tensoring with $V_\mu^1$ we get a resolution of  $\CH_\lambda \otimes V_\mu^1$ as representations of $\hat{\fg}\rtimes \la \sigma \ra$ 
\[ \cdots   \to   F_{\lambda,p}\otimes V_\mu^1  \xrightarrow{\delta^p}  \cdots   F_{\lambda,1}\otimes V_\mu^1    \xrightarrow{\delta^1} F_{\lambda,0}\otimes V_\mu^1   \xrightarrow{\delta^0}   0, \]
As $\fg$-modules, we have 
\[ (\hat{M}(V_{w\star \lambda }  ) \otimes V_\mu^1)_{\hat{\fg}^-  }\simeq ( V_{w\star \lambda }\otimes_\bbC U(\hat{\fg}^-))\otimes_{U(\hat{\fg}^-)} V_\mu^1\simeq  V_{w\star \lambda }\otimes V_\mu. \]
Hence the complex 
\[\cdots   \to   (F_{\lambda,p}\otimes V_\mu^1)_{\hat{\fg}^-  }  \xrightarrow{\delta^p}  \cdots   (F_{\lambda,1}\otimes V_\mu^1)_{\hat{\fg}^-  }     \xrightarrow{\delta^1}  (F_{\lambda,0}\otimes V_\mu^1)_{\hat{\fg}^-  }     \xrightarrow{\delta^0}  0 \]
is quasi-isomorphic to 
\[ \cdots   \to   D_p  \xrightarrow{\delta^p}  \cdots   D_1   \xrightarrow{\delta^1} D_0  \xrightarrow{\delta^0}   0.\]
By Proposition \ref{BGG_resolution}, $\sigma$ maps $V_{w\star \lambda}\otimes V_\mu $ to $V_{\sigma(w)\star \lambda}\otimes V_\mu $. In particular if $\sigma(w)=w$, then $\sigma$ acts on the highest weight vectors of  $V_{w\star \lambda}$
  by $\epsilon_w=(-1)^{\ell(w)-\ell_\sigma(w)}$.
\epf

\section{$\sigma$-twisted representation ring and fusion ring}
\label{Sect5}
\subsection{$\sigma$-twisted representation ring}
\label{twist_ring_sect}
Let $V$ be a finite dimensional representation of $\fg$. For any irreducible representation $V_\lambda$ of $\fg$ of highest weight $\lambda$, we denote by $\Hom_{\fg}(V_\lambda, V)$ the multiplicity space of $V_\lambda$ in $V$.
In particular we have the following natural decomposition 
\[V=\bigoplus_{\lambda\in P^+}  \Hom_{\fg}(V_\lambda, V)\otimes V_\lambda. \]

Let $R(\fg,\sigma)$ be the free abelian group with the symbols $[V_\lambda]_\sigma$ as a basis, where $\lambda\in (P^+)^\sigma$. 
Given any finite dimensional representation $V$ of $\fg\rtimes \la \sigma  \ra$, $V$ can be decomposed as follows
\[  V=\bigoplus_{\lambda\in (P^+)^\sigma}   \Hom_{\fg}(V_\lambda, V )\otimes  V_\lambda    \oplus   \bigoplus_{\lambda\not \in (P^+)^\sigma}  \Hom_{\fg}(V_\lambda, V)\otimes V_\lambda, \] 
as a representation of $\fg$. Put
\[  [V]_\sigma:=\sum_{\lambda\in (P^+)^\sigma}  \tr(\sigma| \Hom_{\fg}(V_\lambda,V)) [V_\lambda]_\sigma \in  R(\fg,\sigma).\]

Let $X$ be a finite dimensional representation of the cyclic group $\la \sigma \ra$, and for any representation $V$ of $\fg\rtimes \la \sigma \ra$, $X\otimes V$ is naturally a representation of $\fg\rtimes \la \sigma  \ra$, which is defined as follows
\[  (u, \sigma^i)\cdot x\otimes v=\sigma^i\cdot x \otimes  (u, \sigma^i)\cdot v, \]
where $u\in \fg, x\in X,v\in V$ and $i\in \bbZ$.
Similarly, $V\otimes X$ is also naturally a representation of $\fg\rtimes \la \sigma \ra$. 
The following lemma is obvious.
\bl
\label{trace_lem}
We have  $[X\otimes V]_\sigma= \tr(\sigma| X)  [V]_\sigma$, and  $    [V\otimes X]_\sigma=\tr(\sigma|X)[V]_\sigma$. 
\el

We define a multiplication $\otimes$ on  $R(\mathfrak{g},\sigma)$,
$  [V_\lambda]_\sigma \otimes [V_\mu]_\sigma:=[V_\lambda\otimes V_\mu]_\sigma $, 
 for any $\lambda,\mu\in (P^+)^\sigma$. By definition, we have
\[ [V_\lambda\otimes V_\mu]_\sigma=\sum_{\sigma(\nu)=\nu}  \tr(\sigma|  \Hom_{\fg}(V_\nu, V_\lambda\otimes V_\mu)    )  [V_\nu]_\sigma. \]

\bp
$R(\fg,\sigma)$ is a commutative ring with $[V_0]_\sigma$ as the unit. 
\ep
\bpf
The commutativity is clear. 
We first show that the product $\otimes$ on $R(\fg,\sigma)$ satisfies the associativity, i.e. for any $\lambda, \mu, \nu\in (P^+)^\sigma$,
\[ ( [V_\lambda]_\sigma\otimes [V_\mu]_\sigma)\otimes [V_\nu]_\sigma =[V_\lambda]_\sigma \otimes (  [V_\mu]_\sigma\otimes [V_\nu]_\sigma ). \]
It suffices to show that for any $\lambda\in (P^+)^\sigma$ and any representation $V$ of $\fg\rtimes \la \sigma \ra$,
\[ [V_\lambda]_\sigma\otimes [V]_\sigma=[V_\lambda\otimes V]_\sigma, \text{ and } [V]_\sigma\otimes [V_\lambda]_\sigma=[V\otimes V_\lambda]_\sigma. \]
We have the following equalities
\begin{align*}
[V_\lambda]_\sigma\otimes [V]_\sigma
&=\sum_{\sigma(\mu)=\mu } \tr(\sigma| \Hom_{\fg}(V_\mu,V) ) ( [V_\lambda]_\sigma\otimes [V_\mu]_\sigma )   \\
&=\sum_{\sigma(\mu)=\mu } \tr(\sigma| \Hom_{\fg}(V_\mu,V) ) [V_\lambda\otimes V_\mu]_\sigma\\
&=\sum_{\sigma(\mu)=\mu}  [V_\lambda\otimes V_\mu\otimes  \Hom_{\fg}(V_\mu,V) ]_\sigma\\
&=[ \bigoplus_{\mu}  V_\lambda \otimes  V_\mu\otimes \Hom_{\fg}(V_\mu,V)      ]_\sigma=[V_\lambda\otimes V]_\sigma,
\end{align*}
where  the third equality follows from  Lemma \ref{trace_lem}, and others follows from definition of the multiplication $\otimes$. 
The equality $[V]_\sigma\otimes [V_\lambda]_\sigma=[V\otimes V_\lambda]_\sigma$ can be proved similarly. 
In the end $[V_0]_\sigma$ is the unit since for any $\lambda\in (P^+)^\sigma$, 
\[ [V_\lambda]_\sigma\otimes [V_0]_\sigma=[V_0]_\sigma\otimes [V_\lambda]_\sigma=[V_\lambda\otimes V_0]_\sigma=[V_\lambda]_\sigma.\]
\epf

Recall that  $W_\lambda$ denotes the representation of $\fg_\sigma$ of highest weight $\iota(\lambda)$, and $W_\lambda(\mu)$ is the $\iota(\mu)$-weight space of $W_\lambda$, where $\iota$ is defined in Section \ref{root_system}. The following theorem is due to Jantzen \cite{Ja}.
\bt
\label{Jantzen}
Let $\lambda\in (P^+)^\sigma$ and $\mu\in P^\sigma$. We have 
$ \tr(\sigma| V_\lambda(\mu)  )=\dim W_{\lambda}(\mu)$.
\et
For any finite dimensional representation $V$ of $\fg\rtimes \la \sigma\ra$,
we define the $\sigma$-twisted character $\Ch_\sigma(V)$ of $V$ as follows
\[\Ch_\sigma(V):=\sum_{\mu\in P^\sigma} \tr(\sigma|V(\mu)) e^\mu,\]
where $V(\mu)$ denotes the $\mu$-weight space of $V$. The following lemma is obvious. 
\bl
\label{sigma_tensor}
For any two finite dimensional $\fg\rtimes \la\sigma \ra$-representations $V,V'$, we have 
\[ \Ch_\sigma(V\otimes V')=\Ch_\sigma(V)\Ch_\sigma(V'). \]
\el

  \bl
 \label{tensor_inv_multiplicity}
 Let $\vec{\lambda}$ be a tuple of $\sigma$-invariant dominant weights of $\fg$ and let $\nu$ be another $\sigma$-invariant dominant weight of $\fg$. The following equality holds
 \[  \tr(\sigma|  \Hom_{\fg}(V_\nu, V_{\vec{\lambda}}  )) =\tr(\sigma|  (V_{\vec{\lambda}} \otimes V_{\nu^*} )^\fg  ). \]
 \el
 \bpf
Let $w_0$ be the longest element in the Weyl group $W$ of $\fg$. There exists a representative $\bar{w}_0$ of $w_0$ in  $G$ such that $\sigma(\bar{w}_0)=\bar{w}_0$ (see \cite[Section 2.3]{HS}). Hence $\sigma(\overline{w}_0\cdot v_\nu)=\overline{w}_0\cdot v_\nu$, where $v_\nu\in V_{\nu}$ is the  highest weight vector. The vector $\overline{w}_0\cdot v_\nu$ is of the lowest weight $w_0(\nu)$. 
Let $V_{\nu}^*$ denote the  dual  representation of  $V_\nu$, and let $\sigma^{*}$ be the action on $V^*_\nu$ induced by the action $\sigma$ on $V_\nu$. Then $\sigma^*$ keeps the highest weight vectors in $V_\nu^*$ invariant. 
 
As representations of $\fg$, there is an isomorphism $V^*_\nu\simeq V_{-w_0(\nu)}=V_{\nu^*}$ which is unique up to a scalar. It intertwines the action of $\sigma^*$ on $V_\nu^*$ and the action of $\sigma$ on $V_{\nu^*}$. 
Note that there is a natural isomorphism $\Hom_\fg(V_\nu, V_{ \vec{\lambda}}  ) \simeq (V_{\vec{\lambda}} \otimes V^*_\nu)^\fg$, which is $\sigma$-equivariant. This concludes the proof.
 \epf

The following theorem was proved in \cite{HS}. We give a simple proof here using Jantzen formula directly.
\bt[\cite{HS}]
\label{HS}
Let $\vec{\lambda}$ be a tuple of dominant weights of $\fg$. We have 
$ \tr(\sigma|   V^\fg_{\vec{\lambda}}) = \dim W_{\vec{\lambda}} ^{\fg_\sigma}.$
\et
\bpf
On one hand, from the decomposition
$ V_{\vec{\lambda}}\simeq  \bigoplus_{\nu\in P^+}  \Hom_{\fg}(V_\mu, V_{\vec{\lambda}})\otimes  V_\mu $, 
we have
\[  \Ch_\sigma(V_{\vec{\lambda}})=\sum_{\mu\in (P^+)^\sigma}  \tr(\sigma| \Hom_{\fg}(V_\mu, V_{\vec{\lambda}})) \Ch_\sigma(V_\mu). \]
On the other hand, we have the following equalities
\begin{align*} 
 \Ch_\sigma(V_{\vec{\lambda}})&=\Ch_\sigma(V_{\lambda_1})\cdots  \Ch_\sigma(V_{\lambda_k})= \Ch(W_{\lambda_1})\cdots  \Ch(W_{\lambda_k})\\
 &=\Ch(W_{\vec{\lambda}})=\sum\dim \Hom_{\fg_\sigma}(W_\mu, W_{\vec{\lambda}}) \Ch(W_\mu),
\end{align*}
where the first equality follows from Lemma \ref{sigma_tensor} and  the second equality follows from Theorem \ref{Jantzen}. In view of Lemma \ref{tensor_inv_multiplicity}, the theorem follows.
\epf

 

Let $R(\fg_\sigma)$ denote the representation ring of $\fg_\sigma$. 
\bp
\label{iso_twist}
There is a natural ring isomorphism $ R(\fg,\sigma)\simeq R(\fg_\sigma)$
by sending $[V_\lambda]_\sigma\mapsto [W_{\lambda   }]$ for any $\sigma$-invariant dominant weight $\lambda$.
\ep
\bpf

For any $\lambda,\mu\in  (P^+)^\sigma$, consider the following two decompositions
\[ [V_\lambda]_\sigma\otimes [V_\mu]_\sigma=\sum_{\sigma(\mu)=\mu} \tr(\sigma| \Hom_{\fg}(V_\nu,V_\lambda\otimes V_\mu))   [V_\nu]_\sigma, \]

\[ [W_\lambda]\otimes [W_\mu ] =\sum_{\sigma(\mu)=\mu}  \dim \Hom_{\fg}(W_{\nu},W_\lambda \otimes W_\mu    )    [W_\nu ]. \]

In view of Theorem \ref{HS} and Lemma \ref{tensor_inv_multiplicity}, we have
\[  \tr(\sigma| \Hom_{\fg}(V_\nu, V_\lambda\otimes V_\mu))=\dim \Hom_{\fg_\sigma}(W_\nu, W_\lambda  \otimes W_\mu    ). \]
Hence the proposition follows. 
\epf

\subsection{A new definition of $\sigma$-twisted fusion ring via Borel-Weil-Bott theory}
\label{new_def_fusion_sect}
\bl
\label{Euler_Poincare}
The operation $[\cdot]_\sigma$ satisfies Euler-Poincar\'e property, i.e. 
for any complex of  finite dimensional $\fg\rtimes \la \sigma \ra$-representations
\[V^\bullet:= \cdots   \xrightarrow{d_{i-1}}   V^i  \xrightarrow{d_i}     V^{i+1}   \xrightarrow{d_{i+1}}   \cdots  \]
such that only finite many $V^i$ are nonzero, 
we have 
\[ [V^\bullet]_\sigma=\sum_{i} (-1)^i [H^i(V^\bullet)]_\sigma, \]
where $[V^\bullet]_\sigma:=\sum_{i}  (-1)^i  [V^i]_\sigma$, and $H^i(V^\bullet)$ is the $i$-th cohomology of this complex. 
\el
\bpf
First of all, we have Euler-Poincar\'e property in the representation ring $R(\fg\rtimes \la \sigma \ra)$ of $\fg\rtimes \la \sigma \ra $, i.e.
\[  \sum_{i}  (-1)^i  [V^i]=\sum_{i} (-1)^i  [H^i(V^\bullet)]. \]
Secondly we can define a linear map  $  R(\fg\rtimes \la \sigma \ra)\to   R(\fg, \sigma) $
given by  $[V]\mapsto  [V]_\sigma$. It is well-defined and additive, since any finite dimensional representation of $\fg\rtimes \la \sigma \ra$ is completely reducible. 
Hence the lemma follows. 
\epf

Recall the $\sigma$-twisted fusion ring $R_\ell(\fg,\sigma)$ defined in Section \ref{twisted_fusion_sect}. 
We  embed $R_\ell(\fg,\sigma)$ into $R(\fg,\sigma)$ as free abelian groups by simply sending $\lambda$ to $[V_\lambda]_\sigma$ for any $\lambda\in P_\ell^\sigma$. From now on we view $R_\ell(\fg,\sigma)$ as a free abelian group with basis $\{[V_\lambda]_\sigma\li  \lambda\in P_\ell^\sigma\}$. The fusion product $\lambda\cdot \mu$ in $R_\ell(\fg,\sigma)$ will be written as $[V_\lambda]_\sigma\cdot [V_\mu]_\sigma$. 

Given any integrable representation $\CH$ of $\hat{\fg}$, we denote by $\CH_{\hat{\fg}^-}$ the coinvariant space of $\hat{\fg}^-$ on $\CH$. If $\CH$ is a representation of $\hat{\fg}\rtimes \la \sigma \ra$, then  the space $\CH_{\hat{\fg}^-}$ is naturally a representation of $\fg\rtimes \la \sigma \ra$. For any $\lambda,\mu\in P_\ell^\sigma$, we define
\be
\label{new_def_prod}
[V_\lambda]_\sigma\otimes_\ell [V_\mu]_\sigma:=    [ (H^*(\Gr_G, \CL_\ell(V_\lambda\otimes V_\mu) )^\vee)_{\hat{\fg}^-  }]_\sigma \in R_{\ell}(\fg,\sigma),
\ee
where we view  $(H^*(\Gr_G, \CL_\ell(V_\lambda\otimes V_\mu) )^\vee)_{\hat{\fg}^-  }$  as a  complex of $\fg\rtimes \la \sigma \ra$-representations with zero differentials. 

Note that all representations of $\hat{\fg}$ appearing in 
$H^*(\Gr_G, \CL_\ell(V_\lambda\otimes V_\mu) )^\vee$  are of level $\ell$, and only finite many cohomology groups are nonzero. Hence the above definition makes sense.

Recall the representation $\CH_\nu\otimes V^z_\mu$ defined in Section \ref{affine_BBG_Kostant}. The following is a vanishing theorem of Lie algebra cohomology due to Teleman \cite{Te}.
\bt
\label{Teleman}
For any $\lambda, \mu, \nu\in P_\ell$ and for any $i\geq 1$, $V_\lambda$ does not occur in $H_i(\hat{\fg}^-, \CH_\nu\otimes  V^z_\mu  )$ as a $\fg$-reprepsentation.
\et


We now  show that the product defined in (\ref{new_def_prod}) is exactly the fusion product. 
\bt
\label{new_old}
Two products on $R_\ell(\fg,\sigma)$ coincide, i.e. for any $\lambda,\mu\in P_\ell^\sigma$, we have
$ [V_\lambda]_\sigma \otimes_\ell [V_\mu]_\sigma=[V_\lambda]_\sigma\cdot [V_\mu]_\sigma.$
\et
\bpf
Consider the following decomposition
\[V_\lambda\otimes V_\mu=\bigoplus_{\nu}  \Hom_{\fg}(V_\nu, V_\lambda\otimes V_\mu)  \otimes   V_\nu. \]

By the fact (\ref{Affine_Weyl_Dom}), we may in further write 
  \begin{equation}
  \label{Decom_ten}
  V_\lambda\otimes V_\mu \simeq \bigoplus_{w\in \Wmin, \nu\in P_\ell }  \Hom_{\fg}(V_{w\star \nu}, V_\lambda\otimes V_\mu)     \otimes   V_{w\star \nu}.\end{equation}
We have the following chain of equalities
\begin{align*}
 [V_\lambda]_\sigma \otimes_\ell [V_\mu]_\sigma 
 &=\sum_i (-1)^i [(H^i(\Gr_G, \CL_\ell(V_\lambda\otimes V_\mu))^\vee)_{\hat{\fg}^-} ]_\sigma \\
 &=\sum_i(-1)^i  \sum_{\substack{ w\in (\Wmin)^\sigma \\ \ell(w)=i, \nu\in P_\ell ^\sigma } }   [  \Hom_{\fg}(V_{w\star \nu}, V_\lambda\otimes V_\mu)   \otimes (H^i(\Gr_G, \CL_\ell(V_{w\star \nu}))^\vee)_{\hat{\fg}^-} ]_\sigma \\
  &= \sum_{w\in (\Wmin)^\sigma, \nu\in P_\ell ^\sigma } (-1)^{\ell(w)}   [  \Hom_{\fg}(V_{w\star \nu}, V_\lambda\otimes V_\mu)  \otimes V_{\nu, \epsilon_w}) ]_\sigma \\
 &= \sum_{w\in (\Wmin)^\sigma, \nu\in P_\ell ^\sigma }   (-1)^{\ell_\sigma(w)} \tr(\sigma| \Hom_{\fg}(V_{w\star \nu}, V_\lambda\otimes V_\mu)  ) [V_{\nu} ]_\sigma,
  \end{align*}
where the second equality follows from the decomposition (\ref{Decom_ten}), the third equality follows from Corollary \ref{sign_cor} and the fourth equality follows from Lemma \ref{trace_lem}.
  By Lemma \ref{sigma_tau} and Proposition \ref{Propagation_prop}, we have the following $\sigma$-equivariant isomorphisms:
\begin{align*}
V_{\fg, \ell, \lambda, \mu, \nu^*}(\bbP^1, 0,1,\infty)
&\simeq V_{\fg, \ell, \lambda^*, \mu^*, \nu}( \bbP^1, 0,1,\infty )\\
&\simeq   (\CH_\nu\otimes V^\infty_{\lambda^*} \otimes  V^1_{\mu^*} )_{\fg[t^{-1}]}     \\  
&\simeq \Hom_\fg(V_\lambda, H_0(\hat{\fg}^-, \CH_\nu\otimes V^1_{\mu^*} )).
\end{align*}
The following formula follows immediately from 
Theorem \ref{Teleman}
\[
   \tr(\sigma|  \Hom_{\fg}(V_\lambda, H_0(\hat{\fg}^-, \CH_\nu\otimes V^1_{\mu^*} ) ) )=\sum_{i} (-1)^i\tr(\sigma|  \Hom_{\fg}(V_\lambda, H_i(\hat{\fg}^-, \CH_\nu\otimes V^1_{\mu^*} )) ). \]
By Lemma \ref{Euler_Poincare} and Theorem \ref{BGG_Kostant}, we have
\begin{align*}
& \sum (-1)^i \tr(\sigma|  \Hom_{\fg}(V_\lambda, H_i(\hat{\fg}^-, \CH_\nu\otimes  V^1_\mu  ) ) ) \\
 &= \sum_{i}  (-1)^i\sum_{w\in (\Wmin)^\sigma, \ell(w)=i}   \tr(\sigma| \Hom_{\fg}(V_\lambda, V_{w\star \nu, \epsilon_w } \otimes V_{\mu^*}     ) )\\
 &=\sum_{w\in (\Wmin)^\sigma}  (-1)^{\ell_\sigma(w)}    \tr(\sigma|   \Hom_{\fg}(V_\lambda, V_{w\star \nu } \otimes V_{\mu^*} )   ).
\end{align*}  
  It follows that
\begin{align*}
[V_\lambda]_\sigma \cdot [V_\mu]_\sigma 
&=\sum_{\nu\in P_\ell^\sigma}  \tr(\sigma| V_{\fg, \ell, \lambda, \mu, \nu^*}(\bbP^1, 0,1,\infty)  )  [V_\nu]_\sigma\\
&=\sum_{w\in (\Wmin)^\sigma, \nu\in P^\sigma_\ell}   (-1)^{\ell_\sigma(w)}  \tr(\sigma|  \Hom_{\fg}(V_\lambda, V_{w\star \nu } \otimes V_{\mu^*} ) ) [V_\nu]_\sigma.
\end{align*}
In the end, we need to check that 
\[  \tr(\sigma| \Hom_{\fg}(V_{w\star \nu}, V_\lambda\otimes V_\mu) )=    \tr(\sigma|  \Hom_{\fg}(V_{\lambda}, V_{w\star\nu } \otimes V_{\mu^*})).\]
In view of Lemma \ref{tensor_inv_multiplicity}, it  reduces to show that the trace of $\sigma$ on $V^\fg_{w\star \nu^*, \lambda, \mu}$ and $V^\fg_{\lambda^*, w\star \nu, \mu^*}$ are equal. This is a consequence of Lemma \ref{inv_intertwiner}.
\epf

From the proof of Theorem \ref{new_old}, we get the following twisted analogue of Kac-Walton formula (in the usual setting, see \cite{Ka,Wa}).
\bt
\label{trace_conformal_tensor}
 For any  $\lambda,\mu, \nu\in P_\ell^\sigma$, we have 
\[  \tr(\sigma| V_{\fg, \lambda,\mu, \nu}(\bbP^1, 0,1,\infty) )= \sum_{w\in (\Wmin)^\sigma}  (-1)^{\ell_\sigma(w)} \tr(\sigma |  V^{\fg}_{\lambda,\mu, w\star \nu}). \]
\et

\br
The proofs of Theorem \ref{new_old}, \ref{trace_conformal_tensor} do not rely on the fact that the trace on conformal blocks is a fusion rule. In fact Theorem \ref{trace_conformal_tensor} is used to show that the trace on conformal blocks gives a fusion rule, see Lemma \ref{trace_integer}.
\er
\subsection{Ring homomorphism from $\sigma$-twisted representation ring to $\sigma$-twisted fusion ring}
\label{ring_hom_sect}

 We first construct a $\bbZ$-linear map  $$\pi_\sigma:  R(\mathfrak{g},\sigma)\to  R_\ell(\mathfrak{g}, \sigma).$$
 For any finite dimensional $\fg\rtimes \la \sigma \ra$-representation $V$, we define 
 \[ \pi_\sigma([V]_\sigma):=  [ (H^*(\Gr_G, \CL_\ell(V))^\vee)_{\hat{\fg}^- }  ]_\sigma  \in R_\ell(\fg, \sigma). \]

\bl
\label{lem_change}
For any $w\in (\Wmin)^\sigma$ and $\lambda\in (P^+)^\sigma$, we have
\[[(H^*(\Fl_G, \CL_\ell(w\star \lambda))^\vee)_{\hat{\fg}^- }  ]_\sigma=(-1)^{\ell_\sigma(w)}[(H^*(\Fl_G, \CL_\ell(\lambda))^\vee)_{\hat{\fg}^- }  ]_\sigma. \]
\el
\bpf
We can write $\lambda=y\star \lambda_0$ where $y\in (\Wmin)^\sigma$ and $\lambda_0\in (P_\ell)^\sigma$. Then $w\star \lambda=(wy)\star \lambda_0$.
In view of Theorem \ref{BWB} and Theorem \ref{sign_thm1_equiv}, we have 
\begin{align*}
[(H^*(\Fl_G, \CL_\ell(w\star \lambda))^\vee)_{\hat{\fg}^- }  ]_\sigma
&= (-1)^{\ell_\sigma(wy)}[(H^*(\Fl_G, \CL_\ell( \lambda_0))^\vee)_{\hat{\fg}^- }  ]_\sigma\\
&=(-1)^{\ell_\sigma(w)} [(H^*(\Fl_G, \CL_\ell( \lambda))^\vee)_{\hat{\fg}^- }  ]_\sigma.
\end{align*}
Hence the lemma follows.
\epf

\bp
\label{sign_change}
Given a  finite dimensional representation $V$ of $\fg\rtimes \la \sigma \ra$. 
For any $\lambda\in P_\ell^\sigma$ and $w\in (W^+_{\ell+\check{h}})^\sigma$, the following equality holds in $R_\ell(\fg,\sigma)$
\[   [ (H^*(\Gr_G, \CL_\ell(V_{w\star \lambda}\otimes V)  )^\vee)_{\hat{\fg}^-}]_\sigma =(-1)^{\ell_\sigma(w)}  [ (H^*( \Gr_G, \CL_\ell(V_{ \lambda} \otimes V)  )^\vee)_{\hat{\fg}^-} ]_\sigma  .\]
\ep

\bpf

In view of Lemma \ref{gr_fl}, it suffices to show that 
\[   [ (H^*(\Fl_G, \CL_\ell(\bbC_{w\star \lambda} \otimes V|_{B\rtimes \la \sigma \ra})  )^\vee)_{\hat{\fg}^-  }]_\sigma =(-1)^{\ell_\sigma(w)}  [ (H^*( \Fl_G, \CL_\ell(\bbC_{\lambda} \otimes V|_{B\rtimes \la \sigma \ra})  )^\vee)_{\hat{\fg}^-}   ]_\sigma. \]

Note that there exists a filtration  of $B\rtimes \la \sigma \ra$-representations 
\[ 0=V_0\subset V_1\subset  V_2\subset \cdots   \subset  V_k=V \]
on $V$,
such that for each $i$, 
\[V_{i}/V_{i-1} \simeq  
\begin{cases}
V(\mu)  \quad    \text{ if }  \sigma(\mu)=\mu \\
\bigoplus_{i=0}^{r-1}  V(\sigma^i(\mu))  \quad   \text{ otherwise}
\end{cases}, \]
where $V(\mu)$ denotes the $\mu$-weight space of $V$.
By Lemma \ref{trace_lem}, it is easy to check that 
\begin{align*}
& [(H^*(\Fl_G, \CL_\ell(\bbC_\lambda\otimes  \bigoplus_{i=0}^{r-1} V(\sigma^i(\mu)))^\vee)_{\hat{\fg}^-} ]_\sigma \\
& =\begin{cases}
\tr(\sigma|V(\mu) )[(H^*(\Fl_G, \CL_\ell(\lambda+\mu))^\vee)_{\hat{\fg}^- } ]_\sigma     \quad   \text{ if } \sigma(\mu)=\mu  \\
0  \quad  \quad  \quad  \quad  \text{ otherwise } 
\end{cases}.
 \end{align*}

Hence we get the following isomorphisms
\begin{align*}
 [(H^*(\Fl_G, \CL_\ell(\bbC_\lambda\otimes  V|_{B\rtimes \la\sigma \ra} )   )^\vee)_{\hat{\fg}^-  } ]_\sigma&=\sum_{i} [(H^*(\Fl_G, \CL_\ell(\bbC_\lambda\otimes V_i/V_{i-1}))^\vee)_{\hat{\fg}^-}  ]_\sigma  \\
 &= \sum_{\mu\in P^\sigma}  \tr(\sigma| V(\mu) ) [ (H^*(\Fl_G, \CL_\ell(\lambda+\mu))^\vee)_{\hat{\fg}^-}  ]_\sigma.
 \end{align*}
Similarly, we have 
\[ [(H^*(\Fl_G, \CL_\ell( \bbC_{w\star \lambda}\otimes V|_{B\rtimes \la\sigma \ra} )   )^\vee)_{\hat{\fg}^-  } ]_\sigma\\
 = \sum_{\mu\in P^\sigma}  \tr(\sigma| V(\mu) ) [ (H^*(\Fl_G, \CL_\ell(w\star \lambda+\mu))^\vee)_{\hat{\fg}^-}  ]_\sigma. \]

We can write $w$ as $ w=\tau_{\beta}y^{-1}$,
where $y\in W^\sigma$ and $\tau_\beta$ is the translation for $\beta\in (\ell+\check{h} )Q^\sigma$. It is easy to check that $ w\star \lambda+ \mu=w\star ( \lambda+  y\cdot \mu )$.

Since $V$ is a representation of $\fg\rtimes \la \sigma \ra$, for any $y\in W^\sigma$ we have
$  \tr(\sigma|  V(\mu) )=\tr(\sigma|   V(y\cdot \mu)  ) $, 
 where $V(\mu)$ and $V(y\cdot\mu)$ denote the weight spaces of $V$ as representation of $\fg$.
 We have the following chain of equalities
\begin{align*}
 & [(H^*(\Fl_G, \CL_\ell(\bbC_{w\star \lambda}\otimes  V|_{B\rtimes \la\sigma \ra})   )^\vee)_{\hat{\fg}^-  } ]_\sigma \\
  &=\sum_{\mu\in P^\sigma}  \tr(\sigma| V(\mu)) [ (H^*(\Fl_G, \CL_\ell(w\star (\lambda+\mu) ))^\vee)_{\hat{\fg}^-}  ]_\sigma    \\
  &=\sum_{\mu\in P^\sigma}  \tr(\sigma| V(\mu)) (-1)^{\ell_\sigma(w)}[ (H^*(\Fl_G, \CL_\ell( \lambda+\mu) )^\vee)_{\hat{\fg}^-}  ]_\sigma    \\
  &=(-1)^{\ell_\sigma(w)} [ (H^*(\Fl_G, \CL_\ell( \bbC_\lambda\otimes V|_{B\rtimes \la\sigma \ra}  )  )^\vee)_{\hat{\fg}^-}  ]_\sigma,
 \end{align*}
where the second isomorphism follows from Lemma \ref{lem_change}. 
This finishes the proof.

\epf

\bp
\label{wall_case}
If $\lambda\in (P^+)^\sigma$ and $\lambda+\rho$ is in an affine wall of $\Wl$, then 
\[  [(H^*(\Gr_G, \CL_\ell( V_\lambda  \otimes  V))^\vee)_{\hat{\fg}^-}  ]_\sigma =0. \]
\ep
\bpf
By Part (3) of Proposition \ref{alcove_all_diagram}, $\lambda+\rho$ is in an affine wall of $\Wl^\sigma$, where by (\ref{affine_Weyl_inv}), $\Wl^\sigma\simeq W^\sigma\ltimes (\ell+\check{h})\iota(Q^\sigma)$. Hence in view of
 Lemma \ref{lattice_tran}, we can assume that
$\lambda+\rho$ is in the following affine wall of $\Wl^\sigma$ in $P^\sigma\otimes \bbR$,
\[  H_{\alpha_\sigma,a}=\{ \lambda+\rho \in P^\sigma\otimes\bbR \li  \la \lambda+\rho, \check{ \alpha}_\sigma  \ra=a  \}, \]
 where $\check{\alpha}_\sigma$ is the coroot of a root $\alpha_\sigma$ of $\fg_\sigma$, and 
 \[ a\in \begin{cases} 
 (\ell+\check{h}) \mathbb{Z}  \quad  \text{ if $\fg$ is not of type } A_{2n}  \\ 
   \frac{\ell+\check{h}}{2}\mathbb{Z}    \quad   \text{ if } \fg=A_{2n} 
 \end{cases}.\]
Equivialently, 
\be
\label{affine_wall}
 (\tau_{a \alpha_\sigma   } \cdot s_{\alpha_\sigma})\star (\lambda)=(s_{\alpha_\sigma}\cdot \tau_{-a \alpha_\sigma})\star (\lambda) =\lambda.
 \ee
where $s_{\alpha_\sigma}$ is the reflection with respect to $\alpha_\sigma$ in   $\Wl^\sigma$ and $\tau_{a\alpha_\sigma}$ is the translation by $a\alpha_\sigma$. Moreover, 
\[ (-1)^{\ell_\sigma(  \tau_{a \alpha_\sigma   }  \cdot s_{\alpha_\sigma}) }  = (-1)^{\ell_\sigma( \tau_{ a\alpha_\sigma  }  ) } (-1)^{\ell_{\sigma}(s_{\alpha_\sigma}) }= -1, \]
since by Lemma \ref{even}, $\ell_\sigma( \tau_{ a\alpha_\sigma  }  ) $ is an even integer.

By Proposition \ref{sign_change} we have
\[  [(H^*(\Gr_G, \CL_\ell( V_\lambda  \otimes  V))^\vee)_{\hat{\fg}^-}  ]_\sigma=- [(H^*(\Gr_G, \CL_\ell( V_\lambda  \otimes  V))^\vee)_{\hat{\fg}^-}  ]_\sigma. \]
Hence $  [(H^*(\Gr_G, \CL_\ell( V_\lambda  \otimes  V))^\vee)_{\hat{\fg}^-}  ]_\sigma =0$. 
\epf

\bt
\label{ring_hom}

The linear map $\pi_\sigma:  R(\fg, \sigma)\to R_\ell(\fg, \sigma)$ is a ring homomorphism. 
\et

\bpf

By Theorem \ref{new_old}, we can use the product $\otimes_\ell$ for $R_\ell(\fg,\sigma)$. 
We need to check that for any $\lambda,\mu\in (P^+)^\sigma$,
\be
\label{hom}
\pi_\sigma( [V_\lambda\otimes V_\mu]_\sigma  )=  \pi_\sigma([V_\lambda]_\sigma)\otimes_\ell \pi_\sigma([V_\mu]_\sigma). 
\ee

If $\lambda+\rho$ or $\mu+\rho$ is in an affine Wall, then by Proposition \ref{wall_case}, 
both sides of (\ref{hom}) are zero. Hence (\ref{hom}) holds. 

If $\lambda+\rho$ and $\mu+\rho$ are not in any affine Wall, let $\lambda_0\in P_\ell^\sigma$ such that $w_\lambda\star \lambda_0=\lambda$ and let $\mu_0\in P_\ell^\sigma$ such that $w_\mu\star \mu_0=\mu$ where $w_\lambda,w_\mu\in (\Wmin)^\sigma$, then
\begin{align*}  \pi_\sigma( [V_\lambda\otimes V_\mu]_\sigma  )&=[(H^*(\Gr_G, \CL_\ell(V_\lambda\otimes V_\mu))^\vee)_{\hat{\fg}^-} ]_\sigma \\
&=(-1)^{\ell_\sigma(w_\lambda)}[(H^*(\Gr_G, \CL_\ell(V_{\lambda_0}\otimes V_\mu))^\vee)_{\hat{\fg}^-} ]_\sigma\\
&=(-1)^{\ell_\sigma(w_\lambda) + \ell_\sigma(w_\mu)}[(H^*(\Gr_G, \CL_\ell(V_{\lambda_0}\otimes V_{\mu_0}))^\vee)_{\hat{\fg}^-} ]_\sigma\\
&=(-1)^{\ell_\sigma(w_\lambda) + \ell_\sigma(w_\mu)}[V_{\lambda_0}]_\sigma \otimes_\ell [V_{\mu_0}]_\sigma\\
&=\pi_\sigma([V_\lambda]_\sigma)\otimes_\ell \pi_\sigma([V_\mu]_\sigma),
\end{align*}
where the second, the third and the fifth equalities follows from Proposition \ref{sign_change}, and the fourth equality is the definition (\ref{new_def_prod}). 
This finishes the proof of the theorem.
\epf

We  can explicitly describe the map $\pi_\sigma$. 
\bco  
\label{descrip_map}
The map
$\pi_\sigma:   R(\fg,\sigma)\to  R_\ell(\fg, \sigma)$ can be described as follows, for any $\lambda\in (P^+)^\sigma$ we have
\[ \pi_\sigma([V_{ \lambda}]_\sigma)=
\begin{cases}
0  \quad  \text{ if }  \lambda+\rho \text{ belongs to an affine Wall of } \Wl  \text{ in } P_\bbR\\
(-1)^{\ell_\sigma(w)}[V_{w^{-1}\star \lambda}]_\sigma    \quad  \text{ if } w^{-1}\star \lambda \in P^\sigma_\ell \text{ for some } w\in (\Wmin)^\sigma.
\end{cases}
\]
\eco
\bpf
The corollary is an immediate consequence of Corollary \ref{sign_cor}, Proposition \ref{sign_change} and Proposition \ref{wall_case}. 
\epf

\subsection{Characters of the $\sigma$-twisted fusion ring}
\label{character_fusion_sect}

In Section \ref{character_fusion_sect} and Section \ref{proof_section} we basically follow the arguments in \cite[Section 9]{Be}. However our arguments of Lemma \ref{cardinality} and Proposition \ref{Casimir_cal} are substantially different, since in our setting there is no natural identification between $P_\sigma/(\ell+\check{h}) \iota(Q^\sigma)   $ and $T_{\sigma,\ell}$.

Recall that $P_\sigma$ (resp.\,$Q_\sigma$) is the weight lattice (resp.\,root lattice) of $\fg_\sigma$, and the bijection map $\iota: P_\sigma\simeq  P^\sigma$ defined in Section \ref{root_system}. 

Let $\bbZ [P_\sigma]$ be the group ring of $P_\sigma$; we denote by $(e^\lambda)_{\lambda\in P_\sigma}$ its  basis so that the multiplication in $\bbZ[P_\sigma]$ obeys the rule $e^\lambda e^\mu=e^{\lambda+\mu}$. The action of $W_\sigma$ and $W^\sigma_{\ell+\check{h}}\simeq W_\sigma\ltimes (\ell+\check{h})\iota(Q^\sigma)$ on $P_\sigma$ extends to $\bbZ[P_\sigma]$. We denote by $\bbZ[P_{\sigma}]_{W^\sigma}$ (resp.\,$\bbZ[P_{\sigma}]_{\Wl^\sigma}$) the quotient of $\bbZ[P_\sigma]$ by the sublattice spanned by $e^\lambda-(-1)^{\ell_\sigma(w)} e^{w\star \lambda}$ for any $w\in W_\sigma$  (resp.\,$w\in W^\sigma_{\ell+\check{h} }$). Let $p: \bbZ[P_{\sigma}]_{W_\sigma}\to\bbZ[P_{\sigma}]_{\Wl^\sigma} $ be the projection map. 

\bl
\label{kernel_describ}
The kernel $\ker(p)$ is spanned by the class of $e^{\lambda+\alpha}-e^{\lambda}$ in $\bbZ[P_{\sigma}]_{W_\sigma}$, for $\lambda\in P_\sigma$ and $\alpha\in (\ell+\check{h})\iota(Q^\sigma)$. 
\el
\bpf
We first define a group action $\bullet$ of $W^\sigma_{\ell+\check{h}}$ on $\bbZ[P_{\sigma}]$. For any $e^\lambda\in \bbZ[P_{\sigma}] $ and $w\tau_\alpha\in W^\sigma_{\ell+\check{h}}$ where $w\in W_\sigma$ and $\alpha\in (\ell+\check{h})\iota(Q^\sigma)$, define
\[ w\tau_\alpha\bullet e^\lambda:= (-1)^{\ell_\sigma(  w\tau_\alpha  ) }e^{w\star (\lambda+ \alpha) }. \]
It is easy to see that this gives a group action of $W^\sigma_{\ell+\check{h}}$ on $\bbZ[P_{\sigma}]$. Note that in the above formula, $(-1)^{\ell_\sigma(  w\tau_\alpha  ) }=(-1)^{\ell_\sigma(  w ) }$, since by Lemma \ref{even}, $\ell_\sigma(\tau_\alpha)$ is even.

Let $\mathbb{Z}[P_\sigma]_{(\ell+\check{h})\iota(Q^\sigma) }$ denote the space of coinvariants of $\mathbb{Z}[P_\sigma]$ with respect to the translation action of  $(\ell+\check{h})\iota(Q^\sigma)$. Consider the following short exact sequence
\[ \xymatrix{    0 \ar[r] &  K   \ar[r]  &  \mathbb{Z}[P_\sigma]  \ar[r]  &  \mathbb{Z}[P_\sigma]_{(\ell+\check{h})\iota(Q^\sigma) }  \ar[r] & 0   }, \]
where $K$ is the sublattice of $\mathbb{Z}[P_\sigma] $ spanned by $e^{\lambda+\alpha}-e^\lambda$ for $\lambda\in P_\sigma$ and $\alpha\in (\ell+\check{h})\iota(Q^\sigma) $. With respect to the action $\bullet$ of $W_\sigma$, we apply the functor of $W_\sigma$-coinvariants to the above short exact sequence. Since coinvariant functor is right exact, we get the following exact sequence
\[ \xymatrix{   K_{W_\sigma}   \ar[r]  &  \mathbb{Z}[P_\sigma]_{W_\sigma}  \ar[r]  &  (\mathbb{Z}[P_\sigma]_{(\ell+\check{h})\iota(Q^\sigma) } )_{W_\sigma} \ar[r] & 0   }. \]
Observe that
\[ \bbZ[P_{\sigma}]_{\Wl^\sigma} =  (\mathbb{Z}[P_\sigma]_{(\ell+\check{h})\iota(Q^\sigma) } )_{W_\sigma}.\]
This concludes the proof of the lemma.
\epf

By Proposition \ref{iso_twist} and Theorem \ref{ring_hom} we get a ring homomorphism  $\tilde{\pi}_\sigma: R(\fg_\sigma)\simeq R(\fg,\sigma)\to  R_\ell(\fg,\sigma)$. Let $\phi_\sigma$ be the map $R(\fg_\sigma)\to \bbZ[P_\sigma]_{W_\sigma} $ sending $[W_\lambda]$ to the class of $e^{\lambda}$. Similarly, let $\phi_{\sigma,\ell}$ be the map  $R_\ell(\fg,\sigma)\to  \bbZ[P_\sigma]_{W^\sigma_{\ell+\check{h}}}$ sending $[V_\lambda]_\sigma$ to the class $e^\lambda$ for any $\lambda\in P^\sigma_\ell$. By the same arguments as in \cite[Section 8]{Be}, $\phi_\sigma$ and $\phi_{\sigma,\ell}$ are bijections. As a consequence of Corollary \ref{descrip_map}, the following diagram commutes
\be
\label{comm_lem}
\xymatrix{
R(\fg_\sigma) \ar[r]^{\tilde{\pi}_\sigma}  \ar[d]^{\phi_\sigma} & R_\ell(\fg,\sigma)  \ar[d]^{\phi_{\sigma,\ell}} \\
\bbZ[P_\sigma]_{W_\sigma}   \ar[r]^{p}  & \bbZ[P_\sigma]_{W^\sigma_{\ell+\check{h}  }}  
}.
 \ee

For any $\lambda\in P_\sigma$, put  
\begin{equation} \label{def_J} J(e^{ \lambda+\rho})=\sum_{w\in W_\sigma} (-1)^{\ell_\sigma(w)} e^{w(\lambda+\rho_\sigma)}, \end{equation}
where $\rho_\sigma$ is the sum of all fundamental weights of $\fg_\sigma$. Recall that  $\iota(\rho)=\rho_\sigma$ via the bijection 
 $\iota: P^\sigma\simeq P_\sigma$. 
By Weyl character formula, for any $\lambda\in P^+_\sigma$ and $t\in T_\sigma$, we have
$  \tr(t| W_\lambda)  =\frac{J(e^{\lambda+\rho_\sigma})(t) }{J(e^{\rho_\sigma}  ) (t) } $. 
Let $T_{\sigma,\ell}$ be the finite subgroup of $T_\sigma$ given by 
\[ T_{\sigma, \ell}:=\{t\in T_\sigma \li   e^\alpha(t)=1, \alpha\in (\ell+\check{h}) \iota(Q^\sigma )\}. \]
\bp
\label{factor_through}
 For any $t\in T_{\sigma,\ell}$, the character $\tr(t|\cdot)$ factors through $\tilde{\pi}_\sigma: R(\fg_\sigma)\to R_\ell(\fg, \sigma)$.
\ep
\bpf
Let $j_t: \bbZ[P_\sigma]_{W^\sigma}\to \bbC$ be the additive map such that for any $\lambda\in P_\sigma$,
$ j_t(e^\lambda)=\frac{J(e^{\lambda+\rho_\sigma})(t)  }{ J(e^{\rho_\sigma} )(t) }$. By the definition of $\bbZ[P_\sigma]_{W^\sigma}$ and $J(\cdot)$, it is easy to check that $j_t$ is well-defined. 
By Weyl character formula, the following  diagram commutes:
\[ 
\xymatrix{
R(\fg_\sigma) \ar[r]^{\phi_\sigma}   \ar[rd]_{\tr(t|\cdot)}  &  \bbZ[P_\sigma]_{W_\sigma}  \ar[d]^{j_t}  \\
 &   \bbC
}.
\]
By the commutativity of the diagram (\ref{comm_lem}) and Lemma \ref{kernel_describ}, to show $\tr(t|\cdot)$ factors through $\tilde{\pi}_\sigma$, 
we need to check that $j_t$   takes the value zero on $e^{\lambda+\alpha}-e^{\lambda}$ for any $\lambda\in P_\sigma$ and $\alpha\in (\ell+\check{h})\iota(Q^\sigma) $. Since $t$ satisfies that $e^\alpha(t)=0$ for any $\alpha\in (\ell+\check{h}) \iota(Q^\sigma)$, it is clear that $j_t$ takes the value zero on  $e^{\lambda+\alpha}-e^{\lambda}$. This concludes the proof.
\epf
An element $t\in T_\sigma$ is  regular if the stabilizer of $W_\sigma$ at $t$ is trivial. We denote by $T^{\rm reg}_{\sigma, \ell}$ the set of regular elements in $T_{\sigma, \ell}$.
Let $\check{\rho}_\sigma$ denotes the sum of all fundamental coweights of $\fg_\sigma$. Consider the short exact sequence 
$$0\to 2\pi i \check{Q}_{\sigma}\to \mathfrak{t}_\sigma \to T_{\sigma}\to 1, $$
where $\check{Q}_{\sigma}$ denote the dual root lattice of $\fg_\sigma$ and $\ft_\sigma$ denotes the Cartan subalgebra of $\fg_\sigma$. 
 Let  $\check{L}_\sigma$ be the dual lattice of $\iota(Q^\sigma)$ in $\ft_\sigma$. We have the following natural isomorphism 
 \be
 \label{T_lattice}
 T_{\sigma,\ell}\simeq  (\frac{1}{\ell+\check{h}}\check{L}_\sigma) /\check{Q}_{\sigma}\simeq  \check{L}_\sigma/  (\ell+\check{h}) \check{Q}_{\sigma}.
 \ee
 For any $\check{\mu}\in \check{L}_\sigma$, we denote by $t_{\check{\mu}}$ the associated element of $\check{\mu}+\check{\rho}_\sigma$ in $T_{\sigma,\ell}$.

We put  $\check{P}_{\sigma, \ell}:= \{ \check{\mu}\in   \check{P}^+_\sigma    \li   \la \check{\mu}, \theta_\sigma  \ra_\sigma  \leq \ell  \}$,
 where $\theta_\sigma$ denotes the highest root of $\fg_\sigma$ and $\check{P}_\sigma^+$ denotes the set of dominant coweights of $\fg_\sigma$. 

\bl
\label{ident}
Assume that $\fg\neq A_{2n}$. There exists a bijection  $  \check{P}_{\sigma, \ell} \simeq T^{\rm reg}_{\sigma, \ell}/W_{\sigma} $
with the map given by $\check{\mu}\mapsto  t_{\check{\mu} }$,
\el
\bpf
When $\fg\neq A_{2n}$, by Lemma \ref{lattice_tran} $\iota(Q^\sigma)= Q_{\sigma}$. Thus $\check{L}_\sigma=\check{P}_\sigma$. 
We observe that $\la \check{\rho}_\sigma, \theta_\sigma \ra=\check{h}-1$ where $\check{h}$ is the dual Coxeter number of $\fg$. 
  This can be read from \cite[Table 2,p.66]{Hu2}).
It follows that  
\[ \check{P}_{\sigma, \ell}= \{ \check{\mu}\in   \check{P}^+_\sigma    \li   \la \check{\mu}+\check{\rho}_\sigma, \theta_\sigma  \ra_\sigma  < \ell+\check{h}  \}, \]
i.e. $ \check{P}_{\sigma, \ell}$ consists of all  points of $\check{P}^+_\sigma  $ sitting in the interior of the fundamental alcove with respect to the action of the affine Weyl group $W_\sigma\ltimes  (\ell+\check{h})\check{Q}_\sigma  $. From the isomorphism (\ref{T_lattice}), 
we can see that any $W_\sigma$-orbit in $T^{\rm reg}_{\sigma,\ell}$
has a unique representative in $ \check{P}_{\sigma, \ell}$. Hence the lemma follows. 
\epf

\begin{lemma}
\label{cardinality}
The cardinality of $T^{\rm reg}_{\sigma, \ell}/W_{\sigma}$ is the equal to the cardinality of $P_\ell^\sigma$.
 \end{lemma}
\bpf
When $\fg$ is of type $A_{2n}$, by Lemma \ref{lattice_tran}, $\iota(Q^\sigma)=\frac{1}{2}Q_{\sigma,\ell}$ where $Q_{\sigma,\ell}$  is the lattice spanned by long roots of $G_\sigma$. The proof of this lemma is exactly the same as the proof of \cite[Lemma 9.3]{Be}. We omit the detail.

Now we assume $\fg\neq  A_{2n}$. 
Put
\[   P_{\sigma,\ell}:=\{  \lambda\in P^+_\sigma \li  \la \lambda, \check{\theta}_\sigma \ra_\sigma   \leq  \ell   \},\]
where $ \check{\theta}_\sigma$ denotes the highest coroot of $\fg_\sigma$. In view of (\ref{compatibility}) and Lemma \ref{highest_root}, the map $\iota$ induces a natural bijection $\iota:    P^\sigma_\ell  \simeq P_{\sigma, \ell}$.

In view of Lemma \ref{ident}, 
we are reduced to show that $\check{P}_{\sigma, \ell}$ and $P_{\sigma,\ell}$ have the same cardinality. If $\fg_\sigma$ is not of type $B_n$ or $C_n$, it is  true, since in this case  weight lattice and  coweight lattice, root lattice and coroot lattice can be identified.
Otherwise, if $\fg_\sigma$ is of type $B_n$ or $C_n$, by comparing the highest roots of $B_n$ and $C_n$ (see \cite[Table 2,\,p.66]{Hu2}), we conclude that $\check{P}_{\sigma, \ell}$ and $P_{\sigma,\ell}$  indeed have the same cardinality.

\epf

The following proposition completely describes all characters of $R_{\ell}(\fg, \sigma)$.
\bp
\label{all_characters}
 $\{ \tr(t|\cdot)  \li   t\in  T^{\rm reg}_{\sigma,\ell} /W_\sigma \}$ gives a full set of  characters of $R_\ell(\fg,\sigma)$.
\ep
\bpf
This is an immediate consequence of Proposition \ref{factor_through} and Lemma \ref{cardinality}.
\epf
\subsection{Proof of Theorem \ref{main}}
\label{proof_section}
Let $\check{T}_{\sigma, \ell}$ denote the finite abelian subgroup
$ \check{T}_{\sigma,\ell}:=  P_\sigma/ (\ell+\check{h})\iota(Q^\sigma) $. 
For any $\lambda\in P_\sigma$, we denote by $\check{t}_\lambda$ the element in $\check{T}_{\sigma,\ell}$ associated to $\lambda+ \rho_\sigma$. 

Recall that $\Phi_\sigma$ is the set of roots of $\fg_\sigma$. 
In the following lemma we determine  $\chi(\omega_\sigma)$ for each $\chi=\tr(t|\cdot )$, where $\omega_\sigma$ is the Casimir element defined in (\ref{Casimir}).
\bp
\label{Casimir_cal}
For any $t\in T^{\rm reg}_{\sigma,\ell}$, we have 
$ \sum_{\lambda\in  P_\ell^\sigma}  | \tr(t|W_\lambda) |^2 =\frac{|T_{\sigma,\ell} |   }{\Delta_\sigma(t)}$,
where $\Delta_\sigma= \prod_{\alpha\in \Phi_\sigma}  (e^{\alpha}-1)$. 
\ep
\bpf
When $\fg=A_{2n}$, the proof of this lemma is identical to the proof of \cite[Lemma 9.7]{Be}. We omit the detail. 

Now we assume $\fg\neq  A_{2n}$. In this case, we have
\[ \check{T}_{\sigma,\ell}=P_\sigma/(\ell+\check{h})Q_\sigma,\quad  \text{ and } \quad T_{\sigma,\ell}\simeq \check{P}_\sigma/(\ell+\check{h})\check{Q}_\sigma.\] 

For any $ \lambda\in P_\sigma$ and  $ \check{\mu}\in \check{L}_\sigma=\check{P}_\sigma$, we have
\[ J(e^{\lambda+\rho_\sigma})(t_{\check{\mu}})=\sum_{w\in W_\sigma}  (-1)^{\ell_\sigma(w)  }  e^{ 2\pi i  \frac{   \la  \lambda+\rho_\sigma, w(\check{\mu}+\check{\rho}_\sigma)  \ra_\sigma  }  {\ell+\check{h}}   }   =J(e^{\check{\mu} + \check{\rho}_\sigma })(\check{t}_\lambda  ), \]
where we put $J(e^{\check{\mu} + \check{\rho}_\sigma })=\sum_{w\in W_\sigma} (-1)^{\ell_\sigma(w)}e^{w(\check{\mu}+ \check{\rho}_\sigma  )} $.
By Weyl character formula, we have
\[ \sum_{\lambda\in  P_\ell^\sigma} | \tr ( t_{\check{\mu}  } | W_\lambda )     |^2=  \frac{1}{\Delta_\sigma(t_{\check{\mu}})}  \sum_{\lambda\in P_\ell^\sigma }  | J(e^{\check{\mu}+\check{\rho}_\sigma  })(\check{t}_\lambda )  |^2. \]

We now introduce an inner product $(\cdot,\cdot )$ on the space $L^2(\check{T}_{\sigma,\ell}  )$ of functions on the finite abelian group $\check{T}_{\sigma,\ell}$, 
\[ (\phi, \psi) :=\frac{1}{|\check{T}_\ell|}\sum_{\check{t}\in \check{T}_{\sigma,\ell} }  \phi(\check{t})\overline{\psi(\check{t})}, \quad \text{for any functions $\phi,\psi$ on $\check{T}_{\sigma,\ell}$. }  \]    
The function $J(e^{\check{\mu}+\check{\rho}_\sigma  })$ on $\check{T}_{\sigma,\ell}$ is $W_\sigma$-antisymmetric, i.e. $ J(e^{ w\cdot (\check{\mu}+\check{\rho}_\sigma ) })=(-1)^{\ell_\sigma(w)}  J(e^{\check{\mu}+\check{\rho}_\sigma  })$. 
It shows that if $t$ is not regular, then for any $\check{t}\in \check{T}_{\sigma,\ell}$, $J(e^{\check{\mu}+\check{\rho}_\sigma  })(\check{t})=0 $. 
It follows that 
\[  \sum_{\lambda\in P_\ell^\sigma }  | J(e^{\check{\mu}+\check{\rho}_\sigma  })(\check{t}_\lambda )  |^2  = \frac{|\check{T}_{\sigma,\ell} |}{|W_\sigma|} || J(e^{\check{\mu}+\check{\rho}_\sigma  }) ||, \]
where $|| J(e^{\check{\mu}+\check{\rho}_\sigma  }) ||=\sqrt{( J(e^{\check{\mu}+\check{\rho}_\sigma  }), J(e^{\check{\mu}+\check{\rho}_\sigma  }) ) }$.

If $t$ is regular, in view of Lemma \ref{ident} we can assume $t=t_{\check{\mu}}$ where  $\check{\mu}\in \check{P}_{\sigma, \ell}$.
Now we show that the restriction of $e^{w\cdot (\check{\mu}+\check{\rho}_\sigma)}$ on $\check{T}_{\sigma,\ell}$ are all distinct. 
For any two distinct elements $w,w'\in W_\sigma$, if $e^{w(\check{\mu}+\check{\rho}_\sigma)}$ and $e^{w'(\check{\mu}+\check{\rho}_\sigma)}$ are equal on $\check{T}_{\sigma,\ell}$, it means that the pairing $\la w(\check{\mu}+\check{\rho}_\sigma)-w'(\check{\mu}+\check{\rho}_\sigma), \lambda  \ra_\sigma\in (\ell+\check{h} )\bbZ$ for any $\lambda\in P_\sigma$. Equivalently, $w(\check{\mu}+\check{\rho}_\sigma)-w'(\check{\mu}+\check{\rho}_\sigma)\in (\ell+\check{h})\check{Q}_\sigma$. It is impossible as $\check{\mu}+\check{\rho}_\sigma$ is in the fundamental alcove of the affine Weyl group $W_\sigma \ltimes  (\ell+\check{h})\check{Q}_\sigma$.

By the orthogonality relation for the characters of $\check{T}_{\sigma,\ell}$, we have  $ || J(e^{\check{\mu}+\check{\rho}_\sigma  }) ||=|W_\sigma|  $.
Hence, 
\[ \sum_{\lambda\in  P_\ell^\sigma}  | \tr(  t_{\check{\mu}  } |  W_\lambda )     |^2=\frac{|\check{T}_{\sigma,\ell}|}{\Delta_\sigma(t_{\check{\mu}})}. \]
 From the non-degeneracy of the pairing
$  \check{T}_{\sigma,\ell}\times  T_{\sigma,\ell}\to  \bbC^\times$ 
given by $(\check{t}_{\lambda}, t_{\check{\mu} })\mapsto  e^{ 2\pi i \frac{ \la \lambda+\rho_\sigma,\check{\mu}+\check{\rho}_\sigma \ra_\sigma}{\ell+\check{h}}  }  $,
   we have $|T_{\sigma,\ell}|=|\check{T}_{\sigma,\ell}|$. This concludes the proof of the proposition. 
\epf

Finally Theorem \ref{main}  follows from  Proposition \ref{machine}, Proposition \ref{Casimir_cal}, and Proposition \ref{all_characters}.

\subsection{A corollary of Theorem \ref{main}}
\label{stange_corollary_sect}
 Let $\sigma $ be a nontrivial diagram automorphism on $\fg=sl_{2n+1}$. Then the orbit Lie algebra $\fg_\sigma$ is isomorphic to $ sp_{2n}$. 
 
\begin{theorem}
\label{strange_formula}
With the same setting as in Theorem \ref{main}. 
If   $\ell$ is an odd positive integer, then we have the following formula
\[  \tr(\sigma| V_{sl_{2n+1}, \ell, \vec{\lambda}}(C,\vec{p}))=\dim V_{sp_{2n}, \frac{\ell-1}{2}, \vec{\lambda}}(C,\vec{p}). \]
\end{theorem}
\bpf
By assumption, $\la\lambda_i, \check{\theta}\ra  \leq \ell$ for any $\lambda_i$. In view of (\ref{compatibility}) and Lemma \ref{highest_root}, 
we have $\la \iota(\lambda_i), \check{\theta}_{\sigma,s}  \ra_\sigma \leq  \ell/2 $, where $\check{\theta}_{\sigma,s}$ is the coroot of the highest root $\theta_\sigma$ of $\fg_\sigma$. Since $\ell$ is odd and $\la \iota(\lambda), \check{\theta}_{\sigma,s}  \ra_\sigma $ is an integer, it follows that  $\la \iota (\lambda), \check{\theta}_{\sigma,s}  \ra_\sigma \leq  \frac{\ell-1}{2}$. 

Note that   $P_\sigma=\frac{1}{2}Q_{\sigma,\ell}$ where $Q_{\sigma,\ell}$ is the lattice spanned by long roots of $\fg_\sigma$. Moreover, $\check{h}=2n+1$ and   $\check{h}_{\sigma}=n+1$ where $\check{h}_{\sigma}$ is the dual Coxeter number of $\fg_\sigma$. Combining the  Verlinde formula (\ref{verlinde}) and Theorem \ref{main}, the  corollary follows. 

\epf

 \end{document}